\documentclass[11pt,a4]{article}
\usepackage[T2A]{fontenc}
\usepackage[cp866]{inputenc}
\usepackage[mathscr]{eucal}
\usepackage{amssymb}
\usepackage{amsmath}
\renewcommand{\leq}{\leqslant}

\newcommand{\Q}{{\mathbb Q}}
\renewcommand{\C}{{\mathbb C}}
\newcommand{\N}{{\mathbb N}}
\newcommand{\R}{{\mathbb R}}
\newcommand{\Z}{{\mathbb Z}}

\newcommand{\F}{{\frak F}}
\newcommand{\Ll}{{\frak L}}
\newcommand{\PP}{{\frak P}}
\renewcommand{\k}{\rule{0.7em}{0.7em}}
\usepackage{fancyhdr}
\begin{document}

\sloppy

{\normalsize 

\thispagestyle{empty}

\mbox{}
\\[-1.5ex]
\centerline{
{\bf TOPOLOGICAL, SMOOTH AND HOLOMORPHIC CLASSIFICATIONS}
}
\\[0.5ex]
\centerline{
{\bf OF NONAUTONOMOUS LINEAR DIFFERENTIAL SYSTEMS }
}
\\[0.5ex]
\mbox{}\hfill
{\bf AND PROJECTIVE MATRIX RICCATI EQUATIONS}
\hfill\mbox{}
\\[2.75ex]
\centerline{
\bf 
V.N. Gorbuzov$\!{}^{\ast}$, V.Yu. Tyshchenko$\!{}^{\ast\ast}$
}
\\[1.75ex]
\centerline{
\it 
Department of Mathematics and Computer Science, 
Yanka Kupala Grodno State University,
}
\\[1ex]
\centerline{
\it 
Ozeshko 22, Grodno, 230023, Belarus
}
\\[1ex]
\centerline{
$\!{}^{\ast}$E-mail: gorbuzov@grsu.by
}
\\[1ex]
\centerline{
$\!{}^{\ast\ast}$E-mail: valentinet@mail.ru
}
\\[5.25ex]
\centerline{{\large\bf Abstract}}
\\[1ex]
\indent
The questions of global topological, smooth and holomorphic classifications of the differential systems, defined by covering foliations, are considered. The received results are applied to non\-autonomous linear differential systems and projective matrix Riccati equations.
\\[2ex]
\\[1.25ex]
\\[5.5ex]
\centerline{{\large\bf Contents}}
%{\large\bf Contents} 
\\[1.5ex]
Introduction\footnote[1]{The main results of this paper has been published in the articles
{\it Buletinul AS Moldova.Matematica}, 1998, No. 3 (28), 49-56; 
{\it Vestnik of the Yanka Kupala Grodno State Univ.}, 2002, Ser. 2, No. 1, 14-19;
2006, Ser. 2, 
\linebreak
No. 1, 20-28; 
{\it Differential Equations}, 2003, Vol. 39, No. 4, 565-567;  2003, Vol. 39, No. 12, 1596-1599; 
{\it Vestnik of the Byelorussian State Univ.}, 2005, Ser. 1, No. 3, 74-79; 
2007, Ser. 1, No. 1, 82-86; 
2007, Ser. 1, No. 3, 
\linebreak
96-101; 
2010, Ser.\! 1, No.\! 1, 109-113; 
{\it Izvestiya of the Gomel Sate Univ.}, 2006, No. 5, 3-6; 
and in the mo\-no\-g\-ra\-p\-hy Gorbuzov V.N. {\it Integrals of differential systems} (Russian), Grodno State Univ., Grodno, 2006.

}
                                                                                                                   \dotfill\ 2
\\[0.75ex]
1. Covering foliations and their classifications 
                                                                                                                   \dotfill\ 2
\\[0.5ex]
2. Phase groups of covering foliations, defined by complex nonautonomous 
\\
\mbox{}\hspace{1em}
linear differential systems
                                                                                                                     \dotfill\ 4
\\[0.5ex]
3. Conjunctions of linear actions on $\C^n$
                                                                                                                     \dotfill\ 5
\\[0.5ex]
4. Applications to the complex nonautonomous linear differential systems 
                                                                                                                     \dotfill\ 9
\\[0.5ex]
5. Phase groups of covering foliations, defined by complex nonautonomous
\\
\mbox{}\hspace{1em}
projective matrix Riccati equations 
                                                                                                                     \dotfill\ 11
\\[0.5ex]
6. Conjunctions of linear-fractional actions on $\C P^n$
                                                                                                                     \dotfill\ 12
\\[0.5ex]
7. Applications to the complex nonautonomous projective matrix Riccati  equations
                                                                                                                     \dotfill\ 18
\\[0.5ex]
8. Phase groups of covering foliations, defined by real nonautonomous 
\\
\mbox{}\hspace{1em}
linear differential systems 
                                                                                                                     \dotfill\ 19
\\[0.5ex]
9. Conjunctions of linear actions on $\R^n $
                                                                                                                     \dotfill\ 20
\\[0.5ex]
10. Applications to real nonautonomous linear differential systems
                                                                                                                     \dotfill\ 26
\\[0.5ex]
11. Phase groups of covering foliations, defined by real nonautonomous Riccati equations 
                                                                                                                     \dotfill\ 27
\\[0.5ex]
12. Conjunctions of linear-fractional actions on $\overline{\R} $
                                                                                                                     \dotfill\ 27
\\[0.5ex]
13. Applications to real nonautonomous Riccati equations 
                                                                                                                     \dotfill\ 33
\\[0.75ex]
References 
                                                                                                                     \dotfill\ 34
\\[-2ex]

\newpage

\mbox{}
\\[-1.75ex]
\centerline{
\large\bf  
Introduction
}
\\[1.35ex]
\indent
The questions of global topological classification of differential systems (i.e. defined by them foliations) have been considered for the first time in the work [1]. In it the criterion of global topological equivalence for real autonomous linear ordinary differential systems of general situation, in particular, has been received. Further the given problem was considered in [2] and has definitively been solved in [3].
Similar problems for real completely solvable autonomous linear differential systems in case of two independent variables were studied in [4], and in a case, when number of independent variables on 1 less then numbers of dependent variables, were studied in [5] and [6].

In a complex case global topological classification of autonomous linear ordinary differential systems of general situation has been spent in [7 -- 11] and in [12] given problem has been considered in a case of completely solvable autonomous linear differential systems. Besides, in the works 
[13 -- 15] this question was studied for complex autonomous polynomial ordinary differential systems of the second order.

For nonautonomous differential systems the problem of global topological classification was considered only for the scalar complex linear ordinary differential equation [16; 17].

It is necessary to notice, that all criteria of global topological equivalence of corresponding differential systems received only for integrated in the quadratures cases (that has essentially facilita\-ted reception of these criteria).

In given article we will spend global topological [12; 18; 20 -- 26], smooth [21; 22; 24 -- 26], 
\linebreak
$\R$--holomorphic [19; 24; 25] (in a complex case) and holomorphic [24; 26] classifications of nonautonomous linear differential systems and projective matrix Riccati equations [27], generally speaking, not integrated in quadratures.
\\[2.5ex]
\centerline{
\large\bf  
1. Covering foliations and their classifications
}
\\[1.5ex]
\indent
{\bf Definition 1.1.} 
{\it 
Let $A$ and $B$ be path connected smooth varieties of dimensions ${\rm dim}\ A=n$ and ${\rm dim}\ B=m.$ 
Smooth foliation $\F$ of dimension $m$ on the variety $A\times B,$ locally transversal to $A\times\{b\}$ for all $b\in B,$ we will name 
\textit{\textbf{a covering foliation}}, if the projection $p\colon A\times B\to B$ on the second factor defines for each layer of it foliation covering of the variety $B.$ Thus variety $A$ we will name 
\textit{\textbf{a phase layer}}, and variety $B$ we will name \textit{\textbf{a base}} of covering foliation $\F$.
}
\vspace{0.5ex}

{\bf Definition 1.2.} 
{\it 
Let $\F_c^{}$ 
\vspace{0.25ex}
be a layer of the covering foliation $\F,$ containing the point $c\in A\times B.$ 
\textit{\textbf{The phase group}} $Ph(\F,b_0^{}),\ b_0^{}\in B,$ 
\vspace{0.35ex}
of the covering foliation $\F$ we will name the group of the diffeomorphisms 
\vspace{0.25ex}
${\rm Diff}(A,\pi_1^{}(B,b_0^{}))$ of the actions on the phase layer $A$ by fundamental group $\pi_1^{}(B,b_0^{})$ 
with noted point $b_0^{},$ defined under formulae 
\\[1.5ex]
\mbox{}\hfill
$
\Phi^{\gamma}(a)=q\circ r\circ s (1)
$ 
\ for all $a\in A,
$ 
\ for all 
$
\gamma\in \pi_1^{}(B,b_0^{}),
\hfill
$ 
\\[1.5ex]
where $r$ 
\vspace{0.25ex}
is a lifting of one of ways $s(\tau)\subset B$ for all $\tau\in [0,1],$ corresponding to the element $\gamma$ of the group 
\vspace{0.35ex}
$\pi_1^{}(B,b_0^{}),$ on the layer $\F_{(a,s(0))}^{}$ of the covering foliation $\F$ in the point $(a,s(0)),$ and $q\colon A\times B\to A$ is a projection to the first factor.
}
\vspace{0.5ex}

It is easy to see, that owing to path connectivity and smoothness of the variety $B,$ then phase groups $Ph (\F, b_1^{})$ and 
$Ph (\F, b_2^{}) $ are smoothly conjugated for any two points $b_1^{}$ and $b_2^{}$ of the base $B.$ 
Therefore further, as a rule, we will speak simply about of \textit{\textbf{the phase group}} $Ph (\F) $ of the covering foliation $\F,$ 
not connecting it with any point of the base $B.$

Further by consideration of the questions connected with topological (smooth,
$\!\!\R\!$-ho\-lo\-morphic, holomorphic) classifications of the covering foliations, we will believe everywhere, that their phase layers are homeomorphically (diffeomorphically, $\R\!$-holomorphically, holomorphically) equivalent each other, and in two last cases 
covering foliations we will consider, accordingly, holomorphic and $\R\!$-holomorphic. Thus under $\R\!$-holomorphism (holomorphism) 
we will understand bijective $\R\!$-holomor\-phic (holomorphic) map, having $\R$-holomor\-phic (holomorphic) to themselves inverse map.

{\bf Definition 1.3.} 
{\it 
We will say that the covering foliation $\F^1$ on the variety $A_1^{}\times B_1^{}$ 
is \textit{\textbf{topologically}} {\rm(}\textit{\textbf{smoothly}}, $\R\!$-\textit{\textbf{holomorphically}}, 
\textit{\textbf{holomorphically}}{\rm)} 
\textit{\textbf{equivalent}} to the covering foliation $\F^2$ on the variety $A_2^{}\times B_2^{}$ 
\vspace{0.25ex}
if exists the homeomorphism {\rm(}the diffeomorphism, the $\R\!$-holomorphism, the holo\-morphism{\rm)} 
$h\colon A_1^{}\times B_1^{}\to A_2^{}\times B_2^{}$ such that 
\\[1.5ex]
\mbox{}\hfill
$
q_2^{}\circ h (A_1^{}\times B_1^{}) =A_2^{},
\quad 
h\bigl (\F^1 _ {c_1^{}} \bigr) = \F^2 _ {h (c_1^{})}$ 
\ for all $c_1^{}\in A_1^{}\times B_1^{},
\hfill
$ 
\\[1.5ex]
where $q_2^{}$ is a projection to the first factor.
}
\vspace{0.75ex}

{\bf Definition 1.4.} 
{\it 
We will say that the covering foliation $\F^1$ on the variety $A_1^{}\times B_1^{}$ is 
\textit{\textbf{embeddable}} {\rm(}\textit{\textbf{smoothly embeddable}}, 
$\R\!$-\textit{\textbf{holomorphically embeddable}}, \textit{\textbf{holomorphically embeddable}}{\rm)} 
in the covering foliation $\F^2$ on the variety $A_2^{}\times B_2^{}$ if there is such embedding {\rm(}smooth embedding, 
$\R\!$-holomorphic embedding, holomorphic embedding{\rm)} 
\\[1.5ex]
\mbox{}\hfill
$
h\colon A_1^{}\times B_1^{} \hookrightarrow A_2^{}\times B_2^{},
\hfill
$ 
\\[1.5ex]
that $q_2^{}\circ h (A_1^{}\times B_1^{}) =A_2^{}$ and 
$h\bigl (\F^1 _ {c_1^{}} \bigr) \hookrightarrow \F^2 _ {h (c_1^{})}$ for all $c_1^{}\in A_1^{}\times B_1^{}.$
}
\vspace{0.75ex}

{\bf Definition 1.5.} 
\vspace{0.25ex}
{\it 
We will say that the covering foliation $\F^1$ on the variety $A_1^{}\times B_1^{}$ 
\textit{\textbf{covers}} {\rm(}\textit{\textbf{smoothly covers}}, $\R\!$-\textit{\textbf{holomorphically covers}}, 
\textit{\textbf{holomorphically covers}}{\rm)} the covering foliation $\F^2$ on the variety $A_2^{}\times B_2^{}$ 
\vspace{0.35ex}
if exists such covering {\rm(}smooth covering, $\R\!$-holomorphic covering, holomorphic covering{\rm)} 
\vspace{0.35ex}
$h\colon A_1^{}\times B_1^{}\to A_2^{}\times B_2^{},$ that $q_2^{}\circ h (A_1^{}\times B_1^{}) =A_2^{}$ and 
$h\bigl (\F^1 _ {c_1^{}}\bigr) \to \F^2 _ {h (c_1^{})}$ for all $c_1^{}\in A_1^{}\times B_1^{}.$
}
\vspace{0.75ex}

{\bf Theorem 1.1.} 
{\it 
For topological {\rm(}smooth, $\R\!$-holomorphic, holomorphic{\rm)} equivalence of the covering foliations $\F^1$ and $\F^2$ 
it is necessary and enough existence of the isomorphism $\mu$ of the fundamental groups $\pi_1^{}(B_1^{})$ and 
\vspace{0.35ex}
$\pi_1^{}(B_2^{}),$ generated by the homeomorphism {\rm(}diffeomorphism, $\R\!$-holo\-morphism, holomorphism{\rm)} 
\vspace{0.35ex}
$g _ {\mu}^{}\colon B_1^{}\to B_2^{}$ of the bases, and existence of the homeomorphism {\rm(}diffeomorphism, 
$\R\!$-holomorphism, holomorphism{\rm)} $f\colon A_1^{}\to A_2^{}$ of phase layers  such that
\\[1ex]
\mbox{}\hfill                                                  % (1.1)
$
f\circ \Phi_1 ^ {\gamma_1^{}} = \Phi_2 ^ {\mu (\gamma_1^{})} \circ f 
$
\ for all 
$
\gamma_1^{}\in\pi_1^{} (B_1^{}), 
$
\hfill {\rm(1.1)}
\\[1.65ex] 
where $\Phi _ {\xi} ^ {\gamma _ {\xi}^{}} \in Ph\bigl (\F ^ {\xi}\bigr),\  
\gamma _ {\xi}^{} \in\pi_1^{} (B _ {\xi}^{}),\  \xi =1,2.$
}
\vspace{0.9ex}

{\sl Proof.} 
At first we will notice, that for the smooth $(\R\!$-holomorphic, holomorphic) varieties, which are bases of the covering foliations, definitions of fundamental groups and their actions on smooth $(\R\!$-holomorphic, holomorphic) phase layers by means 
of continuous and by means of smooth $(\R\!$-holomorphic, holomorphic) ways are equivalent.
\vspace{0.25ex}

{\sl The necessity.} Let the map 
\vspace{0.25ex}
$h\colon A_1^{}\times B_1^{}\to A_2^{}\times B_2^{}$ defines topological (smooth, 
$\R\!$-ho\-lo\-mo\-r\-p\-hic, holomorphic) equivalence of the covering foliations $\F^1$ and $ \F^2.$
Map $h$ induces the homeomorphism (the diffeomorphism, the $\R\!$-holo\-morphism, the holo\-morphism) 
\vspace{0.35ex}
$g _ {\mu}^{} =p_2^{}\circ h $ of bases $B_1^{}$ and $B_2^{}$ which, in turn, induces the isomorphism 
\vspace{0.35ex}
$\mu\colon \pi_1^{} (B_1^{}) \to\pi_1^{} (B_2^{})$ of their fundamental groups. 
\vspace{0.25ex}
Let $b_1^0$ is noted point of the variety $B_1^{}.$ Then on the basis of definition 1.2 and that fact that map 
$h$ translates the layer $\F^1 _ {(a_1^{}, b_1^0)} $ in the layer $ \F^2 _ {h (a_1^{}, b_1^0)},$ we come to relations (1.1), 
where $f=q_2^{}\circ h.$
\vspace{0.25ex}

{\sl The sufficiency.} 
Let for operations of phase groups $Ph (\F^1) $ and $Ph (\F^2) $ relations (1.1) are fulfilled at the homeomorphism 
(the diffeomorphism, the $\R\!$-holo\-morphism, the holo\-morphism) 
\vspace{0.5ex}
$g_ {\mu}^{}$ of bases $B_1^{}$ and $B_2^{}.$ We take the way 
$s_1^{}\colon [0,1] \to B_1^{}$ such that $s_1^{}(0) =b_1^{},$ $s_1^{} (1) =b_1^0.$ 
Also we will suppose
\\[1.75ex]
\mbox{}\hfill                                                  % (1.2)
$
h (a_1^{}, b_1^{}) =( l_2 ^ {{}-1} \circ f\circ l_1^{} (a_1^{}),\, g _ {\mu}^{} (b_1^{}))
$
\ for all $a_1^{}\in A_1^{},
$
\ for all $b_1^{}\in B_1^{}, 
$
\hfill {\rm(1.2)}
\\[1.75ex] 
where $l_1^{}(a_1^{})=q_1^{}\circ r_1^{}(1),\ q_1^{}$ is a projection on the first  factor, 
\vspace{0.35ex}
$r_1^{}(\tau)$ is a raising of the way $s_1^{}$ on a layer of foliation $\F^1$ in the point 
\vspace{0.5ex}
$(a_1^{},b_1^{}),\ l_2^{{}-1}(a_2^{})=q_2^{}\circ r_2^{{}-1}(1),\ r_2^{{}-1}(\tau)$ is an outcome of a raising of the way 
\vspace{0.5ex}
$s_2^{{}-1}=g_{\mu}^{}\circ s_1^{{}-1}$ on a layer of foliation $\F^2$ in the point 
$(f\circ l_1^{}(a_1^{}),g_{\mu}^{}(b_1^0)),\ s_1^{{}-1}$ is a way, inverse to the way $s_1^{}.$ 

Now directly we come to a conclusion, that fiber bijective map (1.2) sets topological (smooth, $\R\!$-holomorphic,
holomorphic) equivalence of covering foliations $\F^1$ and $\F^2.$ \k
\vspace{0.5ex}

Similarly to given theorem it is proved following two assertions.
\vspace{0.5ex}

{\bf Theorem 1.2.} 
{\it 
For embedding {\rm(}smooth embedding, $\R\!$-holomorphic em\-bedding, holomorphic em\-bedding{\rm)} 
covering foliation $\F^1$ in covering foliation $\F^2$ it is necessary and enough existence of the homomorphism $\mu$ 
\vspace{0.25ex}
of fundamental group $\pi_1^{}(B_1^{})$ in fundamental group $\pi_1^{} (B_2^{}),$ generated by the embedding 
{\rm(}the smooth embedding, the $\R\!$-holomorphic embedding, the holomorphic embedding{\rm)} 
\vspace{0.25ex}
$g_{\mu}^{}\colon B_1^{} \hookrightarrow B_2^{}$ bases, and existence of the homeomorphism 
{\rm(}the diffeomorphism, the $\R\!$-holomorphism, the holomorphism{\rm)} 
$f\colon A_1^{}\to A_2^{}$ of phase layers such that relations {\rm (1.1)} are carried out.
}
\vspace{0.35ex}

{\bf Theorem 1.3.} 
{\it 
That covering foliation $\F^1$ covered {\rm(}smoothly covered, $\R\!$-holomorphic covered, holomorphic covered{\rm)} 
covering foliation $\F^2$ it is necessary and enough existence of the monomorphism $\mu$ fundamental group $\pi_1^{}(B_1^{})$ in \vspace{0.25ex}
fundamental group $\pi_1^{}(B_2^{}),$ generated by the covering {\rm(}the smooth covering, the $\R\!$-holomorphic covering, 
the holomorphic covering{\rm)} $g_{\mu}^{}\colon B_1^{} \hookrightarrow B_2^{}$ of the bases, 
\vspace{0.25ex}
and existence of the homeomorphism {\rm(}the diffeomorphism, the $\R\!$-holomorphism, the holomorphism{\rm)} 
$f\colon A_1^{}\to A_2^{}$ of phase layers such that relations {\rm (1.1)} are carried out.
}
\vspace{0.35ex}

Since theorems 1.1 -- 1.3 reduce problems of topological, smooth, $\R\!$-holomorphic, and holomorphic classifications 
of covering foliations to the problems of corresponding classifications of their phase groups at morphisms (isomorphisms,  homomorphisms,  monomorphism) we will consider further questions of topological, smooth, $\R\!$-holomorphic, and holomorphic 
conjunctions of phase groups of covering foliations, defined by the nonautonomous linear differential systems 
and projective matrix Riccati equations.
\\[3.75ex]
\centerline{
\large\bf  
2. Phase groups of covering foliations, 
}
\\[0.35ex]
\centerline{
\large\bf  
defined by complex nonautonomous  linear differential systems
}
\\[1.5ex]
\indent
We will consider linear differential systems
\\[1.5ex]
\mbox{}\hfill                                               % (2.1)
$
\displaystyle
dw =\sum\limits _ {j=1} ^mA_j^{} (z_1^{}, \ldots, z_m^{})\;\! w\, dz_j^{} 
$
\hfill (2.1)
\\[1ex]
and
\\[1.5ex]
\mbox{}\hfill                                               % (2.2)
$
\displaystyle
dw =\sum\limits _ {j=1} ^mB_j^{} (z_1^{}, \ldots, z_m^{})\;\! w\, dz_j^{}\;\!, 
$
\hfill (2.2)
\\[1.5ex]
ordinary at $m=1$ 
\vspace{0.5ex}
and completely solvable [28] at $m> 1,$ where $w = (w_1^{}, \ldots, w_n^{}),$ 
square matrices $A_j^{} (z_1^{}, \ldots, z_m^{}) = \| a _ {ikj}^{} (z_1^{}, \ldots, z_m^{})\|$ and 
\vspace{0.75ex}
$B_j^{} (z_1^{}, \ldots, z_m^{}) = \| b _ {ikj}^{} (z_1^{}, \ldots, z_m^{})\|$ of the size $n$ consist  
from holomorphic functions 
\vspace{0.5ex}
$a _ {ikj}^{}\colon A\to\C $ and $b _ {ikj}^{}\colon B\to\C,\ i,k =1,\ldots,  n,$ 
$j =1,\ldots, m,$ multiplication of matrixes we will carry out by multiplication of columns 
of the first matrix for corresponding lines of the second, path connected holomorphic varieties $A $ and $B $ 
are holomorphically equivalent each other, fundamental groups $\pi_1^{}(A)$ and $\pi_1^{}(B)$ have final number 
$ \nu\in\N $ of the forming.

The common solutions of linear differential systems (2.1) and (2.2) define covering foliations $\Ll^1$ and $ \Ll^2,$ accordingly, on varieties 
$\C^n\times A$ and $\C^n\times B.$ 

We will say that linear differential systems (2.1) and (2.2) are 
{\it topologically {\rm(}smoothly, $\R\!$-holomorphically, holomorphically{\rm)} equivalent} if exists the homeomorphism (the diffeomorphism, the $\R\!$-holomorphism, the holo\-morphism) $h\colon\C^n\times A\to\C^n\times B,$ translating layers of the covering foliation 
$\Ll^1$ in layers of the covering foliation $\Ll^2.$ 

Similarly we introduce the concepts of 
{\it embedding {\rm(}smooth embedding, $\R\!$-holomorphic embedding, holomorphic embedding}) and 
{\it covering {\rm(}smoothly covering, $\R\!$-holomorphically covering, holomorphically covering}) of linear differential systems.

The phase group $Ph (\Ll^1) $ of the covering foliation $ \Ll^1$ is generated by the forming nondegenerate linear transformations 
\vspace{0.25ex}
$P_r^{}w$ for all $w\in\C^n,\ P_r^{}\in GL (n, \C),\ r =1,\ldots, \nu,$ and the phase group $Ph (\Ll^2) $ of the covering foliation $\Ll^2$
\vspace{0.25ex}
is generated by the forming nondegenerate linear transformations $Q_r^{}w$ for all 
$w\in\C^n,\ Q_r^{}\in GL (n, \C),\ r =1,\ldots, \nu.$
\\[3.75ex]
\centerline{
\large\bf  
3. Conjunctions of linear actions on $\C^n$
}
\\[1.5ex]
\indent
We will consider a problem about a finding of necessary and sufficient conditions of existence such 
homeomorphism (diffeomorphism, $\R\!$-holo\-morphism, holomorphism) 
$f\colon \C^n\to\C^n,$ that identities 
\\[1ex]
\mbox{}\hfill
$ 
f (P_r^{}w) =Q_r^{}f (w)
$
\ for all $w\in\C^n,$ 
\ for all $r\in I, 
$
\hfill (3.1) 
\\[1.75ex]
take place, 
\vspace{0.25ex}
where $f (w)\! =\! (f_1^{}(w),\ldots, f_n^{} (w)),$ square matrices 
$P_r^{}\!\in GL (n, \C),\  \Q_r^{}\!\in GL (n, \C)$ for all $r\in I,\ I $ are some set of indexes. 
\vspace{0.25ex}
Group of linear actions on $\C^n,$ formed by the matrices $P_r^{}$ for all $r\in I,$ we will designate by $L^1,$ 
and through $L^2$ we will designate the similar group, formed by the matrices $Q_r^{}$ for all $r\in I $. 
Besides, further everywhere a set $\{\lambda_1^{},\ldots,\lambda_n^{}\}$ of nonzero complex numbers we will name 
{\it simple} if 
\\[1.5ex]
\mbox{}\hfill
$
\dfrac{\lambda_k^{}}{\lambda_l^{}}\not\in s_{lk}^{\pm 1},
\quad 
s_{lk}^{}\in\N,
\ \ 
l\neq k,
\ \ 
k=1,\ldots, n,
\ \
l=1,\ldots, n,
\hfill
$ 
\\[1.5ex]
and a square matrix of the size $n> 1$ we will name {\it simple} if it has simple structure and a simple collection of eigenvalues.

Consider at first a topological conjunction of the Abelian linear groups $L^1$ and $L^2$. In this case we will notice that if all matrices 
$P_r^{}$ (all matrices $Q_r^{})$ have simple structure, for all $r\in I,$ they are reduce to a diagonal kind by the common transformation of similarity.
\vspace{0.75ex}

{\bf Theorem 3.1} [16; 17]. 
{\it 
For the topological conjunction {\rm (3.1)} at $n=1$ of linear groups $L^1$ and $L^2$ it is necessary and enough, that either 
\\[1.5ex]
\mbox{}\hfill                                      % (3.2)
$
q_r^{}=p_r^{}\;\!|p_r^{}|^{\alpha},
\quad 
{\rm Re}\,\alpha >{}-1
$
\ for all $r\in I, 
$
\hfill {\rm (3.2)}
\\[0.5ex]
or}
\\[0.5ex]
\mbox{}\hfill                                      % (3.3)
$
q_r^{} =\overline {p} _r^{}\;\!|p_r^{}|^{\alpha},
\quad 
{\rm Re}\, \alpha>{}-1$
\ {\it for all} $r\in I. 
$
\hfill {\rm (3.3)}
\\[2ex]
\indent
{\sl Proof. The necessity.} 
We will assume, that conjugating homeomorphism $f$ keeps orientation (a case, when homeomorphism $f$ 
changes orientation, it is considered similarly).

As rotations of a complex plane around of an origin of coordinates on the angles $\varphi$ and $\psi,$ where 
${}-\pi <\varphi\leq\pi,\ {}-\pi <\psi\leq\pi, $ are topologically conjugated, if and only if $\varphi =\psi,$ that for all $r\in I $ with a 
condition $ |p_r^{} | = 1,$ we have $q_r^{}=p_r^{},$ and, so, the relations (3.2) are carried out at any $\alpha.$

If $|p _ {r_1^{}}^{}| \neq 1$ and $ |p _ {r_2^{}}^{}| \neq 1$ for some $r_1^{}$ and $r_2^{}$ from $I,$ then
\\[1.75ex]
\mbox{}\hfill                                      % (3.4)
$
{\rm sgn}\;\!\ln |p _ {r_1^{}}^{} | =
{\rm sgn}\;\!\ln |p _ {r_2^{}}^{}|. 
$
\hfill {\rm (3.4)}
\\[1.75ex]
It is easy to see, that there are such sequences $ \{l_s^{} \} $ and $ \{m_s^{}\} $ of integers, that
\\[1.75ex]
\mbox{}\hfill                                      % (3.5)
$
\lim\limits _ {s\to {}+\infty}\ p _ {{}_{\scriptstyle r_1^{}}}^ {\, {}^{\scriptstyle l_s^{}}}\, 
p _ {{}_{\scriptstyle r_2^{}}}^ {\,{}^{\scriptstyle m_s^{}}} =1, 
$
\hfill {\rm (3.5)}
\\[2ex]
and moreover $\lim\limits _ {s\to {}+\infty} |l_s^{} | =\lim\limits _ {s\to {}+\infty} |m_s^{} | = {}+ \infty $.
\vspace{0.75ex}

From identity (3.1) follows that 
\\[2ex]
\mbox{}\hfill
$
f\bigl(\;\!p _ {{}_{\scriptstyle r_1^{}}}^ {\,{}^{\scriptstyle l_s^{}}}\;\! 
p _ {{}_{\scriptstyle r_2^{}}}^ {\,{}^{\scriptstyle m_s^{}}}w\bigr) =\;\!
q _ {{}_{\scriptstyle r_1^{}}}^ {\,{}^{\scriptstyle l_s^{}}}\;\! 
q _ {{}_{\scriptstyle r_2^{}}}^ {{}^{\scriptstyle m_s^{}}} f (w)$ 
\ for all $w\in\C,$
\ for all $s\in\N.
\hfill
$ 
\\[2ex]
From here taking into account (3.5) it is had, that
\\[1.5ex]
\mbox{}\hfill                                      % (3.6)
$
\lim\limits _ {s\to {}+\infty}\, 
q _ {{}_{\scriptstyle r_1^{}}} ^ {{}^{\scriptstyle l_s^{}}}\;\! 
q _ {{}_{\scriptstyle r_2^{}}} ^ {{}^{\scriptstyle m_s^{}}} =1. 
$
\hfill {\rm (3.6)}
\\[1.75ex]
Therefore from (3.5) and (3.6) for some values of logarithms we receive equalities
\\[1.75ex]
\mbox{}\hfill                                      % (3.7)
$
\dfrac{\ln\;\! |p _ {r_1^{}}^{} |}{\ln\;\! |p _ {r_2^{}}^{} |} =
{}-\lim\limits _ {s\to {}+\infty}\ \dfrac{m_s^{}}{l_s^{}} 
$
\hfill (3.7)
\\[1ex]
and
\\[1.25ex]
\mbox{}\hfill                                      % (3.8)
$
l_s^{}\ln p _ {r_1^{}}^{} +
m_s^{}\ln p _ {r_2^{}}^{} =\;\!
l_s^{}\ln q _ {r_1^{}}^{} + m_s^{}\ln q _ {r_2^{}}^{}.
$
\hfill (3.8)
\\[2ex]
Let's divide the left and right parts of equality (3.8) on $l_s^{}$ and we will pass to a limit at $s\to {}+\infty $. 
Then taking into account (3.7) we will receive that 
\\[1.75ex]
\mbox{}\hfill                                      % (3.9)
$
\dfrac{\ln\;\! q _ {r_1^{}}^{}-\;\!\ln\;\! p _ {r_1^{}}^{}}{\ln\;\! |p _ {r_1^{}}^{}|} = 
\dfrac{\ln\;\! q _ {r_2^{}}^{}-\;\!\ln\;\! p _ {r_2^{}}^{}}{\ln\;\! |p _ {r_2^{}}^{} |}\,. 
$
\hfill (3.9)
\\[1.75ex]
\indent
Now, believing, that 
\\[1.25ex]
\mbox{}\hfill
$
\alpha _ {r_1^{}}^{} = 
\dfrac{\ln\;\! q _ {r_1^{}}^{}-\ln\;\! p _ {r_1^{}}^{}}{\ln\;\! |p _ {r_1^{}}^{} |},
\hfill
$ 
\\[1.5ex]
we receive equality (3.2), where $\alpha =\alpha _ {r_1^{}}^{}.$
\vspace{0.5ex}
Besides, from (3.9) follows that $\alpha _ {r_1^{}}^{} = \alpha _ {r_2^{}}^{} = \alpha$ for all $r_1^{}$ and $r_2^{}$ from $I$ taking into account that $ |p _ {r_1^{}}^{} | \neq 1$ and $|p _ {r_2^{}}^{}| \neq 1.$ 
\vspace{0.5ex}
And, at last, the inequality ${\rm Re}\, \alpha>{}-1$ follows from equalities 
${\rm Re}\, \alpha={\rm Re}\, \alpha _ {r_1^{}}^{} =\dfrac{\ln\;\! |q _ {r_1^{}}^{}|}{\ln\;\! |p _ {r_1^{}}^{}|}-1$ and (3.4).  
\vspace{0.5ex}

{\sl The sufficiency} is proved by construction of conjugating homeomorphism 
\vspace{0.35ex}
$f(w) = \gamma w|w |^ {\alpha}$ for all $w\in\C$ at performance of relations (3.2); 
and $f (w) = \gamma\;\! \overline{w}\;\! |w |^{\alpha}$ for all $w\in\C$ at per\-for\-man\-ce of relations (3.3). \k
\vspace{0.75ex}

{\bf Theorem 3.2.} 
{\it 
Let at $n> 1$ the matrices 
\\[1.5ex]
\mbox{}\hfill
$
P_r^{}=S\;\! {\rm diag}\{p_{1r}^{},\ldots, p_{nr}^{}\}\;\! S ^ {{}-1},
\quad 
Q_r^{}=T\;\! {\rm diag}\{q_{1r}^{}, \ldots, q_{nr}^{}\}\;\! T^{{}-1},
\hfill
$ 
\\[1.5ex]
and the matrixes $\ln P_r^{}$ and $\ln Q_r^{}$ be simple for all $r\in I.$
\vspace{0.25ex}
Then for the topological conjunction {\rm (3.1)} of linear groups $L^1$ and $L^2$ it is necessary and enough existence 
of such permutation $\varrho\colon (1,\ldots, n) \to (1, \ldots, n) $ and complex numbers $\alpha_k^{}$ 
with ${\rm Re}\, \alpha_k>{}-1,\ k =1,\ldots,  n,$ that either 
\\[0.5ex]
\mbox{}\hfill                                        % (3.10)
$
q_{\varrho(k)\;\! r}^{} =p_{kr}^{} |p_{kr}^{}| ^ {\alpha_k^{}}$ 
\ for all 
$
r\in I, 
\quad
k =1,\ldots,  n,
$
\hfill {\rm (3.10)}
\\[1ex]
or}
\\[0.5ex]
\mbox{}\hfill                                        % (3.11)
$
q_{\varrho(k)\;\! r}^{} = 
\overline{p}_{kr}^{}\;\! |p_{kr}^{} | ^ {\alpha_k^{}}$
\ for all 
$
r\in I, 
\quad
k =1,\ldots,  n.
$
\hfill {\rm (3.11)}
\\[2ex]
\indent
{\sl Proof.} With the help of replacement $\xi (w) =T^{{}-1} f (Sw)$ for all $w\in\C^n, $ from identities (3.1) we pass to identities 
\\[1.5ex]
\mbox{}\hfill                                        % (3.12)
$
\xi \bigl({\rm diag} \{p _ {1r}^{}, \ldots, p _ {nr}^{} \}\;\! w\bigr) =
{\rm diag}\{q _ {1r}^{}, \ldots, q _ {nr}^{}\}\, \xi (w)
$
\ for all 
$w\in\C^n,
$
\ \ for all 
$
r\in I.
$
\hfill (3.12)
\\[1.75ex]
Therefore the topological conjunction of linear groups $L^1$ and $L^2$ is equivalent to performance of identities (3.12).

{\sl The necessity.} 
Let identities (3.12)  are carried out. Holomorphism $u_r^{}(w) =P_r^{}w$ for all $w\in\C^n $ 
(holomorphism $v_r^{}(w) =Q_r^{}w$ for all $w\in\C^n)$ defines on space $\C^n $ the invariant holomorphic foliation 
${\frak U}_r^{}$ (the invariant holomorphic foliation ${\frak V}_r^{})$ 
\vspace{0.5ex}
of complex dimension $1,$ defined by the basis of nondegenerate absolute invariants [29] 
\vspace{0.75ex}
$w_k^{{}-\ln p_{nr}^{}} w_n^{\ln p_{kr}^{}},\ k =1,\ldots, n-1$ 
(by the basis of nondegenerate absolute invariants 
\vspace{0.5ex}
$w_k^{{}-\ln q_{nr}^{}} w_n^{\ln q_{kr}^{}},\, k\!=\! 1,\ldots, n\!-\!1)\!$ for all $\!r\!\in\! I.$ 

We will designate through $\C_k $ the coordinate complex plane $w_l^{}=0,\ l\neq k,\ l =1,\ldots, n,$ and through 
$\stackrel{{\small o}}{\C}_k$ we will designate the coordinate complex plane $\C_k^{}$ with the pricked out an origin of coordinates, 
$k =1,\ldots, n.$ 

As matrixes $\ln P_r^{}$ (matrixes $\ln Q_r^{})$ are simple for all $r\in I,$ that at: 

1) $\dfrac{\ln p_{kr}^{}}{\ln p_{nr}^{}}\not\in\R\ 
\biggl(\dfrac{\ln q_{kr}^{}}{\ln q_{nr}^{}}\not\in\R\biggr),$ 
\vspace{1.25ex}
the closure of each of hypersurfaces $w_k^{{}-\ln p_{nr}^{}}\;\!w_n^{\ln p_{kr}^{}}=C_k^{}$ and 
$w_k^{{}-\ln q_{nr}^{}}\;\! w_n^{\ln q_{kr}^{}}=C_k^{}$ at $C_k^{}\neq 0$ contained the hyperplane $w_k^{}=0;$ 
\vspace{1.25ex}

2) $\dfrac{\ln p_{kr}^{}}{\ln p_{nr}^{}}\not\in\Q\ 
\biggl(\dfrac{\ln q_{kr}^{}}{\ln q_{nr}^{}}\not\in\Q\biggr),$ 
\vspace{1.25ex}
closure of each of hypersurfaces $w_k^{{}-\ln p_{nr}^{}}\;\!w_n^{\ln p_{kr}^{}}=C_k^{}$ and 
$w_k^{{}-\ln q_{nr}^{}}\;\! w_n^{\ln q_{kr}^{}}=C_k^{}$ at $C_k^{}\neq 0$ 
\vspace{0.5ex}
contain the points, which are doing not belong to this hypersurface; 
\vspace{1.25ex}

3) $\dfrac{\ln p_{kr}^{}}{\ln p_{nr}^{}}\in\Q,\ 
\dfrac{\ln p_{kr}^{}}{\ln p_{nr}^{}}\neq s_{kr}^{{}-1},\ s_{kr}^{}\in\N\ 
\vspace{1.25ex}
\biggl(\dfrac{\ln q_{kr}^{}}{\ln q_{nr}^{}}\in\Q,\ 
\dfrac{\ln q_{kr}^{}}{\ln q_{nr}^{}}\neq s_{kr}^{{}-1},\ s_{kr}^{}\in\N\biggr),$ 
closure of each of hypersurfaces $w_k^{{}-\ln p_{nr}^{}}\;\!w_n^{\ln p_{kr}^{}}=C_k^{}$ and 
\vspace{0.5ex}
$w_k^{{} -\ln q_{nr}^{}}\;\!w_n^{\ln q_{kr}^{}}=C_k^{}$ at $C_k^{}\neq 0$ 
or contain many-sheeted (multivalent) covering of hyperplanes $w_k^{}=0,$ 
\vspace{0.25ex}
or does not contain the fixed point $O\in\C^n$ 
\vspace{0.35ex}
of the holomorphisms $u_r^{}(w)$ for all $w\in\C^n,$ and $v_r^{}(w)$ for all $w\in\C^n; \ k=1,\ldots, n-1,$ for all $r\in I$.

Therefore the layers ${\frak u}_k^{} ={\stackrel{{\small o}}{\C}}_k^{}$ 
\vspace{0.35ex}
(the layers ${\frak v}_k^{} ={\stackrel{{\small o}}{\C}}_k^{}),\ k=1,\ldots, n,$ of the foliation ${\frak U}_r^{}$ 
(of the foliation ${\frak V}_r^{}),$ are homeomorphic to each other and are not homeomorphic 
to any other layers of this foliation, for all $r\in I.$ 
\vspace{0.25ex}

So, conjugating homeomorphism 
\vspace{0.25ex}
$\xi\colon\C^n\to\C^n$ takes the layer ${\frak u}_k^{}$ of the foliation 
${\frak U}_r^{}$ to the layer ${\frak v}_{\rho (k)}^{}$ 
\vspace{0.25ex}
of the foliation ${\frak V}_r^{},\ k =1,\ldots, n,$ for all $r\in I.$ 
\vspace{0.25ex}
Besides, from identities (3.12) follows that an origin of coordinates of space $\C^n$ 
\vspace{0.25ex}
is a fixed point of the homeomorphism $\xi.$ 

Then at it projections of the homeomorphism 
\vspace{0.35ex}
$\xi_k^{}\colon\C^n\to\C$ are that, that their restrictions 
$\tilde{\xi}_k^{}\colon\C_k^{}\to\C_{\rho(k)}^{}$ are homeomorphisms and 
\\[1.5ex]
\mbox{}\hfill
$
\tilde{\xi}_{\rho(k)}^{}\bigl(p_{kr}^{}\tilde{w}_k^{}\bigr)\! = \!
q_{\rho(k)\;\! r}^{}\;\! \tilde{\xi}_{\rho (k)}^{}\bigl(\tilde{w}_k^{}\bigr)$ 
for all $\tilde{w}_k^{}\! =\! (0, \ldots, 0, w_k^{}, 0, \ldots, 0),\, k\! =\!1,\ldots, n,$ for all $r\in I.
\hfill
$ 
\\[2ex]
From here on the basis of Theorem 3.1 we come to conclusion, that there are such complex numbers 
\vspace{0.5ex}
$\!\alpha_k^{}\!\!$ with $\!{\rm Re}\;\!\alpha_k^{}\!\!>\!-1,\!$ that one of relations (3.10) or (3.11), $\!\!k\!=\!1,\ldots, n,\!$ is carried out.

{\sl The sufficiency} 
\vspace{0.35ex}
is proved by the construction of conjugating homeomorphism 
$\xi\colon\C^n \to\C^n $ such that its projection 
\vspace{0.75ex}
$\xi_{\rho (k)}^{}(w) = \gamma_k^{}\;\!w_k^{}\;\!|w_k^{}|^{\alpha_k^{}}$ 
for all $w\in\C^n$ if relations (3.10) take place; and its projection 
\vspace{0.35ex}
$\xi_{\rho (k)}^{} (w) = \gamma_k^{}\;\!\overline{w}_k^{}\;\!|w_k^{}|^{\alpha_k^{}}$ for all $w\in\C^n$ if relations (3.11) take place; 
$k =1,\ldots,n.$ \k
\vspace{0.5ex}

Consider now a topological conjunction of the non-Abelian linear groups $L^1$ and $L^2.$
\vspace{0.75ex}

{\bf Theorem 3.3.} 
{\it 
From a topological conjunction of the non-Abelian linear groups $L^1$ and $L^2$ of general situation follows them 
$\R\!$-holomorphic conjunction.
}
\vspace{0.35ex}

{\sl Proof} of the theorem 3.3 directly follows from following two auxiliary statements (the lemmas 3.1 and 3.2).
\vspace{0.5ex}

{\bf Lemma 3.1.} 
\vspace{0.25ex}
{\it 
Let at $n=1$ and $I = \{1,2 \}$ linear groups $L^1$ and $L^2$ are topological conjugated, and the group of subgroup 
\vspace{0.25ex}
$\C^{\ast}$ of nonzero complex numbers on multiplication, generated by numbers $p _ {11}^{}$ and $p _ {12}^{},$ 
is dense in the set $\C$ of complex numbers.
\vspace{0.35ex}
Then conjugating homeomorphism $f $ is set either by the formula 
\vspace{0.5ex}
$f (w) = \gamma\;\! w\;\!|w |^{\alpha}$ for all $w\in\C,$ 
or by the formula $f (w) = \gamma\;\!\overline{w}\;\! |w |^{\alpha}$ for all $w\in\C,$ where} ${\rm Re}\, \alpha>{}-1.$
\vspace{0.5ex}

{\sl Proof.} 
Let conjugating homeomorphism $f $ keeps orientation (the case when it changes orientation, it is considered similarly). 
Owing to the theorem 3.1 relations 
\\[1.5ex]
\mbox{}\hfill
$
q_{1r}^{} =p_{1r}^{}\;\! |p_{1r}^{}|^{\alpha},
\quad 
r =1,2,
\quad 
{\rm Re}\, \alpha>{}-1,
\hfill
$ 
\\[1.75ex]
take place. On the basis of identities (3.1) we come to conclusion about justice of equalities 
\\[2ex]
\mbox{}\hfill
$
f\bigl(p_{{}_{\scriptstyle 11}}^{\,l} p_{{}_{\scriptstyle 12}}^{\,m}\bigr) =
f (1)\;\! p_{{}_{\scriptstyle 11}}^{\,l} p_{{}_{\scriptstyle 12}}^{\,m}\;\! 
|p_{{}_{\scriptstyle 11}}^{\,l} p_{{}_{\scriptstyle 12}}^{\,m} |^{\alpha}
$ 
\ 
for all $l\in\Z,$ 
\ for all $m\in\Z. 
\hfill
$ 
\\[2ex]
\indent
Owing to density of our subgroup of group 
\vspace{0.35ex}
$\C^{\ast}$ in the set $\C$ of complex numbers it is had, 
that for any complex number $w\in\C$ there are such sequences 
\vspace{0.75ex}
$\{l_s^{}(w) \} $ and $ \{m_s^{}(w)\} $ of integers, that 
$ \lim\limits_{s\to {}+\infty}\, p_{{}_{\scriptstyle 11}}^{\, l_s^{} (w)} p_{{}_{\scriptstyle 12}}^{\, m_s^{} (w)} =w,\ 
\lim\limits_{s\to {}+\infty} |l_s^{} (w) | = \lim\limits_{s\to {}+\infty} |m_s^{} (w) | = {}+ \infty.$ 
\vspace{0.5ex}
From here we receive the first representation of the statement of a lemma, where $ \gamma=f (1).$ \k
\vspace{0.75ex}

{\bf Lemma 3.2.} 
{\it 
Let the matrices 
\\[1.25ex]
\mbox{}\hfill
$
P_r^{}=S_r^{}\;\! {\rm diag} \{p_{1r}^{}, \ldots, p_{nr}^{}\}\;\! S_r^{{}-1},
\quad 
Q_r^{}=T_r^{}\;\! {\rm diag}\{q_{1r}^{},\ldots, q_{nr}^{}\}\;\! T_r^{{}-1},
\hfill
$ 
\\[1.75ex]
and the matrixes $\ln P_r^{}$ and $\ln Q_r^{}$ are simple for all $r\in I.$
Then from a topological conjunction of non-Abelian linear groups $L^1$ and $L^2$ of general situation follows them 
$\R\!$-linear conjunction {\rm(}i.e. homeomorphism $f\colon \C^n\to \C^n$ in identities {\rm (3.1)} is non\-de\-ge\-nerate 
linear $\R\!$-holomorphic transformation{\rm)}.
}
\vspace{0.35ex}

{\sl  Proof.} 
Let identities (3.1) take a place. By means of replacement $\xi(w) = T_1^{{}-1} f (S_1^{}w)$ for all $w\in\C^n $ 
from them we pass to identities
\\[1.75ex]
\mbox{}\hfill                                                             % (3.13)
$
\xi\bigl({\rm diag}\{p_{11}^{}, \ldots, p_{n1}^{}\}\;\! w\bigr) =
{\rm diag}\{q_{11}^{}, \ldots, q_{n1}^{}\}\;\! \xi (w)$ 
\ for all $w\in\C^n.
$
\hfill (3.13)
\\[1.75ex]
\indent
Now similarly, as well as at the proof of the theorem 3.1, we come to conclusion, 
that an origin of coordinates of space $\C^n$ is a fixed point of the homeomorphism $\xi,$ and its projections 
$\xi_k^{}\colon \C^n\to\C$ are that, that relations 
\\[1.5ex]
\mbox{}\hfill
$
\xi_{\rho (k)}^{} (0, \ldots, 0, p_{k1}^{} w_k^{}, 0, \ldots, 0) =
q_{\rho(k)\;\! 1}^{}\, \xi_{\rho (k)}^{} (0, \ldots, 0, w_k^{}, 0, \ldots, 0),
\hfill
$ 
\\[2.25ex]
\mbox{}\hfill
$
\xi_l^{} (0, \ldots, 0, w_k^{}, 0, \ldots, 0) =0
$
\ for all $w_k^{}\in\C,
\quad 
l\neq \rho (k),\ \ 
l =1,\ldots,  n,\ \ 
k =1,\ldots, n,
\hfill
$ 
\\[1.5ex]
take place; where $ \rho\colon (1, \ldots, n) \to (1, \ldots, n)$  is some permutation.
\vspace{0.25ex}

From here in case of general situation on the basis of relations (3.1), the theorem 3.2 and lemma 3.1 
we come to conclusion about justice of identities 
\\[1.5ex]
\mbox{}\hfill
$
\xi_{\rho (k)}^{} (0, \ldots, 0, w_k^{}, 0, \ldots, 0) = 
\gamma_k^{}w^{\ast}_k \;\!|w_k^{} |^{\alpha_k^{}}$ 
\ for all $w_k^{}\in\C,
\quad 
{\rm Re}\,\alpha_k^{}>{}-1,
\ \ k =1,\ldots, n,
\hfill
$ 
\\[2ex]
where $w^{\ast}_k=w_k^{}\vee\overline{w}_k^{},\ k =1,\ldots, n. $ 
From these identities further we come to relations 
\\[2ex]
\mbox{}\hfill                                                % (3.14)
$
\begin{array}{c}
\xi_{\rho (k)}^{} (w) = \gamma_k^{}w^{\ast}_k\;\! |w_k^{} |^{\alpha_k^{}} 
\bigl(1 + \varphi_{\rho (k)}^{} (w)\bigr)\ \ \text{for all}\ w\in\C^n,
\quad {\rm Re}\,\alpha_k^{}>{}-1,
\\[2.5ex]
\varphi_{\rho (k)}^{} (0, \ldots, 0, w_k^{}, 0, \ldots, 0) =0
\ \ \text{for all}\  w_k^{}\in\C\backslash \{0\},
\ \ 
k =1,\ldots,  n.
\end {array}
$
\hfill (3.14)
\\[2ex]
\indent
Now on the basis of identities (3.13) it is had, that 
$\varphi_k^{}\bigl({\rm diag}\{p_{11}^{}, \ldots, p_{n1}^{}\}\;\! w\bigr) = 
\varphi_k^{}(w)$ for all $w\in\C^n,\ k =1,\ldots, n. 
$ 
Considering relations (3.14) and basis of nondegenerate absolute invariants of foliation ${\frak U}_1^{}$ from the proof 
of the theorem 3.2, we receive representations 
\\[2ex]
\mbox{}\hfill                                                % (3.15)
$
\xi_{\rho (k)}^{} (w) = \gamma_k^{}w^{\ast}_k\;\! |w_k^{} |^{\alpha_k^{}} 
\Bigl(1 +\psi_{\rho (k)}^{}\Bigr(w_1^{}w_k^{{}-\ln p_{11}^{}/\ln p_{k1}^{}}, \ldots, 
w_{k-1}^{} w_k^{{}-\ln p_{k-1,1}^{}/\ln p_{k1}^{}}, 
\hfill
$
\\[0.5ex]
\mbox{}\hfill (3.15)
\\[0.5ex]
\mbox{}\hfill
$
w_{k+1}^{} w_k^{{}-\ln p_{k+1,1}^{}/\ln p_{k1}^{}}, \ldots, w_n^{}w_k^{{}-\ln p_{n1}^{}/\ln p_{k1}^{}}\Bigr)\Bigr)
$
\ for all $w\in\C^n,
\quad 
k =1,\ldots, n, 
\hfill
$ 
\\[2ex]
where functions $\psi_{\rho (k)}$ are continuous on the arguments, $k =1,\ldots, n.$ 

Taking into consideration non\-com\-mu\-tativity of linear groups and identity (3.13), on the basis of relations (3.14) and (3.15) 
taking into account simplicity of matrixes  from a lemma condition we come to its statement. \k
\vspace{0.5ex}

Consider cases of smooth, $\R\!$-holomorphic and holomorphic conjunctions of linear groups $L^1$ and $L^2.$
\vspace{0.25ex}

{\bf Theorem 3.4.} 
{\it 
Linear groups $L^1$ and $L^2$ are smooth $(\R\!$-holomorphic{\rm)} conjugated if and only if they are $\R\!$-linearly conjugated.
}
\vspace{0.25ex}

{\sl Proof. The necessity.} 
Let linear groups $L^1$ and $L^2$ are smooth $(\R\!$-holomorphic) conjugated. Then identities (3.1) 
take place at the diffeomorphism (at the $\R\!$-holomorphism) $f.$ 
Calculating in them full differential in the point $w=0,$ we have 
\\[1.5ex]
\mbox{}\hfill
$
{\sf D}_w^{}f (0) P_r^{}\;\!dw+{\sf D}_{\overline{w}}^{} f (0)\;\! \overline{P}_r^{}\;\! d\;\!\overline{w} \, = \,
Q_r^{}\;\!{\sf D}_w^{} f (0)\, dw+Q_r^{}\;\! {\sf D}_{\overline{w}}^{} f (0)\, d\;\!\overline{w}
$ 
\ 
for all $r\in I. 
\hfill
$ 
\\[1.75ex]
\indent
Therefore $\R\!$-linear transformation 
\vspace{0.25ex}
${\rm D}_w^{} f (0) w+{\sf D}_{\overline{w}}^{}\, f (0)\overline{w}$
for all $w\in\C^n$ defines the conjunction (3.1). As map $f $ 
\vspace{0.25ex}
is a diffeomorphism (an $\R\!$-holomorphism), it is nondegenerate.

{\sl The sufficiency} is proved by direct calculations. \k
\vspace{0.35ex}

Similarly we prove the following statement.
\vspace{0.35ex}

{\bf Theorem 3.5.} 
{\it 
Linear groups $L^1$ and $L^2$ are holomorphic conjugated if and only if they are linearly conjugated.
}
\\[2.5ex]
\centerline{
\large\bf  
4. Applications to the complex nonautonomous linear differential systems
}
\\[1.5ex]
\indent
Theorems 3.1, 3.2, 3.4, 3.5 and Lemma 3.2 allow on the basis of Theorems 1.1 -- 1.3 to spend topological, smooth, 
$\R\!$-holomorphic and holomorphic classifications of complex nonautonomous linear differential systems of a kind (2.1).
Besides, from the theorem 3.3 it is had such statement.
\vspace{0.35ex}

{\bf Theorem 4.1.} 
{\it 
From topological equivalence of complex nonautonomous linear differential systems 
with non-Abelian phase groups of general situation follows them $\R\!$-holomorphic equivalence.
}
\vspace{0.35ex}

From the theorem 4.1 follows that in the case of two and more dependent variables complex nonautonomous linear differential systems with coefficients, holomorphic on path connected holomorphic varieties with non-Abelian fundamental groups, are structurally unstable.

Let's result now concrete examples.

Consider the following linear equations: ordinary differential
\\[1.75ex]                           % (4.1)
\mbox{}\hfill
$
\displaystyle
\dfrac{dw}{dz} \ =\  \sum\limits_{r=1}^{\nu}\, \sum\limits_{l=1}^{\eta_{{}_{r1}}^{}} 
\, \dfrac{\Lambda_{rl1}^{}}{\left(z - z_{r1}^{}\right)^l}\, w
$
\hfill (4.1)
\\[1ex]
and 
\\[1.75ex]                           % (4.2)
\mbox{}\hfill
$
\displaystyle
\dfrac{dw}{dz} \ = \ 
\sum\limits_{r=1}^{\nu} \, 
\sum\limits_{l=1}^{\eta_{{}_{r2}}^{}} \,
\dfrac{\Lambda_{rl2}^{}}{\left(z - z_{r2}^{}\right)^l}\, w\,;
$
\hfill (4.2)
\\[1.75ex]
integral 
\\[1.75ex]                           % (4.3)
\mbox{}\hfill
$
\displaystyle
w(z) = 
\int\limits_{0}^{1}\prod\limits_{r=1}^{\nu}
\biggl(\dfrac{z -z_{r1}^{}}{\zeta - z_{r1}^{}}\biggr)^{\!\Lambda_{r11}^{}}
\exp\Biggl(\
\sum_{l=2}^{\eta_{{}_{r1}}}\frac{\Lambda_{rl1}}{1 - l}
\biggr(
\dfrac{1}{\left(z - z_{r1}\right)^{l-1}} -
\dfrac{1}{\left(\zeta - z_{r1}\right)^{l-1}}\biggr)\!\Biggr)
w(\zeta)\; d\zeta\vspace{-4mm}
$
\hfill (4.3)
\\[4ex]
and 
\\[3ex]                           % (4.4)
\mbox{}\hfill
$
\displaystyle
w(z) = 
\int\limits_{0}^{1}\prod_{r=1}^{\nu}
\biggl(\dfrac{z - z_{r2}^{}}{\zeta - z_{r2}^{}}\biggr)^{\!\Lambda_{r12}^{}}
\exp\Biggl(\
\sum_{l=2}^{\eta_{{}_{r2}}}
\dfrac{\Lambda_{rl2}}{1 - l}\biggl(
\dfrac{1}{\left(z - z_{r2}^{}\right)^{l-1}} -
\dfrac{1}{\left(\zeta - z_{r2}^{}\right)^{l-1}}\biggr)\!\Biggr)\;\!
w(\zeta)\; d\zeta\,;
$
\hfill (4.4)
\\[2.5ex]
integral-differential
\\[2.5ex]
\mbox{}\quad                                        %  (4.5)
$
\displaystyle
\dfrac{dw}{dz} \, = \ 
\int\limits_{0}^{1}\;\!
\sum_{r=1}^{\nu}\ \sum_{l=1}^{\eta_{{}_{r1}}}\
\dfrac{\Lambda_{rl1}^{}}{\left(z - z_{r1}^{}\right)^l} \
\prod_{r=1}^{\nu}
\Biggl(\biggl(\dfrac{z - z_{r1}^{}}{\zeta -z_{r1}^{}}\biggr)^{\!\Lambda_{r11}^{}}\cdot
\hfill
$
\\[-0.25ex]
\mbox{}\hfill  (4.5)
\\[0.75ex]
\mbox{}\hfill
$
\displaystyle
\cdot\,
\exp\Biggl(\
\sum_{l=2}^{\eta_{{}_{r1}}}
\dfrac{\Lambda_{rl1}^{}}{1 - l}\,
\biggl(\,\dfrac{1}{\left(z - z_{r1}^{}\right)^{l-1}} \,-\,
\dfrac{1}{\left(\zeta - z_{r1}^{}\right)^{l-1}}\biggr)\Biggr)\Biggr)
w(\zeta)\; d\zeta
\qquad\quad\mbox{}
$
\\[2.5ex]
and
\\[2.5ex]
\mbox{}\quad                                        %  (4.6)
$
\displaystyle
\dfrac{dw}{dz} \, = \ 
\int\limits_{0}^{1}\;\! 
\sum_{r=1}^{\nu}\
\sum_{l=1}^{\eta_{{}_{r2}}}\
\dfrac{\Lambda_{rl2}^{}}{\left(z - z_{r2}^{}\right)^l}\
\prod_{r=1}^{\nu}\Biggl(\biggl(\dfrac{z - z_{r2}^{}}{\zeta -z_{r2}^{}}\biggr)^{\!\Lambda_{r12}^{}}
\cdot
\hfill
$
\\[-0.25ex]
\mbox{}\hfill  (4.6)
\\[0.75ex]
\mbox{}\hfill
$
\displaystyle
\cdot\,
\exp\Biggl(\
\sum_{l=2}^{\eta_{{}_{r2}}}\, 
\dfrac{\Lambda_{rl2}^{}}{1 - l}\, 
\biggl(\dfrac{1}{\left(z - z_{r2}^{}\right)^{l-1}} \, -\, 
\dfrac{1}{\left(\zeta - z_{r2}^{}\,\right)^{l-1}}\biggr)\Biggr)\Biggr)
w(\zeta)\; d\zeta\,,
\qquad\quad\mbox{}
$
\\[2.5ex]
where 
\vspace{0.25ex}
$ \Lambda_{rl1}^{}$ and $ \Lambda_{rl2}^{}$ are complex numbers, 
the path from integrals (4.3) --- (4.6), connecting points $0$ and $1,$ are homeomorphic to segments and do not pass through points 
$z_{1r}^{} $ and $z_{2r}^{}, \ r =1,\ldots,  \nu.$ 

The common solutions of the linear equations (4.1), (4.3), and (4.5) define the same covering foliation
$\Ll^3$ on the variety $\C\times\Gamma_1^{\nu},$ 
and the common solutions of the linear equations (4.2), (4.4), and (4.6) 
\vspace{0.25ex}
define the same covering foliation $\Ll^4$ on the variety $ \C\times\Gamma_2^{\nu},$ 
where $\Gamma_1^{\nu}$ (accordingly, $\Gamma_2^{\nu})$ 
\vspace{0.35ex}
is an open complex plane with $\nu$ eliminated points 
$z_{r1}^{},\ r =1,\ldots, \nu$ (with $\nu$ eliminated points $z_{r2}^{},\ r =1,\ldots,\nu).$ 
\vspace{0.35ex}
The phase group $Ph (\Ll^3) $ of covering foliation $ \Ll^3$ is generated by nondegenerate linear transformations 
$p_r^{}w$ for all $w\in\C,\ r =1,\ldots, \nu,$ and the phase group $Ph (\Ll^4) $ of covering foliation $ \Ll^4$ is generated by nondegenerate linear \vspace{0.35ex}
transformations $q_r^{}w$ for all $w\in\C,\  r =1,\ldots, \nu,$ where 
$p_r^{} =\exp\bigl(2\pi\;\!i\;\!\Lambda_{r11}^{}\bigr),\
q_r^{} = \exp\bigl(2\pi\;\! i\;\!\Lambda_{r12}^{}\bigr), \ r=1,\ldots, \nu, \ i^2={}-1.$
\vspace{0.75ex}

Owing to Theorems 3.1 and 1.1 it is received the following statement.
\vspace{0.75ex}

{\bf Theorem 4.2.} 
{\it 
For topological equivalence of covering foliations $\Ll^3$ and $ \Ll^4$ it is necessary and enough existence of such permutation 
$ \varkappa\colon (1\ldots\nu) \to (1\ldots\nu)$ and complex number
\vspace{0.35ex}
$\alpha$ with ${\rm Re}\, \alpha\ne {}- 1,$ that either 
$q_{\varkappa (r)}^{} = p_r^{}\;\!|p_r^{}|^{\alpha},\ r =1,\ldots,\nu,$ or 
\vspace{1ex}
$q_{\varkappa (k)}^{} = \overline{p}_r^{}\;\!  |p_r^{}|^{\alpha},\ r =1,\ldots, \nu.$
}

Notice that the linear ordinary differential equations (4.1) and (4.2) in special case 
$ \eta_{r1}^{} =$ 
\linebreak
$= \eta_{r2}^{} =1, \ z_{r1}^{} =z_{r2}^{}, \ r =1,\ldots, \nu,$ are considered in [17]. 
Besides, set example shows possibility of application of the device of covering foliations for the mathematical objects 
which are distinct from differential systems.
\vspace{0.5ex}

Consider linear ordinary differential systems 
\\[2ex]
\mbox{}\hfill      % (4.7)
$
\dfrac{dw}{dz} =A_1^{}(z)\;\! w
$
\hfill (4.7)
\\[1ex]
and
\\[1.75ex]
\mbox{}\hfill      % (4.8)
$
\dfrac{dw}{dz} =B_1^{}(z)\;\! w,
$
\hfill (4.8)
\\[2.5ex]
where square matrices 
\vspace{0.75ex}
$A_1^{}(z) = \| a_{ik1}^{}(z)\|$ and $B_1^{}(z) = \| b_{ik1}^{}(z)\|$ of the size $n,\ n>1,$ consist from 
1-periodic holomorphic functions $a_{ik1}^{}\colon \C\to\C$ and 
$b_{ik1}^{}\colon \C\to\C,\ i, k =1,\ldots, n.$ 
\vspace{0.5ex}

The common solutions 
\vspace{0.35ex}
of linear ordinary differential systems (4.7) and (4.8) define covering foliations 
$\Ll^5$ and $ \Ll^6,$ accordingly, on the variety $\C^n\times Z,$ 
\vspace{0.35ex}
where $Z $ is the cylinder $S^1\times\R,$ $S^1$ is an unit circle.
\vspace{0.35ex}
The phase group $Ph (\Ll^5) $ of covering foliation $\Ll^5$ is generated by nondegenerate linear transformation 
\vspace{0.35ex}
$P_1^{}w$ for all $w\in\C^n,\ P_1^{}\in GL (n, \C),$ and the phase group $Ph (\Ll^6) $ 
\vspace{0.35ex}
covering foliation $\Ll^6$ is generated by nondegenerate linear transformation 
$Q_1^{}w$ for all $w\in\C^n,\ Q_1^{}\in GL (n, \C).$
\vspace{0.75ex}

From here on the basis of Theorems 3.2 and 1.1 it is had such statement.
\vspace{0.75ex}

{\bf Theorem 4.3.} 
{\it 
Let at $n> 1$ the matrices 
\\[1.5ex]
\mbox{}\hfill
$
P_1^{}=S\;\! {\rm diag}\{p_{11}^{},\ldots, p_{n1}^{}\}\;\! S^{{}-1},
\quad 
Q_1^{}=T\;\! {\rm diag}\{q_{11}^{},\ldots, q_{n1}^{}\}\;\! T^{{}-1},
\hfill
$ 
\\[1.75ex]
the matrixes $\ln P_1^{}$ and $\ln Q_1^{}$ are simple.
\vspace{0.5ex}
Then for topological equivalence of linear ordinary differential systems
\vspace{0.5ex}
{\rm (4.7)} and {\rm (4.8)} it is necessary and enough existence of such permutation 
$\varrho\colon (1,\ldots, n) \to (1, \ldots, n) $ and complex numbers $ \alpha_k^{}$ with 
\vspace{0.75ex}
${\rm Re}\, \alpha_k^{}\ne {}-1,\ k =1,\ldots, n,$ that either
$
q_{\varrho (k)\;\! 1}^{} =p_{k1}^{}\;\! |p_{k1}^{} |^{\alpha_k^{}},$
or
$
q_{\varrho (k)\;\! 1}^{} = \overline{p}_{k1}^{}\;\! |p_{k1}^{} |^{\alpha_k^{}},\ k =1,\ldots, n.$
}
\vspace{1ex}

Now we will consider 
\vspace{0.25ex}
linear ordinary differential systems (4.7) and (4.8) in a case when square matrices 
\vspace{0.5ex}
$A_1^{}(z)$ and $B_1^{}(z)$ of the size $n$ consist from holomorphic functions 
$a_{ik1}^{}\colon \Gamma_1^{\nu} \to\C $ and 
$b_{ik1}^{}\colon \Gamma_2^{\nu} \to\C,\ i =1,\ldots, n,\ k =1,\ldots, n.$ 
\vspace{0.5ex}
In this case the common solutions of linear ordinary differential systems (4.7) and (4.8) 
\vspace{0.35ex}
define covering foliations 
$\Ll^7$ and $ \Ll^8,$ accordingly, on the varieties $\C^n\times\Gamma_1^{\nu} $ and $ \C^n\times\Gamma_2^{\nu}.$
\vspace{0.35ex}
The phase group $Ph (\Ll^7) $ of covering foliation $ \Ll^7$ 
\vspace{0.35ex}
is generated by nondegenerate linear transformations 
$P_r^{}w$ for all $w\in\C^n,\ P_r^{}\in GL (n, \C),\ r =1,\ldots, \nu,$ 
\vspace{0.35ex}
and the phase group $Ph (\Ll^8) $ 
\vspace{0.35ex}
of covering foliation $ \Ll^8$ is generated by nondegenerate linear transformations 
$Q_r^{}w$ for all $w\in\C^n,\ P_r^{}\in GL (n, \C),\ r =1,\ldots, \nu. $
\vspace{0.5ex}

In a commutative case on the basis of Theorems 3.2 and 1.1 we receive the statement.
\vspace{0.5ex}

{\bf Theorem 4.4.} 
{\it 
Let at $n> 1$ the matrices 
\\[1.5ex]
\mbox{}\hfill
$
P_r^{}=S\;\! {\rm diag}\{p_{1r}^{}, \ldots, p_{nr}^{}\}\;\! S^{{}-1} 
\quad
\bigl(Q_r^{}=T\;\! {\rm diag}\{q_{1r}^{}, \ldots, q_{nr}^{} \}\;\! T^{{}-1}\bigr),
\hfill
$ 
\\[1.75ex]
the matrixes $\ln P_r^{}$ and $\ln Q_r^{}$ are simple, $r =1,\ldots,  \nu.$
\vspace{0.35ex}
Then for topological equivalence of linear ordinary differential systems
\vspace{0.25ex}
{\rm (4.7)} and {\rm (4.8)} it is necessary and enough existence of such permutations 
\vspace{0.35ex}
$\varkappa\colon (1,\ldots,\nu) \to (1,\ldots, \nu),\  
\varrho\colon (1, \ldots, n)\to (1, \ldots, n)$ and complex numbers 
$\alpha_k^{}$ with ${\rm Re}\, \alpha_k^{}\ne {}-1,\ k =1,\ldots, n,$ that either
\\[1.5ex]
\mbox{}\hfill
$
q_{\varrho (k)\;\! \varkappa (r)}^{} =p_{kr}^{}\;\! |p_{kr}^{} |^{\alpha_k^{}},
\quad 
\ k =1,\ldots, n, 
\ \ 
r =1,\ldots, \nu, 
\hfill
$
\\[0.75ex]
or
\\[0.75ex]
\mbox{}\hfill
$
q_{\varrho (k)\;\! \varkappa (r)}^{} = \overline{p}_{kr}^{}\;\! |p_{kr}^{} |^{\alpha_k^{}},
\quad 
k =1,\ldots, n,
\ \ 
r =1,\ldots, \nu.
\hfill
$
}
\\[2.25ex]
\indent
In a noncommutative case on the basis of Lemmas 3.1, 3.2 and Theorem 1.1 we have concrete constructive criteria of topological equivalence of  differential systems (4.7) and (4.8).
\\[3.75ex]
\centerline{
\large\bf  
5. Phase groups of covering foliations, defined by complex nonautonomous
}
\\[0.35ex]
\centerline{
\large\bf  
projective matrix Riccati equations 
}
\\[1.5ex]
\indent
We will consider homogeneous projective matrix Riccati equations [30]
\\[1.75ex]
\mbox{}\hfill                                  % (5.1)
$
\displaystyle
dv =\sum\limits_{j=1}^m 
{\frak A}_j^{} (z_1^{}, \ldots, z_m^{})\;\! v\, dz_j^{} 
$
\hfill (5.1)
\\[0.5ex]
and
\\[0.75ex]
\mbox{}\hfill                                  % (5.2)
$
\displaystyle
dv =
\sum\limits_{j=1}^m {\frak B}_j^{} (z_1^{}, \ldots, z_m^{})\;\! v\, dz_j^{}\;\!, 
$
\hfill (5.2)
\\[2ex]
ordinary at $m=1$ and completely solvable at $m> 1$, where 
\vspace{0.35ex}
$v=(v_1^{}, \ldots, v _ {n+1}^{}) $ are homogeneous coordinates, 
square matrices 
\vspace{0.5ex}
${\frak A}_j^{} (z_1^{}, \ldots, z_m^{}) = \| {\frak a}_{ikj}^{} (z_1^{}, \ldots, z_m^{})\|$ 
and ${\frak B}_j^{}(z_1^{}, \ldots, z_m^{}) = \|{\frak b}_ {ikj}^{} (z_1^{}, \ldots, z_m^{})\|$ of the size $n+1$ 
\vspace{0.75ex}
consist from holomorphic functions 
${\frak a}_{ikj}^{}\colon A\to\C$ and ${\frak b}_{ikj}^{}\colon B\to\C,\ 
i =1,\ldots, n+1,\ k =1,\ldots, n+1,\ j =1,\ldots, m,$ 
\vspace{0.5ex}
path connected holomorphic varieties $A$ and $B$ are holomorphically equivalent 
\vspace{0.35ex}
each other, fundamental groups $\pi_1^{}(A)$ and $\pi_1^{}(B)$ have final number 
\vspace{0.35ex}
$\nu\in\N $ of the forming.

The common solutions of homogeneous projective matrix 
\vspace{0.25ex}
Riccati equations (5.1) and (5.2) define covering foliations 
\vspace{0.25ex}
$\PP^1$ and $\PP^2,$ accordingly, on the varieties $\C P^n\times A $ and $\C P^n\times B $. 
We will say, that homogeneous projective matrix Riccati equations (5.1) and (5.2) are {\it topologically {\rm(}smoothly, 
$\R\!$-holomorphically, holomorphically{\rm)} equivalent}, 
\vspace{0.25ex}
if exists the homeomorphism, 
(the diffeomorphism, the $\R\!$-holomorphism, the holo\-morphism) 
$h\colon\C P^n\times A\to\C P^n\times B,$ translating the layers of the covering foliation $\PP^1$ 
in the layers of the covering foliation $\PP^2.$ 
Similarly we introduce the concepts of 
{\it embedding {\rm(}smooth embedding, $\R\!$-holomorphic embedding, holomorphic embedding}) and 
{\it covering {\rm(}smoothly covering, $\R\!$-holomorphically covering, holomorphi\-cally covering}) 
of homogeneous pro\-jective matrix Riccati equations.
\vspace{0.25ex}
The phase group $Ph (\PP^1) $ of the covering foliation $\PP^1$ is generated [30] 
\vspace{0.35ex}
by the forming nondegenerate linear-fractional transformations 
$P_r^{}v$ for all $v\in\C P^n,\ P_r^{}\in GL (n+1, \C),\ r =1,\ldots, \nu,$ 
\vspace{0.25ex}
and the phase group $Ph (\PP^2)$ of the covering foliation $\PP^2$ 
\vspace{0.25ex}
is generated by the forming nondegenerate linear-fractional transformations 
$Q_r^{}v$ for all $v\in\C P^n,\ Q_r^{}\in GL (n+1, \C),\ r =1,\ldots, \nu.$
\\[3.25ex]
\centerline{
\large\bf  
6. Conjunctions of linear-fractional actions on $\C P^n$
}
\\[1.5ex]
\indent
Now we will consider a problem about a finding of necessary and sufficient conditions of existence 
such homeomorphism (diffeomorphism, $\R\!$-holo\-morphism, holomorphism) 
$f\colon\C P^n\to\C P^n,$ that identities 
\\[1.5ex]
\mbox{}\hfill                                    % (6.1)
$
f (P_r^{}v) =Q_r^{}f (v)
$ 
\ for all 
$
v\in\C P^n,
$
\ for all 
$
r\in I, 
$
\hfill (6.1) 
\\[1.75ex] 
take place, where $f (v) = (f_1^{}(v), \ldots, f_{n+1}^{}(v)),$ 
square matrices 
\\[2ex]
\mbox{}\hfill
$
P_r^{}\in GL (n+1,\C),
\quad  
\Q_r^{}\in GL (n+1,\C)
$ 
\ for all $r\in I.
\hfill
$ 
\\[1.75ex]
Thus group of linear-fractional actions on $\C P^n,$ formed by the matrices 
\vspace{0.35ex}
$P_r^{}$ for all $r\in I$ (by the matrices $Q_r^{}$ for all $r\in I),$ 
we will designate through $PL^1$ (through $PL^2).$
\vspace{0.25ex}

Consider at first a topological conjunction of the Abelian linear-fractional groups $PL^1$ and $PL^2.$

{\bf Lemma 6.1.} 
{\it 
Let at $n=1$ linear-fractional groups $PL^1$ and $PL^2$ be topological conjugated.
Then normal Jordan forms of the matrices $P_r^{}$ and $Q_r^{},$ 
\vspace{0.25ex}
defining by nonidentical is linear-fractional transformations, have identical number of blocks of Jordan}, for all $r\in I $.
\vspace{0.35ex}

{\sl Proof} of the given statement is spent on the basis of that fact, that the quantity of fixed points 
of linear-fractional transformations coincides with the number of eigenvectors of the matrices, defining by these transformations. \k
\vspace{0.5ex}

On the basis of  Lemma 6.1 by direct calculations we come to such statement.
\vspace{0.35ex}

{\bf Lemma 6.2.} 
{\it 
For a topological conjugation at $n=1$ of Abelian linear-fractional groups $PL^1$ and $PL^2$ it is necessary, 
that normal Jordan forms of all matrices $P_r^{}$ and $Q_r^{},$ 
defining by nonidentical linear-fractional transformations, for all $r\in I,$ 
had identical number of blocks of Jordan.
}
\vspace{0.75ex}

{\bf Theorem 6.1.} 
{\it Let at $n=1$ the matrices 
\\[1.5ex]
\mbox{}\hfill
$
P_r^{}=S\;\! {\rm diag}\{p_{1r}^{}, p_{2r}^{}\}\;\! S^{{}-1}$ for all $r\in I,
\quad
Q_r^{}=T\;\! {\rm diag}\{q_{1r}^{}, q_{2r}^{}\}\;\! T^{{}-1}$  for all $r\in I.
\hfill
$
\\[1.75ex]
Then for a topological conjunction of linear-fractional groups $PL^1$ and $PL^2$ it is necessary and enough, that either 
\\[2ex]
\mbox{}\hfill                          % (6.2)
$
\dfrac{q_{1r}^{}}{q_{2r}^{}} = \dfrac{p_{1r}^{}}{p_{2r}^{}}\, 
\Bigl|\dfrac{p_{1r}^{}}{p_{2r}^{}}\Bigr|^{\alpha},
\quad 
{\rm Re}\, \alpha\ne {}-1,
$
\ for all $
r\in I, 
$
\hfill {\rm (6.2)} 
\\[1.25ex]
or}
\\[1.25ex]
\mbox{}\hfill                          % (6.3)
$
\dfrac{q_{1r}^{}}{q_{2r}^{}} = 
\dfrac{\overline{p}_{1r}^{}}{\overline{p}_{2r}^{}}\, 
\Bigl|\dfrac{p_{1r}^{}}{p_{2r}^{}}\Bigr|^{\alpha},
\quad 
{\rm Re}\, \alpha\ne {}-1,
$
\ for all 
$
r\in I.
$
\hfill (6.3)
\\[2.25ex]
\indent
{\sl Proof.} 
With the help of replacement $\xi (v) =T^{{}-1} f (Sv)$ for all $v\in\C P^1,$ 
from identities (6.1) at $n=1$ we pass to identities
\\[1.75ex]
\mbox{}\hfill                          % (6.4)
$
\xi \bigl({\rm diag}\{p_{1r}^{}, p_{2r}^{}\}\;\! v\bigr) =
{\rm diag}\{q_{1r}^{}, q_{2r}^{}\}\;\! \xi (v)
$
\ 
for all 
$
v\in\C P^1,
$
\ for all 
$
r\in I.
$
\hfill (6.4)
\\[2ex]
\indent
Therefore the topological conjunction of linear-fractional groups $PL^1$ and $PL^2$ 
is equivalent to performance of identities (6.4).
\vspace{0.25ex}

{\sl The necessity.} Let identities (6.4) are carried out.
\vspace{0.5ex}

If all $\dfrac{p_{1r}^{}}{p_{2r}^{}} =1$ for all $r\in I$ from (6.4) it is had, as 
\vspace{0.75ex}
$\dfrac{q_{1r}^{}}{q_{2r}^{}} =1$ for all $r\in I.$ 
Therefore in this case relations (6.2) are carried out at any $\alpha $ with ${\rm Re}\,\alpha\ne {}-1.$
\vspace{0.5ex}

Let now $\dfrac{p_{1r}^{}}{p_{2r}^{}} \ne 1,\ r\in I.$ 
\vspace{0.75ex}
On the basis of identities (6.4) it is received, that either 
$\xi (O_1^{}) =O_1^{},$ or $\xi (O_1^{}) =O_2^{},$ 
\vspace{0.5ex}
where $O_{\tau}^{}$ (the origings of coordinates of the affine cards 
$M_{\tau}^{} = \{v,\ v_{\tau}^{}\ne 0\} $ of the atlas $M$ of variety $\C P^1),\ \tau =1,2,$ 
\vspace{0.35ex}
are the common fixed points of linear-fractional transformations 
${\rm diag}\{p_{1r}^{}, p_{2r}^{}\}\;\! v$  for all $v\in\C P^1,$ 
\vspace{0.35ex}
and ${\rm diag}\{q_{1r}^{}, q_{2r}^{}\}\;\! v$ for all $v\in\C P^1$ for all $r\in I.$
\vspace{0.25ex}

If $\xi (O_1^{}) =O_1^{},$ 
\vspace{0.25ex}
then on the basis of Theorem 3.1 we do a conclusion, that at ${\rm Re}\, \alpha>{}-1$ 
take place either equalities (6.2), or equalities (6.3).
\vspace{0.25ex}

If $\xi (O_1^{}) =O_2^{},$ 
\vspace{0.35ex}
then with the help of replacement $\zeta (v) = (\xi_2^{}(v), \xi_1^{}(v))$ for all $v\in\C P^1,$ 
we come to the previous case and as a result it is received either equalities (6.2) at ${\rm Re}\, \alpha <{}-1,$ 
or equalities (6.3) at ${\rm Re}\, \alpha <{}-1.$
\vspace{0.35ex}

{\sl The sufficiency} is proved by construction of conjugating homeomorphism
\\[1.75ex]
\mbox{}\hfill
$
\xi (v) = \bigl(\gamma_1^{}v_1^{}|v_1^{} |^{\alpha},\;\! \gamma_2^{}v_2^{}|v_2^{} |^{\alpha}\bigr)
$ 
\ for all $v\in\C P^1,
\hfill
$ 
\\[1.75ex]
in case of performance of relations (6.2), and by construction of conjugating homeomorphism
\\[1.75ex]
\mbox{}\hfill
$
\xi (v) = \bigl(\gamma_1^{}\overline{v}_1^{}|v_1^{} |^{\alpha},\;\! \gamma_2^{}\overline{v}_2^{}|v_2^{} |^{\alpha}\bigr)
$ 
\ for all $v\in\C P^1,
\hfill
$ 
\\[1.75ex]
in case of performance of relations (6.3). \k
\vspace{1ex}

{\bf Theorem 6.2.} 
{\it 
Let at $n> 1$ the matrices 
\\[1.5ex]
\mbox{}\hfill
$
P_r^{}=S\;\! {\rm diag}\{p_{1r}^{}, \ldots, p_{n+1, r}^{}\}\;\! S^{{}-1},
\quad 
Q_r^{}=T\;\! {\rm diag}\{q_{1r}^{},\ldots, q_{n+1,r}^{}\}\;\! T^{{}-1},
\hfill
$ 
\\[2ex]
sets of numbers
$\biggl\{\ln\dfrac{p_{1r}^{}}{p_{n+1,r}^{}}\,, \ldots, \ln\dfrac{p_{nr}^{}}{p_{n+1, r}^{}}\biggr\}$ 
and 
\vspace{0.75ex}
$\biggl\{\ln\dfrac{q_{1r}^{}}{q_{n+1,r}^{}}\,,\ldots, \ln\dfrac{q_{nr}^{}}{q_{n+1,r}^{}}\biggr\}$ 
are simple, for all $r\in I.$
Then for a topological conjunction of linear-fractional groups $PL^1$ and $PL^2$ 
\vspace{0.35ex}
it is necessary and enough existence of such permutation 
\vspace{0.35ex}
$\varrho\colon (1,\ldots, n+1) \to (1, \ldots, n+1)$ and complex number $\alpha $ with ${\rm Re}\, \alpha>{}-1,$ that either 
\\[1.75ex]
\mbox{}\hfill                        %  (6.5)
$
\dfrac{q_{\varrho (k)\;\! r}^{}}{q_{\varrho(n+1)\;\! r}^{}} = 
\dfrac{p_{kr}^{}}{p_{n+1, r}^{}}\, \biggl|\dfrac{p_{kr}^{}}{p_{n+1, r}^{}}\biggr|^{\alpha}
$
\ for all 
$
r\in I, 
\quad
k =1,\ldots,  n, 
$
\hfill {\rm (6.5)}
\\[1.5ex]
or}
\\[1.5ex]
\mbox{}\hfill                        %  (6.6)
$
\dfrac{q_{\varrho(k)\;\! r}^{}}{q_{\varrho (n+1)\;\! r}^{}} = 
\dfrac{\overline{p}_{kr}^{}}{\overline{p}_{n+1, r}^{}}\, 
\biggl|\dfrac{p_{kr}^{}}{p_{n+1, r}^{}}\biggr|^{\alpha}
$
\ for all 
$
r\in I, 
\quad
k=1,\ldots, n.
$
\hfill {\rm (6.6)}
\\[2.25ex]
\indent
{\bf Proof.} 
With the help of replacement $\xi (v) =T^{{}-1} f (Sv)$ for all $v\in\C P^n,$ from identities (6.1) we pass to identities 
\\[1.5ex]
\mbox{}\hfill                                         % (6.7)
$
\xi \bigl({\rm diag}\{p_{1r}^{}, \ldots, p_{n+1, r}^{}\}\;\! v\bigr) = 
{\rm diag}\{q_{1r}^{}, \ldots, q_{n+1, r}^{}\}\;\! \xi (v)
$ 
\ for all 
$
v\in\C P^n,
$
\ for all 
$
r\in I.
$
\hfill (6.7)
\\[1.5ex]
\indent
So, the topological conjunction of linear-fractional groups $PL^1$ and $PL^2$ 
is equivalent to performance of identities (6.7).

Thus, not belittling a generality, we will consider, that conjugating homeomorphism $\xi$ 
leaves invariant the common fixed points $O_{\tau}^{}$ 
\vspace{0.25ex}
(the origins of coordinates of the affine cards $M_{\tau}^{} = \{v,\;\! v_{\tau}^{} \ne 0\} $ 
\vspace{0.25ex}
of the atlas $M $ of variety $\C P^n), \ \tau =1,\ldots, n+1,$ of nondegenerate linear-fractional transformations 
\\[1.5ex]
\mbox{}\hfill                                         % (6.8)
$
{\rm diag}\{p_{1r}^{}, \ldots, p_{n+1, r}^{}\}\;\! v
$
\ for all 
$
v\in\C P^n,$
\ for all 
$
r\in I, 
$
\hfill (6.8)
\\[1ex]
and
\\[1ex]
\mbox{}\hfill                                         % (6.9)
$
{\rm diag}\{q_{1r}^{}, \ldots, q_{n+1, r}^{}\}\;\! v
$
\ for all 
$
v\in\C P^n,
$
\ for all 
$
r\in I, 
$
\hfill (6.9)
\\[1.75ex]
(that always we can achieve by nondegenerate linear-fractional transformation).
\vspace{0.35ex}

{\sl The necessity.} 
Let identities (6.7) are carried out. Nondegenerate linear-fractional tran\-s\-for\-ma\-ti\-ons (6.8) 
(linear-fractional transformations (6.9)) define on space $\C P^n$ invariant holomorphic foliati\-ons ${\frak C}_r^{}$ 
(invariant holomorphic foliations ${\frak D}_r^{})$ of complex dimension $1,$ defined by basis of nondegene\-rate absolute invariants 
\\[1.75ex]
\mbox{}\hfill
$
v_k^{{}^{\scriptstyle \ln (p_{nr}^{}/p_{n+1, r}^{})}}\;\!
v_n^{{}^{\scriptstyle \ln (p_{n+1, r}^{}/p_{kr}^{})}}\;\!
v_{n+1}^{{}^{\scriptstyle \ln (p_{kr}^{}/p_{nr}^{})}},
\quad 
k =1,\ldots, n-1,
$
\ for all $r\in I 
\hfill
$ 
\\[1.5ex]
(by basis of nondegene\-rate absolute invariants 
\\[1.5ex]
\mbox{}\hfill
$
v_k^{{}^{\scriptstyle \ln (q_{nr}^{}/q_{n+1, r}^{})}} 
v_n^{{}^{\scriptstyle \ln (q_{n+1, r}^{}/q_{kr}^{})}} 
v_{n+1}^{{}^{\scriptstyle \ln (q_{kr}^{}/q_{nr}^{})}},
\quad 
k =1,\ldots, n-1,
$
\ for all $r\in I).
\hfill
$ 
\\[1.5ex]
\indent 
Considering, that conjugating homeo\-morphism $\xi$ takes the layers of the foliations ${\frak C}_r^{}$ 
to the homeomorphic to them the layers of the foliations ${\frak D}_r^{}$ for all $r\in I,$ 
and also taking into consideration simplicity of sets of numbers from a theorem condition, we come to conclusion, 
that homeomorphism $ \xi $ translates Riemann spheres 
\\[1.5ex]
\mbox{}\hfill
$
\overline{\C}_{ls}^{}\colon v_k^{}=0,
\ \ 
k\ne l,
\ \ 
k\ne s,
\ \ 
s\ne l,
\ \ 
k =1,\ldots,  n+1,
\hfill
$ 
\\[1.5ex]
in Riemann spheres of a kind $\overline{\C}_{ij}^{},\  j\ne i,$ 
and all homeomorphisms of Riemann spheres simultaneously either keep, or change orientation.

Let homeomorphisms of Riemann spheres keep orientation 
(the case of homeomorphisms, changing orientation, it is considered similarly).
Then owing to a course of the proof of Theorems 3.2 and 6.2 we come to conclusion, that in the affine card 
$M_{n+1}^{} = \{v,\;\! v_{n+1}^{}\ne 0 \} $ relations 
\\[1.75ex]
\mbox{}\hfill
$ 
\dfrac{q_{\varrho(k)\;\! r}^{}}{q_{\varrho (n+1)\;\! r}^{}} = 
\dfrac{p_{kr}^{}}{p_{n+1, r}^{}}\, 
\biggl| \dfrac{p_{kr}^{}}{p_{n+1, r}^{}}\biggr|^{\alpha_k^{}},
\quad 
r =1,\ldots, \nu,
\quad 
k =1,\ldots, n,
\hfill
$ 
\\[2ex]
are fulfilled and at complex numbers 
\vspace{0.25ex}
$\alpha_k^{}$ the real parts ${\rm Re}\, \alpha_k^{}>{}-1,\ k =1,\ldots,  n,$ are carried out.
\vspace{0.25ex}
Taking into consideration last equalities, and also that fact, that in identities (6.7) vectors 
\vspace{0.35ex}
$\{p_{1r}^{},\ldots, p_{n+1, r}^{}\} $ and 
$\{q_{1r}^{}, \ldots, q_{n+1, r}^{}\} $ are defined to within scalar multipliers, for all $r\in I,$ 
\vspace{0.35ex}
and spending similar reasonings in other affine cards $M_{\tau}^{},\ \tau =1,\ldots, n,$ 
we come to conclusion, that $\alpha_k^{} =\alpha,\ k =1,\ldots, n.$ 
\vspace{0.35ex}
As a result we come to the relations (6.5).

{\sl The sufficiency} is proved by the construction of conjugating homeomorphism 
$\xi\colon\C P^n\to\C P^n$ such that its projection 
\\[1.5ex]
\mbox{}\hfill
$
\xi_{\varrho (k)}^{} (v) =\;\! \gamma_k^{}\;\!v_k^{}\;\!|v_k^{} |^{\alpha}
$ 
\ for all $v\in\C P^n,
\quad 
k =1,\ldots, n+1,
\hfill
$ 
\\[1.75ex]
if relations (6.5) take place; and its projections 
\\[1.75ex]
\mbox{}\hfill
$
\xi_{\varrho (k)}^{}(v) =\;\! \gamma_k^{}\;\!\overline {v}_k^{}\;\!|v_k^{} |^{\alpha}
$ 
\ for all $v\in\C P^n,
\quad 
k =1,\ldots, n+1,
\hfill
$ 
\\[1.75ex]
if parities (6.6) take place. \k
\vspace{0.35ex}

Consider now a topological conjunction of the non-Abelian linear-fractional groups $PL^1$ and $PL^2.$
\vspace{0.5ex}

{\bf Theorem 6.3.} 
{\it 
From a topological conjunction at $n=1$ of non-Abelian  linear-fractional groups $PL^1$ and $PL^2$ of general situation follows them 
$\R\!$-holomorphic conjunction which is carried out by  either nondegenerate linear-fractional, 
or nondegenerate antiholomorphic linear-fractional transfor\-mations.
}
\vspace{0.5ex}

{\bf Proof} of Theorem 6.3 directly follows from the following auxiliary statement.
\vspace{0.75ex}

{\bf Lemma 6.3.} 
{\it 
Let at $n=1$ and $I = \{1,2 \}$ 
\vspace{0.25ex}
linear-fractional groups $PL^1$ and $PL^2$ are topological conjugated, and conjgating homeomorphism 
$f $ is such that{\rm:}
\\[1ex]                                                   % (6.10)
\indent
{\rm 1)} $f (O_{\tau}^{}) =O_{\tau}^{},\ \tau =1,2;$ 
\hfill {\rm (6.10)}
\\[1.25ex]
\indent
{\rm 2)} either
\\[1.5ex]
\mbox{}\hfill                            % (6.11)
$
f (\lambda v_1^{}, v_2^{}) = 
\bigl(\lambda |\lambda |^{\alpha} f_1^{}(v), f_2^{} (v)\bigr)$ 
\ for all $v\in\C P^1,
\quad {\rm Re}\, \alpha>{}-1,
$
\hfill {\rm (6.11)}
\\[0.5ex]
or
\\[0.5ex]
\mbox{}\hfill                            % (6.12)
$
f (\lambda v_1^{}, v_2^{}) = 
\bigl(\;\!\overline{\lambda} | \lambda |^{\alpha} f_1 (v), f_2 (v)\bigr)$ 
\ for all $v\in\C P^1,
\quad 
{\rm Re}\, \alpha>{}-1; 
$
\hfill {\rm (6.12)}
\\[2.25ex]
\indent
{\rm 3)} $\!f (av_1^{}+bv_2^{}, cv_1^{}+dv_2^{}) = 
(Af_1^{} (v) +Bf_2^{} (v), Cf_1^{} (v) +Df_2^{} (v))\!
$
for all 
$
v\!\in\!\C P^1,
\, 
|b | + | c |\!>\! 0;
$
\\[1.75ex]
\indent
{\rm 4)} 
matrices 
$P=\left(\!\!
\begin{array}{cc}
a & b
\\
c & d
\end{array}
\!\!\right)=
S \left(\!\!
\begin{array}{cc}
\mu & 0
\\
0 & 1
\end{array}
\!\!\right) S^{{}-1}
$ 
\vspace{0.75ex}
and 
$
Q=\left(\!\!
\begin{array}{cc}
A & B
\\
C & D
\end{array}
\!\!\right)=
T \left(\!\!
\begin{array}{cc}
\theta & 0
\\
0 & 1
\end{array}
\!\!\right) T^{{}-1}$ are such that 
\vspace{0.75ex}
$|\mu|\ne 1,\ |\theta|\ne 1,\
S=\left(\!\!
\begin{array}{cc}
a_{\ast}^{} & b_{\ast}^{}
\\
c_{\ast}^{} & d_{\ast}^{}
\end{array}
\!\!\right),
\ 
T=\left(\!\!
\begin{array}{cc}
A_{\ast}^{} & B_{\ast}^{}
\\
C_{\ast}^{} & D_{\ast}^{}
\end{array}
\!\!\right),
\ 
\dfrac{B_{\ast}^{}}{D_{\ast}^{}}=
\Bigl(\dfrac{b_{\ast}^{}}{d_{\ast}^{}}\Bigr)^{1+\beta},
$
and
$
{\rm Re}\;\!(\beta-\alpha)\ln\Bigl|\dfrac{b_{\ast}^{}}{d_{\ast}^{}}\Bigr|-
{\rm Im}\;\! (\beta-\alpha)\;\! {\rm arg}\dfrac{b_{\ast}^{}}{d_{\ast}^{}}\ne 0;
$
\vspace{0.75ex}

{\rm 5)} the subgroup of the group $\C^{\ast}$ of nonzero complex numbers 
\vspace{0.5ex}
on the multiplication, formed by numbers $ \lambda $ and $\dfrac{b_{\ast}^{}}{d_{\ast}^{}}\,,$ 
is dense in the set $\C$ of complex numbers.

Then at {\rm (6.11)} this gomeomorphism has a kind 
\vspace{0.75ex}
$f (v) = \bigl(v_1^{}|v_1^{} |^{\alpha},\;\! v_2^{}|v_2^{}|^{\alpha}\bigr)$ for all $v\in\C P^1;$ 
and at {\rm (6.12)} this gomeomorphism has a kind 
$f (v) = \bigl(\;\!\overline{v}_1^{}|v_1^{} |^{\alpha},\;\! \overline {v}_2^{}|v_2^{} |^{\alpha}\bigr)$ 
\vspace{0.75ex}
for all $v\in\C P^1.$
}

{\sl Proof.} 
Owing to relations (6.10) --- (6.12) we come to conclusion, that $\ln |\mu| \ln | \theta |>  0.$

Let the identity (6.11) is carried out (a case, when the identity (6.12) takes place, is considered similarly). 
On the basis of (6.11) and conditions 3 of given lemma it is had following relations in the card $M_2^{}$ 
of the atlas $M$ of variety 
$\C P^1\colon$
\\[1.75ex]
\mbox{}\hfill
$
\psi\bigl(\lambda^k\;\! (P^m)^l\;\! v^2\bigr) = \lambda^k |\lambda |^{k\alpha} (Q^m)^l\;\!\psi (v^2)
$ 
\ for all $v^2\in M_2,
$ 
\ for all $k,l ,m\in\Z.
\hfill
$ 
\\[1.75ex]
Passing in them to a limit at 
\vspace{0.35ex}
$m\to{}-\infty$ if $ | \mu |> 1,$ and passing in them to a limit at 
$m\to {}+\infty $ if $ | \mu | <1,$ we receive, that 
\\[2ex]
\mbox{}\hfill                                        % (6.13)
$
\psi \biggl(\lambda^k \Bigl(\;\!\dfrac{b_{\ast}^{}}{d_{\ast}^{}}\Bigr)^{\!l}\;\!\biggr) = \;\!
\lambda^k\;\! |\lambda |^{k\alpha} \biggl(\dfrac{B_{\ast}^{}}{D_{\ast}^{}}\biggr)^{\!l}
$
\ for all 
$
k\in\Z,
$ 
\ for all 
$
l\in\Z.
$
\hfill (6.13)
\\[2ex]
From the condition 5 of lemma follows, that for any complex number $w\in\C$ there are such sequences 
$ \{k_s^{}(w)\} $ and $\{l_s^{}(w)\} $ of integers, that 
\\[1.75ex]
\mbox{}\hfill
$
\lim\limits_{s\to {}+\infty} \lambda^{k_s^{}(w)} \biggl(\dfrac{b_{\ast}^{}}{d_{\ast}^{}}\biggr)^{\!l_s^{}(w)} =w,
\quad
\lim\limits_{s\to {}+\infty} |k_s^{} (w) | = 
\lim\limits_{s\to {}+\infty} |l_s^{} (w) | = {}+ \infty. 
\hfill
$
\\[1.75ex]
From here on the basis of the relations (6.13) it is had, that
\\[1.75ex]
\mbox{}\hfill
$
\psi (v^2) =v^2\;\!|v^2 |^{\alpha}, 
\qquad
\lim\limits_{s\to {}+\infty} \biggl(\dfrac{b_{\ast}^{}}{d_{\ast}^{}}\biggr)^{(\beta-\alpha)\;\! l_s^{}(v^2)}
$
\ \, for all 
$
v^2\in M_2. 
\hfill
$
\\[1.75ex]
Now from the condition 5 of given lemma follows, that $\Bigl|\dfrac{b_{\ast}^{}}{d_{\ast}^{}}\Bigr| \ne 1,$ 
and from a condition 4 of given lemma follows, that 
\vspace{1ex}
$\lim\limits_{s\to {}+\infty} \biggl(\dfrac{b_{\ast}^{}}{d_{\ast}^{}}\biggr)^{\!(\beta-\alpha)\;\! l_s^{} (v^2)} =1.$ 
As a result it is received the first representation from a condition of the lemma 6.3. \k
\vspace{0.75ex}

{\bf Theorem 6.4.} 
\vspace{0.35ex}
{\it 
From a topological conjunction at $n> 1$ of non-Abelian  linear-fractional groups $PL^1$ and $PL^2$ 
of general situation follows them $\R\!$-holomorphic conjunction.
}
\vspace{0.5ex}

Justice of the given theorem follows from the following auxiliary statement.
\vspace{0.5ex}

{\bf Lemma 6.4.} 
{\it 
Let at $n> 1$ the matrices 
\\[1.5ex]
\mbox{}\hfill
$
P_r^{}=S_r^{}\;\!{\rm diag}\{p_{1r}^{}, \ldots, p_{n+1, r}^{}\}\;\! S_r^{{}-1},
\quad 
Q_r^{}=T_r^{}\;\! {\rm diag}\{q_{1r}^{}, \ldots, q_{n+1, r}^{}\}\;\! T_r^{{}-1},
\hfill
$
\\[1.5ex]
sets of numbers 
\vspace{0.75ex}
$
\biggl\{\ln\dfrac{p_{1r}^{}}{p_{n+1, r}^{}}\,,\ldots, \ln\dfrac{p_{nr}^{}}{p_{n+1,r}^{}}\biggr\}
$ 
and 
$
\biggl\{\ln\dfrac{q_{1r}^{}}{q_{n+1, r}^{}}\,, \ldots, \ln\dfrac{q_{nr}^{}}{q_{n+1, r}^{}}\biggr\}$ 
are simple, for all $r\in I.$
\vspace{0.25ex}
Then from a topological conjunction of non-Abelian linear-fractional groups $PL^1$ and $PL^2$ of general situation follows them 
$ \R\!$-holomorphic conjunction which is carried out by either nondegenerate linear-fractional, 
or nondegenerate antiholomorphic linear-fractional transformations.
}
\vspace{0.35ex}

{\sl Proof.} 
Let identities (6.1) are carried out. With the help of replacement 
$
\xi (v) =T_1^{{}-1} f (S_1^{}v)$ for all $v\in\C P^n,$ from them we pass to identities
\\[1.75ex]
\mbox{}\hfill                                                       % (6.14)
$
\xi\bigl({\rm diag}\{p_{11}^{},\ldots,p_{n+1,1}^{}\}\;\!v\bigr)=
{\rm diag}\{q_{11}^{},\ldots,q_{n+1,1}^{}\}\;\!\xi(v)
$
\ for all $v\in\C P^n.
$
\hfill (6.14)
\\[1.75ex]
\indent
Similarly, as well as at the proof of the theorem 6.2, 
we come to conclusion, that conjugating homeomorphism $\xi$ 
leaves invariant the common fixed points $O_{\tau}^{},\ \tau =1,\ldots, n+1,$ 
of the linear-fractional transformations (6.8) and (6.9) at $r=1.$
\vspace{0.35ex}

Let homeomorphisms of all Riemann spheres $\overline{\C}_{ls}^{},\ s\ne l,$ 
\vspace{0.25ex}
keep orientation (the case of homeomor\-phisms, changing of orientation, it is considered similarly). 
Then on the basis of the theorem 3.2, lemmas 3.1, 3.2 and a course of the proof of the theorem 6.2 
we come to the statement of the lemma 6.4. \k
\vspace{0.5ex}

Consider smooth, $\R\!$-holomorphic and holomorphic conjunctions of Abelian linear-fra\-c\-ti\-o\-nal groups $PL^1$ and $PL^2.$
\vspace{0.5ex}

{\bf Theorem 6.5.} 
{\it
Let the conditions of the theorem {\rm 6.1}  are satisfied.
Then for a smooth {\rm(}an $\R\!$-holomorphic{\rm)} conjunction of fractional-linear groups $PL^1$ and $PL^2$ 
it is necessary and enough, that either 
\\[0.75ex]
\mbox{}\hfill
$
\dfrac{q_{1r}^{}}{q_{2r}^{}} = \biggl(\dfrac{p_{kr}^{}}{p_{n+1, r}^{}}\biggr)^{\!\varepsilon}$
\ \ for all $r\in I,
\hfill
$ 
\\[1ex]
or 
\\[1ex]
\mbox{}\hfill
$
\dfrac{q_{1r}^{}}{q_{2r}^{}} = 
\biggl(\dfrac{\overline{p}_{kr}^{}}{\overline{p}_{n+1, r}^{}}\biggr)^{\!\varepsilon}$
\ \ for all $r\in I;
\quad 
\varepsilon^2=1.
\hfill
$
}
\\[1.5ex]
\indent
{\sl Proof} of the given statement is similar to the proof of the theorem 6.1 and is based on the theorem 3.4. \k
\vspace{0.35ex}

{\bf Theorem 6.6.} 
{\it 
Let the conditions of the theorem {\rm 6.2} are satisfied.
Then for a smooth {\rm(}an $\R\!$-holomorphic{\rm)} conjunction of fractional-linear groups 
$PL^1$ and $PL^2$ it is necessary and enough existence of such permutation 
$\varrho\colon (1,\ldots, n+1) \to (1, \ldots, n+1),$ that either 
\\[1.5ex]
\mbox{}\hfill
$
\dfrac{q_{\varrho(k)\;\! r}^{}}{q_{\varrho(n+1)\;\! r}^{}}=\dfrac{p_{kr}^{}}{p_{n+1,r}^{}}$ 
\ \, for all $r\in I,
\quad 
k=1,\ldots, n,
\hfill
$ 
\\[1.25ex]
or
\\[1.25ex]
\mbox{}\hfill
$
\dfrac{q_{\varrho(k)\;\!r}^{}}{q_{\varrho(n+1)\;\!r}^{}}=
\dfrac{\overline{p}_{kr}^{}}{\overline{p}_{n+1,r}^{}}
$
\ \, for all 
$
r\in I,
\quad 
k=1,\ldots, n.
\hfill
$
}
\\[1.75ex]
\indent
{\sl Proof} of the statements is carried out with use of a course of the proof of the theorem 6.2 by differentiation of identities (6.7). \k
\vspace{0.5ex}

Similarly to Theorems 6.5 and 6.6 it is received such statements.
\vspace{0.5ex}

{\bf Theorem 6.7.} 
{\it 
Let the conditions of the theorem {\rm 6.1}  are satisfied.
Then for a holomorphic conjunction of fractional-linear groups $PL^1$ and $PL^2$ 
it is necessary and enough, that} 
\\[1.75ex]
\mbox{}\hfill
$
\dfrac{q_{1r}^{}}{q_{2r}^{}}=\biggl(\dfrac{p_{kr}^{}}{p_{n+1,r}^{}}\biggr)^{\!\varepsilon}
$ 
\ \, for all 
$
r\in I,
\quad 
\varepsilon^2=1.
\hfill
$
\\[2ex]
\indent
{\bf Theorem 6.8.} 
{\it 
Let the conditions of the theorem {\rm 6.2}  are satisfied.
Then for a holomorphic conjunction of fractional-linear groups $PL^1$ and $PL^2$ 
\vspace{0.25ex}
it is necessary and enough existence of such permutation 
$\varrho\colon (1,\ldots,n+1)\to (1,\ldots,n+1)$, that
\\[1.75ex]
\mbox{}\hfill
$
\dfrac{q_{\varrho(k)\;\!r}^{}}{q_{\varrho(n+1)\;\!r}^{}}=
\dfrac{p_{kr}^{}}{p_{n+1,r}^{}}
$
\ \, for all
$
r\in I,
\quad 
k=1,\ldots, n.
\hfill
$
}
\\[2ex]
\indent
And, at last, we will consider smooth, $\R\!$-holomorphic and holomorphic conjunctions 
of non-Abelian linear-fractional groups $PL^1$ and $PL^2.$
\vspace{0.75ex}

{\bf Theorem 6.9.} 
{\it 
Let the conditions of the lemma {\rm 6.3} are satisfied.
Then for a smooth con\-j\-un\-c\-ti\-on of non-Abelian linear-fractional groups $PL^1$ and $PL^2$ 
of general situation it is ne\-ces\-sa\-ry and sufficient them conjunction, 
which is carried out by either nondegenerate linear-frac\-ti\-o\-nal transformation, 
\vspace{0.5ex}
or nondegenerate antiholomorphic linear-fractional transformation.
}

{\sl Proof} of the theorem 6.9 
is similar to the proof of the theorem 6.3 and is based on the following auxiliary statement.
\vspace{0.5ex}

{\bf Lemma 6.5.} 
{\it 
Let at $n=1$ and $I = \{1 \}$ linear-fractional groups $PL^1$ and $PL^2$ 
are smooth conjugated, and conjugating diffeomorphism $f $ is such that{\rm:}

{\rm 1)} relations {\rm (6.10)} are carried out{\rm;}

{\rm 2)} either
\\[1ex]
\mbox{}\hfill                                       % (6.15)
$
f(p_1^{} v_1^{},v_2^{})=\bigl(p_1^{} f_1^{}(v),f_2^{}(v)\bigr)
$
\ for all 
$
v\in\C P^1,
\quad
p_1^{}\ne 0,
\ \ 
p_1\ne 1;
$
\hfill {\rm(6.15)}
\\[0.5ex]
or
\\[0.5ex]
\mbox{}\hfill                                       % (6.16)
$
f(p_1^{} v_1^{},v_2^{})=
\bigl(\;\!\overline{p}_1^{} f_1^{}(v),f_2^{}(v)\bigr)
$
\ for all 
$
v\in\C P^1,
\quad 
p_1^{}\ne 0,
\ \ 
p_1^{}\ne 1.
$
\hfill {\rm(6.16)}
\\[2ex]
Then at {\rm (6.15)} this homeomorphism looks like 
\vspace{0.5ex}
$ 
f(v)=(a v_1^{},v_2^{})$ for all $v\in\C P^1;$ 
and at {\rm (6.16)} this homeomorphism looks like 
$
f(v)=\bigl(a\;\! \overline{v}_1^{}, \overline{v}_2^{}\bigr)
$ 
for all 
$
v\in\C P^1.
$
}
\vspace{0.75ex}

{\sl Proof.} 
Let the relations (6.10) are carried out. 
Differentiating in the card $M_2^{}$ of the atlas $M$ of the variety $\C P^1$ the identity (6.1) at $v^2=0,$ 
we receive equality 
\\[1.5ex]
\mbox{}\hfill
$
p_1^{}\;\!{\sf D}_{{}_{\scriptstyle v^2}}\psi(0)\, dv^2+\,
\overline{p}_1^{}\;\!{\sf D}_{{}_{\scriptstyle \overline{v}{}^{\;\!2}}}\psi(0)\,d\;\!\overline{v}{}^{\;\!2}=
q_1^{}\bigl( {\sf D}_{{}_{\scriptstyle v^2}}\psi(0)\,dv^2+
{\sf D}_{{}_{\scriptstyle \overline{v}{}^{\;\!2}}}\psi(0)\, d\;\!\overline{v}{}^{\;\!2}\bigr).
\hfill
$ 
\\[1.75ex]
Owing to that $f $ is a diffeomorphism, we have, that 
\vspace{0.5ex}
$\bigl|{\sf D}_{{}_{\scriptstyle v^2}}\psi(0)\bigr|+\bigr|{\sf D}_{{}_{\scriptstyle \overline{v}{}^{\;\!2}}}\psi(0)\bigr|>0.
$ 
Therefore from last equality we come either to the relation $q_1^{}=p_1^{},$ 
or to the relation $q_1^{} =\overline{p}_1^{}.$
\vspace{0.5ex}

In the first case the identity (6.15) takes place. From it we receive identities 
\\[1.5ex]
\mbox{}\hfill
$
\psi\bigl(p_1^{\,l}v^2\bigr)=p_1^{\,l}\;\!\psi(v^2)
$ 
\ 
for all $v^2\in M_2^{},$ 
\ for all $l\in\Z,
\hfill
$ 
\\[1.75ex]
holomorphic differentiating which on $v^2,$ we come to identities 
\\[1.75ex]
\mbox{}\hfill
$
{\sf D}_{{}_{\scriptstyle v^2}}\psi\bigl(p_1^{\,l}v^2\bigr)=
{\sf D}_{{}_{\scriptstyle v^2}}\psi(v^2)
$ 
\ for all $v^2\in M_2^{},$ 
\ for all $l\in\Z.
\hfill
$
\\[1.75ex]
Owing to that transformation $\psi$ is a diffeomorphism, we receive identity 
\\[1.75ex]
\mbox{}\hfill
$
{\sf D}_{{}_{\scriptstyle v^2}}\psi(v^2)=a$ 
\ for all $v^2\in M_2^{}.
\hfill
$ 
\\[1.5ex]
\indent
From here taking into account the relation (6.15) we come to the first representation from the given lemma.

Similarly in the second case 
\vspace{0.75ex}
it is received the second representation from Lemma 6.5. \k

{\bf Theorem 6.10.} 
{\it 
Let the conditions of the lemma {\rm 6.4} are satisfied.
Then for a smooth conjunction of non-Abelian linear-fractional groups $PL^1$ and $PL^2$ 
of general situation  it is necessary and sufficient them conjunction, 
which is carried out by either nondegenerate linear-fractional transformation, 
\vspace{0.5ex}
or nondegenerate antiholomorphic linear-fractional transformation.
}

{\sl Proof. The necessity.} 
As well as at the proof of the lemma 6.4, from identities (6.1) we pass to identities (6.14). 
Further owing to the theorem 6.2 we come to conclusion, that take place either relations (6.5) at $r=1$, or relations (6.6) at $r=1.$

Let relations (6.5) are carried out at $r=1$ (a case of performance of relations (6.6) at $r=1$ is considered similarly). 
We will consider the affine card $M_{n+1}^{} $ of the atlas $M$ of projective space $\C P^n.$ 
Owing to a course of the proof of the theorem 6.2 and theorems 3.4 smooth conjunction 
of narrowings of linear-fractional groups $PL^1$ and $PL^2$ on this card is carried out by nondegenerate linear $\R\!$-holomorphic transformation. By direct calculations on the basis of parities (6.5) at $r=1$ we are convinced, that $ \alpha=0$ and that nondegenerate linear transformation is holomorphic, as keeping eigenvalues of matrices. With the help of linear-fractional functions of transition between affine cards of the atlas $M $ on the basis of the aforesaid nondegenerate linear transformation we receive nondegenerate linear-fractional transformation. 
It also will be required conjugating diffeomorphism.

{\sl The sufficiency} is checked by direct calculations.
\vspace{0.5ex}

Similarly to theorems 6.9 and 6.10 it is received following statements.
\vspace{0.75ex}

{\bf Theorem 6.11.} 
{\it 
Let the conditions of the lemma {\rm 6.3} are satisfied.
Then for a holomorphic conjunction of non-Abelian linear-fractional groups $PL^1$ and $PL^2$ 
of general situation  it is necessary and sufficient them conjunction, which is carried out by nondegenerate linear-fractional transformation.
}
\vspace{0.75ex}

{\bf Theorem 6.12.} 
{\it 
Let the conditions of the lemma {\rm 6.4}  are satisfied.
Then for a holomorphic conjunction of non-Abelian linear-fractional groups 
$PL^1$ and $PL^2$ of general situation  it is necessary and sufficient them conjunction, 
which is carried out by nondegenerate linear-fractional transformation.
}
\\[3.75ex]
\centerline{
\large\bf  
7. Applications to the complex nonautonomous projective 
}
\\[0.35ex]
\centerline{
\large\bf  
matrix Riccati  equations
}
\\[1.5ex]
\indent
Theorems 6.1 --- 6.3, 6.5 --- 6.12,  Lemmas 6.3 and 6.4 give the chance, 
being based on Theorems 1.1 --- 1.3 to make topological, smooth, $\R\!$-holomorphic, 
and holomorphic classifications of complex  nonautonomous homogeneous projective matrix Riccati equations of a kind (5.1), 
and from Theorem 6.4 is received the following statement.
\vspace{0.35ex}

{\bf Theorem 7.1.} 
{\it 
From topological equivalence of complex nonautonomous homogeneous projective matrix Riccati equations 
with non-Abelian phase groups of general situation follows them $\R\!$-holomor\-phic equivalence.
}
\vspace{0.35ex}

As well as in point 4, from this theorem it is had, that complex nonautonomous homogeneous projective 
matrix Riccati equations with coefficients, holomorphic on path connected holomorphic varieties with non-Abelian fundamental groups, 
are structurally unstable.

Consider ordinary homogeneous projective scalar Riccati equations 
(i.e. ordinary homogeneous projective matrix Riccati equations at $n=1)$ 
\\[1.5ex]
\mbox{}\hfill                                                  % (7.1)
$
\dfrac{dv}{dz}={\frak A}_1^{}(z)\;\!v
$
\hfill (7.1)
\\[1.25ex]
and
\\[1.25ex]
\mbox{}\hfill                                                  % (7.2)
$
\dfrac{dv}{dz}={\frak B}_1^{}(z)\;\!v,
$
\hfill (7.2)
\\[2.15ex]
where square matrices ${\frak A}_1^{}(z)=\|{\frak a}_{lk1}^{}(z)\|$ and
\vspace{0.75ex}
${\frak B}_1^{}(z)=\|{\frak b}_{lk1}^{}(z)\|$ of the size $2$ consist from holomor\-phic functions 
${\frak a}_{lk1}^{}\colon\Gamma_1^{\nu_1^{}}\to\C$ and ${\frak b}_{lk1}^{}\colon\Gamma_2^{\nu_2^{}}\to\C,\ l=1,2,\ k=1,2.$ 
\vspace{0.35ex}
The common solutions of ordinary homogeneous projective scalar Riccati equations
\vspace{0.25ex}
(7.1) and (7.2) define covering foliation $ \PP^3$ and $\PP^4$, accordingly, on the varieties 
\vspace{0.35ex}
$\C P^1\times\Gamma_1^{\nu_1^{}}$ and $\C P^n\times\Gamma_2^{\nu_2^{}}.$

The phase group $Ph (\PP^3)$ of covering foliation $\PP^3$ 
\vspace{0.25ex}
is generated by nondegenerate linear-fractional transformations 
\vspace{0.25ex}
$P_r^{}v$ for all $v\in\C P^1,\ P_r^{}\in GL(2,\C),\ r=1,\ldots, \nu_1^{},$ 
and phase group $Ph (\PP^4) $ of covering foliation $ \PP^4$ 
\vspace{0.35ex}
is generated by nondegenerate linear-fractional transformations 
$Q_r^{}v$ for all $v\in\C P^1,\ Q_r^{}\in GL(2,\C),\ r=1,\ldots, \nu_2^{}.$
\vspace{0.5ex}

In a commutative case on the basis of Theorems 6.1 and 1.2 we receive the statement.
\vspace{0.35ex}

{\bf Theorem 7.2.} 
{\it 
Let matrices 
\\[1.5ex]
\mbox{}\hfill
$
P_r^{}=S\;\! {\rm diag} \{p_{1r}^{},p_{2r}^{}\}\;\! S^{{}-1}, \, r=1,\ldots,\nu_1^{},
\ \
Q_r^{}=T\;\! {\rm diag}\{q_{1r}^{},q_{2r}^{}\}\;\! T^{{}-1},\, r= 1,\ldots, \nu_2^{}, \ \nu_1^{}\leq\nu_2^{}.
\hfill
$
\\[1.5ex]
Then for embedding of ordinary homogeneous projective scalar Riccati equation {\rm (7.1)} 
in ordinary homogeneous projective scalar Riccati equation {\rm (7.2)} it is necessary and enough existence of such bijective map 
$\!\varkappa\colon (1,\ldots,\nu_1^{})\!\to\!\Delta\!$ and complex number $\!\alpha\!$ with ${\rm Re}\;\!\alpha\ne\! {}-1,\!$ 
that either
\\[1.75ex]
\mbox{}\hfill
$
\dfrac{q_{1\varkappa(r)}^{}}{q_{2\varkappa(r)}^{}}=
\dfrac{p_{1r}^{}}{p_{2r}^{}}\, \biggl|\dfrac{p_{1r}^{}}{p_{2r}^{}}\biggr|^{\alpha},
\quad 
r=1,\ldots, \nu_1,
\hfill
$
\\[1.5ex]
or 
\\[1.5ex]
\mbox{}\hfill
$
\dfrac{q_{1\varkappa(r)}^{}}{q_{2\varkappa(r)}^{}}=
\dfrac{\overline{p}_{1r}^{}}{\overline{p}_{2r}^{}}\,
\biggl|\dfrac{p_{1r}^{}}{p_{2r}^{}}\biggr|^{\alpha},
\quad 
r=1,\ldots, \nu_1^{},
\hfill
$ 
\\[2ex]
where the set $\Delta$ consists from $\nu_1^{}$ numbers  $\{1, \ldots, \nu_2^{}\}.$
}
\vspace{0.75ex}

In a noncommutative case on the basis of  Lemma 6.3 and Theorem 1.2 it is 
had concrete constructive criteria of embedding of ordinary homogeneous projective scalar  equations.

And, at last, we will consider ordinary homogeneous projective scalar Riccati equations (7.1) and (7.2) 
in a case when square matrices ${\frak A}_1^{}(z)=\|{\frak a}_{lk1}^{}(z)\|$ and
\vspace{0.5ex}
${\frak B}_1^{}(z)=\|{\frak b}_{lk1}^{}(z)\|$ of the size of $2$ consist of $1\!$-periodic holomorphic functions 
${\frak a}_{lk1}^{}\colon\C\to\C$ and ${\frak b}_{lk1}^{}\colon\C\to\C,$ and, besides, 
functions ${\frak b}_{lk1}^{}(z)$ are such that 
\\[1.5ex]
\mbox{}\hfill
${\frak b}_{lk1}^{}(x+i)={\frak b}_{lk1}^{}(x)$ for all $x\in [0,1],\ l=1,2,\ k=1,2.
\hfill
$ 
\\[1.5ex]
\indent
The common solutions of ordinary homogeneous projective scalar 
\vspace{0.25ex}
Riccati equations (7.1) and (7.2) define covering foliations $\PP^5$ and $ \PP^6,$ 
\vspace{0.25ex}
accordingly, on the varieties $\C P^1\times Z$ and $\C P^n\times T^2,$ 
where $T^2$ is a torus, defined by evolvment 
\\[1.5ex]
\mbox{}\hfill
$
K=\{z=x+iy\in\C\colon x\in [0,1],\ y\in [0,1]\}.
\hfill
$
\\[1.5ex]
\indent
The phase group $Ph (\PP^5) $ of covering foliation $ \PP^5$ 
\vspace{0.25ex}
is generated by nondegenerate linear-fractional transformation 
\vspace{0.25ex}
$P_1^{}v$ for all $v\in\C P^1,\ P_1^{}\in GL(2,\C),$ and the phase group $Ph (\PP^6) $ 
of covering foliation $ \PP^6$ is 
\vspace{0.35ex}
generated by nondegenerate linear-fractional transformations 
$Q_r^{}v$ for all $v\in\C P^1,\ Q_r^{}\in GL(2,\C),\ r=1,2.$ 
\vspace{0.35ex}
We will consider, that the cylinder $Z$ covers the torus $T^2,$ 
under a condition, that the band of $0 \leq {\rm Re}\, z\leq 1$ covers the square $K.$
\vspace{0.25ex}

Owing to Theorems 6.1 and 1.3 we have the statement.
\vspace{0.75ex}

{\bf Theorem 7.3.} 
{\it 
Let matrices 
\\[1.25ex]
\mbox{}\hfill
$
P_1^{}=S\;\!{\rm diag} \{p_{11}^{},p_{21}^{}\}\;\! S^{{}-1},
\quad
Q_r^{}=T\;\! {\rm diag}\{q_{1r}^{},q_{2r}^{}\}\;\!T^{{}-1},
\ \ r=1,2.
\hfill
$
\\[1.5ex]
Then for covering of the  projective scalar Riccati equation {\rm (7.2)} 
by the projective scalar Ric\-ca\-ti equation {\rm (7.1)} 
it is necessary and enough existence of such index $r\in \{1,2\} $ 
and cor\-res\-pon\-ding to it complex numbers $\alpha_r^{}$ with ${\rm Re}\, \alpha_r^{}\ne {}-1,$ that either
\\[2ex]
\mbox{}\hfill
$\;\!
\dfrac{q_{1r}^{}}{q_{2r}^{}}=
\dfrac{p_{11}^{}}{p_{21}^{}}\,\biggl|\dfrac{p_{11}^{}}{p_{21}^{}}\biggr|^{\alpha_r^{}},
$
\quad 
or
\quad
$
\dfrac{q_{1r}^{}}{q_{2r}^{}}=
\dfrac{\overline{p}_{11}^{}}{\overline{p}_{21}^{}}\, 
\biggl|\dfrac{p_{11}^{}}{p_{21}^{}}\biggr|^{\alpha_r^{}}.
\hfill
$
}
\\[4.75ex]
\centerline{
\large\bf  
8. Phase groups of covering foliations, defined by real nonautonomous 
}
\\[0.35ex]
\centerline{
\large\bf  
linear differential systems 
}
\\[1.5ex]
\indent
We will consider linear differential systems
\\[1.5ex]
\mbox{}\hfill                                            % (8.1)
$
\displaystyle
dx=\sum\limits_{j=1}^m A_j^{}(t_1^{},\ldots,t_m^{})\;\!x\ dt_j^{} 
$
\hfill (8.1)
\\[1ex]
and
\\[1.25ex]
\mbox{}\hfill                                            % (8.2)
$
\displaystyle
dx=\sum\limits_{j=1}^m B_j^{}(t_1^{},\ldots,t_m^{})\;\!x\ dt_j^{}\;\!, 
$
\hfill (8.2)
\\[1.5ex]
ordinary at $m=1$ and completely solvable at $m> 1,$ 
\vspace{0.35ex}
where $x = (x_1^{}, \ldots, x_n^{}),$ square matrices 
$A_j^{}(t_1^{},\ldots,t_m^{})=\|a_{ikj}^{}(t_1^{},\ldots,t_m^{})\|$ and 
\vspace{0.35ex}
$B_j^{}(t_1^{},\ldots,t_m^{})=\|b_{ikj}^{}(t_1^{},\ldots, t_m^{})\|$ of the size $n$ 
consist from holomorphic functions 
\\[1.5ex]
\mbox{}\hfill
$
a_{ikj}^{}\colon A\to\R
$ 
\, and \,
$
b_{ikj}^{}\colon B\to\R,
\quad 
i=1,\ldots, n,
\ \ k=1,\ldots, n, 
\ \ j=1,\ldots, m,
\hfill
$ 
\\[1.5ex]
path connected holomorphic varieties $A$ and $B$  are holomorphically equivalent each other, 
fundamental groups $\pi_1^{}(A)$ and $\pi_1^{}(B)$ have final number $\nu\in\N$ of the forming.
\vspace{0.25ex}

The common solutions of linear differential systems (8.1) and (8.2) define covering foliations $\Ll^9$ and $ \Ll ^ {10},$ 
accordingly, on varieties $\R^n\times A$ and $\R^n\times B.$ 
\vspace{0.25ex}

We will say, that the linear differential systems (8.1) and (8.2) are 
{\it topologically {\rm(}smoothly,  holomorphically}) {\it equivalent} 
\vspace{0.25ex}
if exists the homeomorphism (the diffeomorphism, the holo\-morphism) 
$h\colon\R^n\times A\to\R^n\times B,$ 
\vspace{0.25ex}
translating the layers of the covering foliation $\Ll^9$ in the layers of the covering foliation $ \Ll^{10}.$ 
Similarly we introduce the concepts of {\it embedding {\rm(}smooth embedding,  holomorphic embedding}) and 
{\it covering {\rm(}smoothly covering,  holomorphically covering}) of linear differential systems.
\vspace{0.25ex}

The phase group $Ph (\Ll^9)$ 
\vspace{0.25ex}
of the covering foliation $ \Ll^9$ is generated by the forming nondegenerate 
\vspace{0.25ex}
linear transformations $P_r^{}w$ for all $x\in\R^n,\ P_r^{}\in GL(n,\R),\ r=1,\ldots, \nu,$ 
and the phase group $Ph (\Ll^{10})$ of the covering foliation $\Ll^{10}$ 
\vspace{0.25ex}
is generated by the forming nondegenerate 
linear transformations $Q_r^{}x$ for all $x\in\R^n,\ Q_r^{}\in GL(n,\R),\ r=1,\ldots, \nu.$
\\[3.75ex]
\centerline{
\large\bf  
9. Conjunctions of linear actions on $\R^n$
}
\\[1.5ex]
\indent
We will consider a problem about a finding of necessary and sufficient conditions of existence such homeomorphism 
(diffeomorphism, holo\-morphism) $f\colon\R^n\to\R^n,$ that identities   
\\[1.5ex]
\mbox{}\hfill                              % (9.1)
$
f(P_r^{}x)=Q_r^{}f(x)
$ 
\ for all $x\in\R^n,
$
\ for all 
$
r\in I, 
$
\hfill (9.1)
\\[1.75ex]
take place, where $f (x)\!=\!(f_1^{}(x), \ldots, f_n^{}(x)),\!$ 
\vspace{0.5ex}
square matrices $P_r^{}, Q_r^{}\!\in GL(n, \R)$ for all $r\in\! I.$ 

Group of linear actions on $\R^n,$ formed by matrices $P_r^{}$ for all $r\in I,$ 
we will designate through $L^3,$ and through $L^4$ we will designate the similar group, 
formed by matrices $Q_r^{}$ for all $r\in I.$ 
Besides, further everywhere {\it strongly hyperbolic} we will name matrices at which all own values 
are various among themselves and on the module are distinct from $1.$
\vspace{0.5ex}

{\bf Theorem 9.1.} 
{\it 
For the topological conjunction {\rm (9.1)} at $n=1$ of linear groups $L^3$ and $L^4$ it is necessary and enough, that}
\\[1.5ex]
\mbox{}\hfill                              % (9.2)
$
q_r^{}=p_r^{}\;\! |p_r^{}|^{\alpha}
$
\ {\it for all}
$
r\in I,
\quad
\alpha>{}-1.
$
\hfill (9.2)
\\[2ex]
\indent
{\sl Proof. The necessity.} 
We will consider at first a case $I = \{1\}.$ 
If $p_1^{}=1,$ that from (9.1) follows, that $q_1^{}=1,$ and, hence, the relation (9.2) is carried out.
\vspace{0.35ex}

If $p_1^{} ={}-1,$ that on the basis of (9.1) it is had, that 
\vspace{0.25ex}
$f\bigl(p_1^{2}x\bigr) =q_1^2f (x)$ for all $x\in\R,$ 
and, therefore $q_1^2=1.$ At $q_1^{}=1$ 
\vspace{0.25ex}
the coordination of orientations of maps in different parts of identity (9.1) 
is broken, owing to what $q_1^{}={}-1.$ Hence, and in this case the relation (9.2) takes place.

Let now $|p_1^{}|\ne 1.$ 
\vspace{0.5ex}
Then from (9.1) follows, as $|q_1^{}|\ne 1.$ And, means, exists 
$ \alpha\ne {}-1,$ that $q_1^{}=p_1^{}|p_1^{}|^{\alpha}$ or $q_1^{}={}-p_1^{}|p_1^{}|^{\alpha}.$
\vspace{0.5ex}

If $q_1^{}=p_1^{}|p_1^{}|^{\alpha},$ that of (9.1) follows, that 
\vspace{0.5ex}
$f\bigl(p_1^{2k}\bigr) =p_1^{2k (1 +\alpha)} f (1)$ for all $k\in\Z.$ 
We will admit, that $ \alpha <{}-1.$ 
\vspace{0.35ex}
Now we will pass in this equality to a limit: in the case 
$ |p_1^{}|> 1$ at $k\to {}+ \infty,$ and in the case $|p_1^{}| <1$ at $k\to {}-\infty.$ 
\vspace{0.35ex}
Every time we will receive the contradiction. Therefore in the case $q_1^{}=p_1^{}|p_1^{}|^{\alpha}$ number $\alpha>{}-1.$
\vspace{0.5ex}

If $q_1^{} ={}-p_1^{}|p_1^{}|^{\alpha},$ 
\vspace{0.5ex}
that is broken a coordination of orientation of maps in different parts of identity (9.1).
Uniting considered above possibility, we come to conclusion about justice of a relation (9.2) at $I = \{1\}.$

Consider now the case $I\ne \{1\}.$ 
\vspace{0.35ex}
Owing to the proof of the previous part of the given statement it is had, that 
\vspace{0.5ex}
$q_r^{}=p_r^{}|p_r^{}|^{\alpha_r^{}},\ \alpha_r^{}>{}-1,$  for all $r\in I,$ and thus 
$ |p_r^{}| = 1$ in only case when when $|q_r^{}| = 1$ for all $I.$ 
\vspace{0.25ex}

Through $I_1^{}$ we will designate set of such indexes $r,$ that $|p_r^{}| = |q_r^{}| = 1.$ 
\vspace{0.35ex}

We will consider 3 logic possibilities, when addition $CI_1^{}$ of sets $I_1^{}$ to set $I\colon$ 
1) is empty; 2) consists of one index; 3) consists of more than one index.
\vspace{0.35ex}

If $CI_1^{} =\O,\ I=I_1^{},$ and $|p_r^{}| = |q_r^{}| = 1$ for all $r\in I,$ 
\vspace{0.5ex}
that relations (9.1) take place at any real $\alpha.$ 
If $CI_1^{}=\{r\},$ that in relations (9.1) can be put $\alpha =\alpha_r^{}.$
\vspace{0.5ex}

Let $|p_{r_1^{}}^{}|\ne 1,\ |p_{r_2^{}}^{}|\ne 1,\ r_1^{}\ne r_2^{},\ r_1^{}\in CI_1^{},\ r_2^{}\in CI_1^{}.$ 
Then from (9.1) follows, that 
\\[1.75ex]
\mbox{}\hfill
$
f\bigl(p_{r_1^{}}^l p_{r_2^{}}^n x\bigr)= q_{r_1^{}}^l q_{r_2^{}}^n f(x)
$ 
\ for all $x\in\R,$
\ for all $l\in\Z,$ 
\ for all $n\in\Z.
\hfill
$ 
\\[1.75ex]
Owing to that the set of rational numbers is everywhere dense on set of real numbers, 
we conclude about existence of such sequences $\{l_s^{}\}$ and
$\{n_s^{}\}$ of integers, that 
\\[2ex]
\mbox{}\hfill
$
\lim\limits_{s\to {}+\infty}p_{r_1^{}}^{2l_s^{}} p_{r_2^{}}^{2n_s^{}}=1,
\qquad 
\lim\limits_{s\to {}+\infty}|l_s^{}|=\lim\limits_{s\to {}+\infty}|n_s^{}|= {}+\infty.
\hfill
$ 
\\[1.5ex]
Then 
\\[1.5ex]
\mbox{}\hfill
$
f(x)= \lim\limits_{s\to {}+\infty}|p_{r_1^{}}|^{2l_s^{}(\alpha_{r_1^{}}- \alpha_{r_2^{}})}f(x)
$ 
\ for all 
$
x\in\R.
\hfill
$ 
\\[2ex]
\indent
From here follows, that $\alpha_{r_1^{}}^{} = \alpha_{r_2^{}}^{} = \alpha^{\ast}>{}-1.$ 
Therefore in relations (9.2) for $\alpha$ it is possible to take $\alpha^{\ast}.$

{\sl The sufficiency} is proved by construction of conjugating homeomorphism 
\\[1.5ex]
\mbox{}\hfill
$
f\colon x\to x|x|^{\alpha}
$  
\ for all 
$
x\in\R,
\quad \alpha>{}-1.\ \k
\hfill
$
\\[2ex]
\indent
{\bf Theorem 9.2.} 
{\it 
Let at $n>1$ and $I = \{1\}$ 
\vspace{0.25ex}
real normal Jordan form of strongly hyperbolic matrix $P_1^{}\in GL(n,\R) $ looks like 
\\[1.5ex]
\mbox{}\hfill
$
J(P_1^{})=
\left(\!\!
\begin{array}{cc}
J_s^{}(P_1^{}) & 0
\\[1ex]
0 & J_u^{}(P_1^{})
\end{array}
\!\!\right),
\hfill
$ 
\\[1.75ex]
where all eigenvalues of the matrix $J_s^{} (P_1^{})$ 
on the module are less than $1,$ and all eigenvalues of the matrix $J_u^{}(P_1^{})$ 
on the module are more than $1,\ p={\rm dim}\, J_s^{}(P_1^{});$ 
\vspace{0.25ex}
real normal Jordan form of strongly hyperbolic matrix 
$Q_1^{}\in GL (n,\R)$ looks like 
\\[1.75ex]
\mbox{}\hfill
$
J(Q_1^{})= 
\left(\!\!
\begin{array}{cc}
J_s^{}(Q_1^{}) & 0
\\[1ex]
0 & J_u^{}(Q_1^{})
\end{array}\!\!\right),
\hfill
$ 
\\[2ex]
where all eigenvalues of the matrix $J_s^{}(Q_1^{})$ 
\vspace{0.35ex}
on the module are less than $1,$ and all eigenvalues of the matrix $J_u^{}(Q_1^{})$ 
\vspace{0.35ex}
on the module are more than $1,\ q={\rm dim}\, J_s^{}(Q_1^{}).$
Then linear groups $L^3$ and $L^4$ are topological conjugated in only if, when 
\\[1.5ex]
\mbox{}\hfill
$
p=q,\quad 
{\rm det}\, J_s^{}(P_1^{})\, {\rm det}\, J_s^{}(Q_1^{})>0,
\quad 
{\rm det}\, J_u^{}(P_1^{})\, {\rm det}\, J_u^{}(Q_1^{})>0.
\hfill
$
}
\\[2ex]
\indent
{\sl Proof.} Let the identities (9.1) take place at $I = \{1\}.$ 

It is easy to see, that for linear actions the dimensions of stable both unstable invariant subspaces and orientations (positive or negative) of narrowings on these invariant subspaces are invariant at a topological conjuction. Therefore for the proof of the theorem it is enough to show a conjunction between real normal Jordan forms 
\vspace{0.25ex}
of the above-stated linear actions and the linear action defined by one of canonical matrices 
\vspace{0.25ex}
${\rm diag}\{\varepsilon e^{{}-1},e^{{}-1},\ldots,e^{{}-1}, e,\ldots,e,\delta e\},$ 
where $\varepsilon^{2}=1,\ \delta^{2}=1$ (such linear action we will name canoni\-cal).
\vspace{0.25ex}

For proof end we will result conjugating to canonical linear actions homeomorphisms of spaces 
$\R$ and $\R^2$ for narrowings on stable and unstable invariant subspaces of the linear action, 
corresponding to strongly hyperbolic matrix:

1) homeomorphism
\\[1ex]
\mbox{}\hfill
$
x_1^{}\to\  x_1^{}|x_1^{}|^{{}^{\scriptstyle {}-\tfrac{1}{\ln\lambda_1^{}}-1}}
$
\ for all 
$
x_1^{}\in\R
\hfill
$ 
\\[2ex]
conjugates the linear action $\lambda_1^{}x_1^{}$ for all $x_1^{}\in\R,\ 0<\lambda_1^{}<1,$ 
with the linear action $e^{{}-1}x_1^{}$ for all $x_1^{}\in\R;$
\vspace{0.25ex}

2) homeomorphism 
\\[1ex]
\mbox{}\hfill
$
x_1^{}\to\ x_1^{}|x_1^{}|^{{}^{\scriptstyle {}-\tfrac{1}{\ln|\lambda_1^{}|}-1}}
$
\ for all $x_1\in\R
\hfill
$ 
\\[2ex]
conjugates the 
\vspace{0.25ex}
linear action $\lambda_1^{}x_1^{}$ for all $x_1^{}\in\R,\ {}-1<\lambda_1^{}<0,$ 
with the linear action ${}-e^{{}-1}x_1^{}$ for all $x_1^{}\in\R;$
\vspace{0.35ex}

3) homeomorphism
\\[1ex]
\mbox{}\hfill
$
x_1^{}\to\  x_1^{}|x_1^{}|^{{}^{\scriptstyle \tfrac{1}{\ln\lambda_1^{}}-1}}
$ 
\ for all 
$
x_1^{}\in\R
\hfill
$ 
\\[2ex]
conjugates the linear action $\!\lambda_1^{}x_1^{}\!$ for all $\!x_1^{}\!\in\!\R, \lambda_1^{}\!\!>\!1,\!$ 
\vspace{0.25ex}
with the linear action $\!ex_1^{}\!$ for all $\!x_1^{}\!\!\in\!\R;$

4) homeomorphism
\\[1ex]
\mbox{}\hfill
$
x_1^{}\to\ x_1^{}|x_1^{}|^{{}^{\scriptstyle \tfrac{1}{\ln |\lambda_1^{}|}-1}}
$ 
\ for all 
$
x_1^{}\in\R
\hfill
$ 
\\[2ex]
conjugates the linear action $\lambda_1^{}x_1^{}$ for all $x_1^{}\in\R,\
\lambda_1^{}<{}-1,$ with the linear action ${}-ex_1^{}$ for all $x_1^{}\in\R;$

5) homeomorphism 
\\[1.5ex]
\mbox{}\hfill
$
(x_1^{},x_2^{})\to \
\Bigl(x_1^{}\cos\bigl(\pi\ln\sqrt{x_1^2+x_2^2}\,\bigr) - x_2^{}\sin\bigl(\pi\ln\sqrt{x_1^2+x_2^2}\,\bigr), \
x_1^{}\sin\bigl(\pi\ln\sqrt{x_1^2+x_2^2}\,\bigr) \ +
\hfill
$
\\[1.75ex]
\mbox{}\hfill
$
+ \ x_2^{}\cos\bigl(\pi\ln\sqrt{x_1^2+x_2^2}\,\bigr)\Bigr)
$
\ for all 
$
(x_1^{},x_2^{})\in\R^2\backslash \{(0,0)\},
\quad 
(0,0)\to (0,0) 
\hfill
$ 
\\[2ex]
conjugates the linear action ${}-e^{{}-1}Ix$ for all $x\in\R^2,$ 
where $I$ there is an identity matrix of the second order, with the linear action 
$e^{{}-1}Ix$ for all $x\in\R^2;$ 
and also conjugates the linear action ${}-eIx$ for all $x\in\R^2$ with the linear action $eIx$ for all $x\in\R^2;$
\vspace{0.35ex}

6) homeomorphism 
\\[1.75ex]
\mbox{}\hfill
$
(x_1^{},x_2^{})\to \ 
\Biggl(
x_1^{}\cos\biggl(\dfrac{{\rm arctg}\dfrac{\beta}{\alpha}}{\ln\sqrt{\alpha^2+\beta^2}}\, \ln\sqrt{x_1^2+x_2^2}\,\biggr) - 
x_2^{}\sin\biggl(\dfrac{{\rm arctg}\dfrac{\beta}{\alpha}}{\ln\sqrt{\alpha^2+\beta^2}}\ \ln\sqrt{x_1^2+x_2^2}\, \biggr),
\hfill
$
\\[2.75ex]
\mbox{}\hfill
$
x_1^{}\sin\biggl(\dfrac{{\rm arctg}\dfrac{\beta}{\alpha}}{\ln\!\sqrt{\alpha^2\!+\!\beta^2}}\,
\ln\!\sqrt{x_1^2\!+\!x_2^2}\biggr)\!   +
x_2\cos\biggl(\dfrac{{\rm arctg}\dfrac{\beta}{\alpha}}{\ln\sqrt{\alpha^2\!+\!\beta^2}}\,
\ln\!\sqrt{x_1^2\!+\!x_2^2}\biggr)\!\!\Biggr)\!\!
\Bigl(\!\sqrt{x_1^2\!+\!x_2^2}\,\Bigr)^{\!\!{}^{\scriptstyle {}-\tfrac{1}{\ln\sqrt{\alpha^2\!+\!\beta^2}} - 1}}
\hfill
$
\\[3ex]
\mbox{}\hfill
for all 
$
(x_1^{},x_2^{})\in\R^2\backslash \{(0,0)\},
\quad 
(0,0)\to (0,0) 
\hfill
$ 
\\[2.5ex]
conjugates the linear action 
\vspace{0.5ex}
$
\left(\!\!\! 
\begin{array}{cc} 
\alpha &\beta
\\[0.25ex]
{}-\beta &\alpha 
\end{array}\!\!\right) x
$
for all 
$x\in\R^2,\
\sqrt{\alpha^2+\beta^2}<1,\ \beta\ne 0,$ 
with the linear action $e^{{}-1}Ix$ for all $x\in\R^2;$
\vspace{0.5ex}

7) homeomorphism 
\\[1.75ex]
\mbox{}\hfill
$
(x_1^{},x_2^{})\to \ 
\Biggl(
x_1^{}\cos\biggl(\dfrac{{\rm arctg}\dfrac{\beta}{\alpha}}{\ln\sqrt{\alpha^2+\beta^2}}\, \ln\sqrt{x_1^2+x_2^2}\,\biggr) - 
x_2^{}\sin\biggl(\dfrac{{\rm arctg}\dfrac{\beta}{\alpha}}{\ln\sqrt{\alpha^2+\beta^2}}\ \ln\sqrt{x_1^2+x_2^2}\, \biggr),
\hfill
$
\\[2.75ex]
\mbox{}\hfill
$
x_1^{}\sin\biggl(\dfrac{{\rm arctg}\dfrac{\beta}{\alpha}}{\ln\!\sqrt{\alpha^2\!+\!\beta^2}}\,
\ln\!\sqrt{x_1^2\!+\!x_2^2}\biggr)\!   +
x_2\cos\biggl(\dfrac{{\rm arctg}\dfrac{\beta}{\alpha}}{\ln\sqrt{\alpha^2\!+\!\beta^2}}\,
\ln\!\sqrt{x_1^2\!+\!x_2^2}\biggr)\!\!\Biggr)\!\!
\Bigl(\!\sqrt{x_1^2\!+\!x_2^2}\,\Bigr)^{\!\!{}^{\scriptstyle \tfrac{1}{\ln\sqrt{\alpha^2+\beta^2}} - 1}}
\hfill
$
\\[3ex]
\mbox{}\hfill
for all 
$
(x_1^{},x_2^{})\in\R^2\backslash \{(0,0)\},
\quad 
(0,0)\to (0,0) 
\hfill
$ 
\\[2.5ex]
conjugates the linear action 
\vspace{0.5ex}
$
\left(\!\!\! 
\begin{array}{cc} 
\alpha &\beta
\\[0.25ex]
{}-\beta &\alpha 
\end{array}\!\!\right) x
$
for all 
$x\in\R^2,\
\sqrt{\alpha^2+\beta^2}>1,\ \beta\ne 0,$ 
with the linear action $eIx$ for all $x\in\R^2.$ \k
\vspace{0.5ex}

Let's consider now a topological associativity of Abelian linear groups $L^3$ and $L^4$ at $n> 1$ and $I\ne \{1\}.$ 
Thus everywhere we will suppose further, that matrixes $P_r^{}$ and $Q_r^{}$ for all $r\in I,$ 
are strongly hyperbolic. In this case all matrixes $P_r^{}$ (matrixes $Q_r^{})$  are reduced to real normal Jordan forms 
$J(P_r^{})$ (real normal Jordan  forms $J(Q_r^{})),$ for all $r\in I,$ by the general 
homothetic transformation and have an identical order of a disposition of real blocks of Jordan corresponding 
\vspace{0.35ex}
to real and complex eigenvalues. Further by means of replacement 
$\xi(x)=T^{{}-1}f(Sx)$ for all $x\in \R^n,$ where 
\vspace{0.35ex}
$P_r^{} =S\;\! J(P_r^{})\;\!S^{{}-1}, \ Q_r^{}=T\;\! J(Q_r^{})\;\!T^{{}-1}$ for all $r\in I,$ 
from identities (9.1) we will pass to identities 
\\[1.5ex]
\mbox{}\hfill                                                         % (9.3)
$
\xi( J(P_r^{})\;\!x)=J(Q_r^{})\;\!\xi(x)
$
\ for all 
$
x\in\R^n,
$
\ for all 
$
r\in I.
$
\hfill (9.3)
\\[1.5ex]
\indent
Therefore the topological associativity of Abelian linear groups $L^3$ and $L^4$ is equivalent to performance of identities (9.3).

Owing to a course of the proof of Theorem 9.2 
\vspace{0.35ex}
we come to conclusion, that the space $\R^n $ can be 
divided into the direct sum of $s$ co-ordinate subspaces 
\vspace{0.5ex}
$\R_k^{l_k^{}}$ dimensions $\dim\R_k^{l_k^{}} =l_k^{}$ such that
they are invariant steady or invariant unstable for contractions 
\vspace{0.35ex}
$J_k^{}(P_r^{})\;\!x^k$ for all $x^k\in \R_k^{l_k^{}}$ 
and $J_k^{}(Q_r^{})\;\!x^k$ for all $x^k\in \R_k^{l_k^{}}$ accordingly, 
\vspace{0.5ex}
on these invariant subspaces, and the given contractions have identical orientation, for all $r\in I;$ 
thus  $x=(x^1,\ldots, x^s),$ 
\\[1.5ex]
\mbox{}\hfill
$
J(P_r^{})\!=\! {\rm diag}\{ J_1^{}(P_r^{}),\ldots, J_s^{}(P_r^{})\}, \ 
J(Q_r^{})\!=\!{\rm diag}\{J_1^{}(Q_r^{}), \ldots, J_s^{}(Q_r^{})\}$ for all 
$
\displaystyle
r\!\in\! I, \ \sum\limits_{k=1}^s\! l_k^{}\!=\!n.
\hfill
$
\\[1.5ex]
\indent
Now let us consider the invariant subspace $\R_k^{l_k^{}}.$ 
\vspace{0.25ex}
For the homeomorphism-narrowing 
$\xi_k^{}\colon \R_k^{l_k^{}}\to \R_k^{l_k^{}}$ of the homeomorphism $\xi\colon\R^n\to \R^n$ 
\vspace{0.55ex}
on this invariant subspace from identities (9.3) follow the identities 
\\[1.25ex]
\mbox{}\hfill
$
\xi _k^{}(J_k^{}(P_r^{})\;\!x^k)=J_k^{}(Q_k^{})\;\!\xi _k^{}(x^k)
$ 
\ for all 
$x^k\in \R_k^{l_k^{}},$ 
\ for all $r\in I.
\hfill
$
\\[1.5ex]
\indent
Not belittling a generality, we will consider, 
\vspace{0.35ex}
that the invariant subspace $\R_k^{l_k^{}}$ is unstable for linear maps 
$J_k^{}(P_1^{})\;\! x^k$ for all $x^k\in \R_k^{l_k^{}}$ and 
\vspace{0.5ex}
$J_k^{}(Q_1^{})\;\! x^k$ for all $x^k\in \R_k^{l_k^{}}$ 
(for in case of a stability of an invariant subspace $\R_k^{l_k^{}}$ it is enough to consider linear maps 
\vspace{0.5ex}
$J_k^{{}-1}(P_1^{})\;\!x^k$ for all $x^k\in \R_k^{l_k^{}}$ and 
$J_k^{{}-1}(Q_1^{})\;\! x^k$ for all $x^k\in \R_k^{l_k^{}}).$
\vspace{0.5ex}
Using the homeomorphisms $u_k^{}\colon\R_k^{l_k^{}}\to \R_k^{l_k^{}}$ 
and $v_k^{}\colon \R_k^{l_k^{}}\to \R_k^{l_k^{}}$ of aspects 1) --- 7) from the proof of the theorem 9.2, 
\vspace{0.35ex}
we will reduce linear maps $J_k^{}(P_1^{})\;\! x^k$ for all $x^k\in \R_k^{l_k^{}}$ and 
\vspace{0.5ex}
$J_k^{}(Q_1^{})\;\! x^k$ for all $x^k\in \R_k^{l_k^{}}$ 
accordingly, in one of two canonical (see the proof of the theorem 9.2) aspects: 
\vspace{0.35ex}
$ex^k$ for all $x^k\in \R_k^{l_k^{}}$ 
(in case of positive orientation of aforementioned linear maps), or 
${\rm diag}\{e,\ldots, e,{}-e\}\;\! x^k$ for all $x^k\in \R_k^{l_k^{}} $ 
\vspace{0.35ex}
(in case of their negative orientation).
Thus, taking into account the disposition of real blocks of Jordan noted before an order, 
by direct evaluations we are convinced, that images of maps 
$J_k^{}(P_r^{})\;\! x^k$ for all $x^k\in \R_k^{l_k^{}}$ 
\vspace{0.35ex}
and $J_k^{}(Q_r^{})\;\! x^k$ for all $x^k\in \R_k^{l_k^{}}$ 
keep linearity, and matrixes thus corresponding to again received linear maps are 
real normal Jordan forms of the same structure, as real normal Jordan forms 
$J_k^{}(P_r^{})$ and $J_k^{}(Q_r^{}),$ accordingly, for all $r\in I\backslash \{1\}.$

Following auxiliary statements allow to solve a problem about a topological associativity for Abelian linear phase 
groups $L^3$ and $L^4$ of general situation.
\vspace{0.75ex}

{\bf Lemma 9.1.} 
{\it 
Let at $I =\{1,2\}$ the matrixes $P_1^{} =Q_1^{}=eI,$ and the matrixes 
\vspace{0.25ex}
$P_2^{}$ and $Q_2^{}$ are represented by real normal Jordan forms. 
\vspace{0.25ex}
Then at a topological conjunction of linear groups $L^3$ and $L^4$ the matrixes $P_2^{}$ and $Q_2^{}$ 
have identical real structure.
}
\vspace{0.5ex}

{\sl Proof.} 
Let linear groups $L^3$ and $L^4$ are topologically conjugated, i.e. identities 
\\[1.5ex]
\mbox{}\hfill                                                % (9.4)
$
\xi(e\;\!x)=e\,\xi(x)
$
\ for all 
$
x\in\R^n,
$
\hfill (9.4)
\\[1.25ex]
and 
\\[1.25ex]
\mbox{}\hfill                                                % (9.5)
$
\xi(P_2^{}\;\!x)=Q_2^{}\;\!\xi(x)
$
\ for all 
$
x\in\R^n.
$
\hfill (9.5)
\\[1.75ex]
take place. On the basis of (9.4) it is received, that 
\\[1.75ex]
\mbox{}\hfill
$
\dfrac{\xi _k^{}(ex)}{\xi _n^{}(ex)}\equiv 
\dfrac{\xi _k^{}(x)}{\xi_n^{}(x)}, 
\quad 
k=1,\ldots, n-1,
\hfill
$ 
\\[1.75ex]
and, means, are fair representations 
\\[1.75ex]
\mbox{}\hfill                                                % (9.6)
$
\dfrac{\xi _k^{}(x)}{\xi_n^{}(x)}\equiv 
g_k^{}\biggl(\;\!\dfrac{x_1^{}}{x_n^{}}\,,\ldots,\dfrac{x_{n-1}^{}}{x_n^{}}\biggr),
\quad 
k=1,\ldots, n-1,
$
\hfill (9.6)
\\[2.25ex]
as $\dfrac{x_k^{}}{x_n^{}}\,, \ k =1,\ldots, n-1,$ 
\vspace{0.5ex}
is a base of nondegenerate absolute invariants of map 
$e\;\!x$ for all $x\in\R^n.$ 
Owing to relations (9.4) --- (9.6) we come to conclusion, that the topological conjunction 
\\[1.5ex]
\mbox{}\hfill                                                % (9.7)
$
h(P_2^{}\;\! v)=Q_2^{}\;\! h(v)
$
\ for all 
$
v\in \R P^{n-1},
$
\hfill (9.7)
\\[1.75ex]
implies from a topological conjunction of linear groups $L^3$ and $L^4$, where 
\\[1.25ex]
\mbox{}\hfill
$
h(v)=\bigl(h_1^{}(v),\ldots, h_n^{}(v)\bigr).
\hfill 
$
\\[1.25ex]
\indent
Now justice of the statement of  Lemma 9.1 
\vspace{0.25ex}
implies from this the fact, that the number of fixed points of linear-fractional map 
\vspace{0.25ex}
$P_2^{}v$ for all $v\in \R P^{n-1}$ (linear-fractional map $Q_2^{}v$ for all $v\in \R P^{n-1})$ 
coincides with number of eigenvectors 
\vspace{1ex}
of the matrix $P_2^{}$ (the matrix $Q_2^{}).$\k

{\bf Lemma 9.2.} 
{\it 
Let  at  $I =\{1, 2\}$ the matrixes $P_1^{}=Q_1^{}=eI,$  and the matrixes 
\\[1.75ex]
\mbox{}\hfill
$
P_2^{}= 
{\rm diag}\left\{ 
\left(\!\! 
\begin{array}{cc}
\alpha _{12}^{} & \beta _{12}^{}
\\[0.75ex]
{}-\beta _{12}^{} & \alpha _{12}^{}
\end{array}\!\!\right),
\ldots ,
\left(\!\!
\begin{array}{cc}
\alpha _{s2}^{} & \beta _{s2}^{}
\\[0.75ex]
{}-\beta _{s2}^{} & \alpha _{s2}^{} 
\end{array}\!\!\right)\!,\,
p_{2s+1,2}^{} ,\ldots ,p_{n2}^{}\right\}
\hfill
$
\\[1.5ex]
and 
\\[1.5ex]
\mbox{}\hfill
$
Q_2^{}=
{\rm diag}\left\{
\left(\!\!
\begin{array}{cc}
\gamma_{12}^{} & \delta_{12}^{}
\\[0.75ex]
{}-\delta_{12}^{} & \gamma_{12}^{}
\end{array}\!\!\right),
\ldots,
\left(\!\!
\begin{array}{cc}
\gamma_{s2}^{} & \delta_{s2}^{}
\\[0.75ex]
{}-\delta_{s2}^{} & \gamma_{s2}^{}
\end{array}\!\!\right)\!,\, 
q_{2s+1,2}^{},\ldots,q_{n2}^{}\right\}
\hfill
$ 
\\[1.5ex]
are strongly hyperbolic, 
\\[1.75ex]
\mbox{}\hfill
$
\pi ^{{}-1}\arg(\alpha _{k2}^{} +i\beta _{k2}^{})\notin\Q, 
\ \
\pi ^{{}-1}\arg(\gamma _{k2}^{} + i\delta _{k2}^{})\notin \Q, 
\ \
\dfrac{\ln(\alpha_{k2}^{} +i\beta _{k2}^{})}{\ln(\alpha _{l2}^{} +i\beta _{l2}^{})}\notin \Q, 
\ \
\dfrac{\ln(\gamma _{k2}^{} +i\delta _{k2}^{})}{\ln(\gamma _{l2}^{} +i\delta _{l2}^{})}\notin \Q, 
\hfill
$
\\[2.25ex]
\mbox{}\hfill
$
l\ne k, 
\ \ k=1,\ldots, s, 
\ \ l=1,\ldots, s.
\hfill
$ 
\\[1.5ex]
Then for a topological conjunction of linear groups $L^3$ and $L^4$ 
\vspace{0.35ex}
it is necessary and enough ex\-is\-ten\-ce of such permutations 
$
\!\rho\colon\! (1,\ldots ,s)\!\to\! (1,\ldots, s), \, 
\chi\!\colon\! (2s\!+\!1,\ldots ,n)\!\to\! (2s\!+\!1,\ldots, n)\!\!$ that}
\\[2ex]
\mbox{}\hfill                                                          % (9.8)
$
\gamma _{\rho (k)\;\! 2}^{} =\alpha _{k2}^{},
\quad
\delta _{\rho (k)\;\! 2}^{} =\beta _{k2}^{},
\quad
k=1,\ldots, s;
\qquad
q_{\chi (k)\;\! 2}^{}=p_{k2}^{},
\quad 
k=2s+1,\ldots, n.
$
\hfill {\rm (9.8)}
\\[2ex]
\indent
{\sl Proof. The necessity.} 
As well as at the proof of the lemma 9.1, we receive relations
\linebreak
 (9.4) --- (9.7). 
On their foundation we come to conclusion, that the conjugating homeomorphism 
$L_k^{} =\bigl\{ v,\, v_\tau^{} =0,\, \tau =1,\ldots, n,\,\tau \ne 2k-1,\,\tau \ne 2k\bigr\},\ k=1,\ldots, s,$ 
\vspace{0.75ex}
of linear-fractional maps $P_2^{}\;\! v$ for all $v\in \R P^{n-1}$ and $Q_2^{}\;\! v$ for all $v\in \R P^{n-1}.$ 
\vspace{0.75ex}
It implies from this, that: 

1) $\overline{O}(v)=L_k^{}$ for all $v\in L_k^{},\ k=1,\ldots, s,$ where $O$ there is an orbit of 
a corresponding point at operations of the previous linear-fractional maps 
(it is proved on the basis of a co\-n\-se\-qu\-en\-ce from the theorem of Kronecker [31] similarly course 
of the proof of Theorem 6.12); 

2) for remaining possible arcwise connected one-dimensional invariant sets $I^1,$ 
homeomorphic to projective straight lines, property $\overline{O}(v)=I^1$ for all $v\in I^1$ 
is not fulfilled (proved by reviewing of operations of the previous linear-fractional maps on bases 
of nondegenerate absolute invariants
\\[1.75ex]
\mbox{}\hfill
$
\left( 
\dfrac{v_{2k-1}^2 +v_{2k}^2}{v_{2s-1}^2+v_{2s}^2}
\right)^{{}-{\rm arctg}\;\!\tfrac{\beta _{s2}^{}}{\alpha _{s2}^{}}}\,
\exp\left( \ln\dfrac{\alpha _{k2}^2+ \beta _{k2}^2}{\alpha _{s2}^2 + \beta _{s2}^2}\ 
{\rm arctg}\;\!\dfrac{v_{2s}^{}}{v_{2s-1}^{}}\right), 
\hfill
$
\\[2.25ex]
\mbox{}\hfill
$
- \ 
{\rm arctg}\;\!\dfrac{\beta_{s\;\!2}^{}}{\alpha _{s\;\!2}^{}}\  
{\rm arctg}\;\!\dfrac{v_{2k}^{}}{v_{2k-1}^{}} 
\ + \
{\rm arctg}\;\!\dfrac{\beta _{k2}^{}}{\alpha _{k2}^{}}\
{\rm arctg}\;\!\dfrac{v_{2s}^{}}{v_{2s-1}^{}}, 
\quad 
k=1,\ldots, s;
\hfill
$
\\[2.25ex]
\mbox{}\hfill
$
\left(\dfrac{v_{2s-1}^2 + v_{2s}^2}{v_{n-1}^2}\right)^{{}-{\rm arctg}\;\!\tfrac{\beta _{s2}^{}}{\alpha _{s2}^{}}}\,
\exp\left(\ln\dfrac{\alpha _{s2}^2+\beta _{s2}^2}{p_{n-1,2}^2}\ 
{\rm arctg}\;\!\dfrac{v_{2s}^{}}{v_{2s-1}^{}}\right); 
\hfill
$
\\[2.25ex]
\mbox{}\hfill
$
\left(\dfrac{v_{2s-1}^2 +v_{2s}^2}{v_n^2}\right)^{{}-{\rm arctg}\tfrac{\beta _{s2}^{}}{\alpha _{s2}^{}}}\,
\exp\left( \ln\dfrac{\alpha _{s2}^2 +\beta _{s2}^2}{p_{n2}^2}\ 
{\rm arctg}\dfrac{v_{2s}^{}}{v_{2s-1}^{}}\right); 
\hfill
$
\\[2.5ex]
\mbox{}\hfill
$
\bigl|v_k^{}\bigr|^{{}^{\scriptstyle \ln\tfrac{| p_{n-1,2}^{}|}{| p_{n2}^{}|}}}\,
\bigl| v_{n-1}^{}\bigr|^{{}^{\scriptstyle \ln\tfrac{| p_{n2}^{}|}{| p_{k2}^{}|}}} \,
\bigr| v_n^{}\bigr|^{{}^{\scriptstyle \ln\tfrac{| p_{k2}^{}|}{| p_{n-1,2}^{}|}}} ,
\quad
k=2s+1,\ldots, n-2;
\hfill
$ 
\\[1.75ex]
and
\\[1.75ex]
\mbox{}\hfill
$
\left( 
\dfrac{v_{2k-1}^2 +v_{2k}^2}{v_{2s-1}^2+v_{2s}^2}
\right)^{{}-{\rm arctg}\;\!\tfrac{\delta _{s2}^{}}{\gamma _{s2}^{}}}\,
\exp\left( \ln\dfrac{\gamma _{k2}^2+ \delta _{k2}^2}{\gamma _{s2}^2 + \delta _{s2}^2}\ 
{\rm arctg}\;\!\dfrac{v_{2s}^{}}{v_{2s-1}^{}}\right), 
\hfill
$
\\[2.5ex]
\mbox{}\hfill
$
- \ 
{\rm arctg}\;\!\dfrac{\delta_{s\;\!2}^{}}{\gamma _{s\;\!2}^{}}\  
{\rm arctg}\;\!\dfrac{v_{2k}^{}}{v_{2k-1}^{}} 
\ + \
{\rm arctg}\;\!\dfrac{\delta _{k2}^{}}{\gamma _{k2}^{}}\
{\rm arctg}\;\!\dfrac{v_{2s}^{}}{v_{2s-1}^{}}, 
\quad 
k=1,\ldots, s;
\hfill
$
\\[2.25ex]
\mbox{}\hfill
$
\left(\dfrac{v_{2s-1}^2 + v_{2s}^2}{v_{n-1}^2}\right)^{{}-{\rm arctg}\;\!\tfrac{\delta _{s2}^{}}{\gamma _{s2}^{}}}\,
\exp\left( \ln\dfrac{\gamma _{s2}^2+\delta _{s2}^2}{q_{n-1,2}^2}\ 
{\rm arctg}\;\!\dfrac{v_{2s}^{}}{v_{2s-1}^{}}\right); 
\hfill
$
\\[2.25ex]
\mbox{}\hfill
$
\left(\dfrac{v_{2s-1}^2 +v_{2s}^2}{v_n^2}\right)^{{}-{\rm arctg}\tfrac{\delta _{s2}^{}}{\gamma _{s2}^{}}}\,
\exp\left( \ln\dfrac{\gamma _{s2}^2 +\delta _{s2}^2}{q_{n2}^2}\ 
{\rm arctg}\dfrac{v_{2s}^{}}{v_{2s-1}^{}}\right); 
\hfill
$
\\[2.5ex]
\mbox{}\hfill
$
\bigl|v_k^{}\bigr|^{{}^{\scriptstyle \ln\tfrac{| q_{n-1,2}^{}|}{| q_{n2}^{}|}}}\,
\bigl| v_{n-1}^{}\bigr|^{{}^{\scriptstyle \ln\tfrac{| q_{n2}^{}|}{| q_{k2}^{}|}}}\, 
\bigr| v_n^{}\bigr|^{{}^{\scriptstyle \ln\tfrac{| q_{k2}^{}|}{| q_{n-1,2}^{}|}}} ,
\quad
k=2s+1,\ldots, n-2,
\hfill
$ 
\\[2ex]
of zero degree of a homogeneity corresponding to them).

Introducing auxiliary variables $z_k^{}=v_{2k-1}^{} +iv_{2k}^{}, \ k=1,\ldots, s,$ 
\vspace{0.25ex}
on the basis of Theorem~3.1 we come to the first part of relations (9.8). 

Besides, the homeomorphism $h$ translates each other points
\\[1ex]
\mbox{}\hfill
$
O_k^{} =\{v,\ v_\tau^{}=0,\ \tau=1,\ldots, n,\  \tau \ne k\}, 
\quad 
k=2s+1,\ldots, n
\hfill
$ 
\\[1.25ex]
(being fixed for linear-fractional maps $Q_2^{}\;\!v$ for all $v\in \R P^{n-1}).$ 
\vspace{0.25ex}
Now by similar reasonings on the basis of (9.5) --- (9.7) and Theorem 9.1 it is received the second part of relations (9.8).

{\sl The sufficiency} 
\vspace{0.25ex}
is proved by direct evaluations application of conjugating homeomorphism 
$\xi(x)=\bigl(x_{2\rho(1)-1}^{},x_{2\rho(1)}^{},\cdots, x_{2\rho(s)-1}^{}, x_{2\rho(s)}^{},
x_{\chi(2s+1)}^{},\ldots ,x_{\chi(n)}^{}\bigr)$ for all $x\in  \R^n.$
\vspace{0.75ex}

Following statements it is received on the basis of two previous, Theorem 9.1 and that fact, that for linear map 
${\rm diag}\,\{e, \ldots, e,-e\}\;\! x$ for all $x\in \R^n $ dimension of a maximum invariant 
subspace on which the contraction of this map has positive orientation, is an invariant at a topological conjunction. \k
\vspace{0.75ex}

{\bf Lemma 9.3.} 
{\it 
Let at $I =\{ 1,2\} $ the matrixes $P_1^{} =Q_1^{} ={\rm diag}\,\{e,\ldots ,e,{}-e\},$ 
and the matrixes $P_2^{}$ and $Q_2^{}$ are represented by real normal Jordan forms. 
Then at a topological conjunction of linear groups $L^3$ and $L^4$ the matrixes $P_2^{}$ and $Q_2^{}$ 
have identical real structure.
}

{\bf Lemma 9.4.} 
{\it 
Let at $I =\{1,2\}$ the matrixes 
$P_1^{} =Q_1^{} ={\rm diag}\,\{e,\ldots ,e,{}-e\},$ and  the matrixes 
\\[1.5ex]
\mbox{}\hfill
$
P_2^{}= {\rm diag}\,\left\{ 
\left(\!\! 
\begin{array}{rc}
\alpha _{12}^{} & \beta _{12}^{}
\\[0.5ex]
{}-\beta _{12}^{} & \alpha _{12}^{}
\end{array}\!\!\right),
\ldots ,
\left(\!\!\begin{array}{rc}
\alpha _{s\;\!2}^{} & \beta _{s\;\!2}^{}
\\[0.5ex]
{}-\beta _{s\;\!2}^{} & \alpha _{s\;\!2}^{}
\end{array}\!\!\right)\!,\,
p_{2s+1,2}^{} ,\ldots ,p_{n2}^{}\right\} 
\hfill
$ 
\\[2ex] 
and 
\\[1.5ex]
\mbox{}\hfill
$
Q_2^{}=
{\rm diag}\,\left\{ 
\left(\!\! 
\begin{array}{rc}
\gamma _{12}^{} & \delta _{12}^{}
\\[0.5ex]
{}-\delta _{12}^{} & \gamma _{12}^{}
\end{array}\!\!\right),
\ldots ,
\left(\!\!\begin{array}{rc}
\gamma _{s\;\!2}^{} & \delta _{s\;\!2}^{}
\\[0.5ex]
{}-\delta _{s\;\!2}^{} & \gamma _{s\;\!2}^{}
\end{array}\!\!\right)\!,\,
q_{2s+1,2}^{} ,\ldots ,q_{n2}^{}\right\} 
\hfill
$ 
\\[2ex] 
are  strongly hyperbolic, 
\\[1.5ex]
\mbox{}\hfill
$
\pi ^{{}-1}\arg (\alpha _{k2}^{} +i\beta _{k2}^{})\not\in \Q, 
\qquad
\pi ^{{}-1}\arg ( \gamma _{k2}^{} +  i\delta _{k2}^{})\not\in \Q, 
\hfill
$
\\[3ex]
\mbox{}\hfill
$
\dfrac{\ln(\alpha _{k2}^{} +i\beta _{k2}^{})}{\ln(\alpha _{l2}^{} +i\beta _{l2})}\not\in\Q, 
\quad
\dfrac{\ln(\gamma _{k2}^{} +i\delta _{k2}^{})}{\ln(\gamma _{l2}^{} +i\delta _{l2}^{})}\not\in \Q, 
\quad
k=1,\ldots, s, \ l=1,\ldots, s,
\ \ 
l\ne k.
\hfill
$ 
\\[2ex]
Then for a topological conjunction of linear groups $L^3$ and $L^4$ it is necessary and enough existence of 
such permutations 
$
\rho\colon (1,\ldots ,s)\to (1,\ldots ,s),\
\chi\colon (2s+1,\ldots ,n)\to (2s+1,\ldots ,n)
$ 
that relations {\rm (9.8)} are fulfilled.
}
\vspace{0.5ex}

In case of smooth and holomorphic conjunctions of linear groups $L^3$ and $L^4$ the following statement takes place.
\vspace{0.5ex}

{\bf Theorem 9.3.} 
{\it 
Linear groups $L^3$ and $L^4$ are smoothly {\rm(}holomorphic{\rm)} conjugated in only case when they are linearly conjugated.
}
\vspace{0.35ex}

{\sl Proof. The necessity.} 
Let identities (9.1) take place. Calculating in them Jacobi matrix in the point $x=0,$ 
we have matrix equalities ${\sf D}_x^{} f(0)P_r^{}=Q_r^{}\;\!{\sf D}_x^{} f(0)$ for all $r\in I$. 
\vspace{0.35ex}

Therefore linear map ${\rm D}_x^{}f (0)x$ for all $x\in\R^n$ satisfies to identities (9.1). 
As $f $ is a diffeomorphism (holomorphism), it is nondegenerate.

{\sl The sufficiency} is checked by direct evaluations. \k
\\[3.75ex]
\centerline{
\large\bf  
10. Applications to real nonautonomous linear differential systems
}
\\[1.5ex]
\indent
Theorems 9.1 --- 9.3, and also the algorithm based on Lemmas 9.1 --- 9.4, allows on the basis of Theorems 1.1 --- 1.3, 
to spend topological, smooth and holomorphic classifications of real nonautonomous linear differential systems of the aspect (8.1). 

Besides, on the basis of the received algorithm it is possible to draw a conclusion, that real completely solvable 
(i.e. at two and more independent variables) linear differential systems with periodic coefficients are structurally unstable 
(the commutativity of phase groups of the given class of differential systems implies from a commutativity of fundamental group of a many-dimensional torus). 

We will notice, that the theorem 9.2 on the basis of the theorem 1.1 gives criterion of a structural 
stability of real linear ordinary differential systems with periodic coefficients.
\\[3.75ex]
\centerline{
\large\bf  
11. Phase groups of covering foliations,
}
\\[0.35ex]
\centerline{
\large\bf  
defined by real nonautonomous Riccati equations 
}
\\[1.5ex]
\indent
We will consider Riccati equations
\\[1.5ex]
\mbox{}\hfill                                                 % (11.1)
$
\displaystyle
dx=\sum\limits_{j=1}^m
\bigl( a_{2j}^{}(t_1^{},\ldots,t_m^{})\;\!x^2+
a_{1j}^{}(t_1^{},\ldots,t_m^{})\;\! x+a_{0j}^{}(t_1^{},\ldots,t_m^{})\bigr)\;\! dt_j^{} 
$
\hfill (11.1)
\\[1.25ex]
and
\\[1.5ex]
\mbox{}\hfill                                                 % (11.2)
$
\displaystyle
dx=\sum\limits_{j=1}^m
\bigl(b_{2j}^{}(t_1^{},\ldots,t_m^{})\;\! x^2+
b_{1j}^{}(t_1^{},\ldots,t_m^{})\;\! x+b_{0j}^{}(t_1^{},\ldots,t_m^{})\bigr)\;\!dt_j^{}
$
\hfill (11.2)
\\[1.5ex]
ordinary at $m\!=\!1\!$ and completely solvable at $\!m> 1,\!$ 
\vspace{0.25ex}
where holomorphic functions $\!a_{ij}^{}\colon\! A\!\to\! {\mathbb R}$ and 
$b_{ij}^{}\colon B\to {\Bbb R},\ i =0,1, 2,\ j =1,\ldots, m,$ 
\vspace{0.25ex}
path connected holomorphic varieties $A$ and $B $ are holomorphic equivalent each other, 
fundamental groups $\pi_1^{}(A)$ and $\pi_1^{}(B)$ have final number $\nu\in\N$ of the forming.

The common solutions of Riccati equations (11.1) and (11.2) define covering foliations 
\vspace{0.25ex}
${\frak P} ^7$ and $ {\frak P}^8,$ 
accordingly, on varieties $ \overline{\R} \times A$ and $\overline{\R} \times B,$ 
where $ \overline{\R}$ is a real straight line $\R,$ supplemented by the point at infinity $\infty.$

We will say, that real Riccati equations (11.1) and (11.2) are {\it topologically {\rm(}smoothly,  holomorphically{\rm)} equivalent} 
\vspace{0.25ex}
if exists a homeomorphism (a diffeomorphism, a holo\-morphism) $h\colon\overline {\R} \times A\to\overline {\R}\times B,$ 
\vspace{0.25ex}
translating the layers of the covering foliation ${\frak P}^7$ in the layers of the covering foliation ${\frak P}^8.$ 
Similarly we enter concepts of {\it embedding {\rm(}smooth embedding,  holomorphic embedding}) and 
\vspace{0.25ex}
{\it covering {\rm(}smoothly covering,  holomorphically covering}) of Riccati equations.

The phase group $Ph ({\frak P}^7) $ of the covering foliation ${\frak P}^7$ is generated by the forming nondegenerate 
linear-fractional transformations 
\\[2ex]
\mbox{}\hfill                                           % (11.3)
$
P_r^{}(x)=
\dfrac{p_{1r}^{}x+p_{2r}^{}}{p_{3r}^{}x+p_{4r}^{}}
$
\ \ for all
$x\in\overline{\R},
\quad 
r=1,\ldots, \nu,
$
\hfill (11.3)
\\[1.5ex]
to which we will put in accordance nondegenerate matrices 
\vspace{0.5ex}
$
P_r^{}=
\left(\!\!
\begin{array}{cc}
p_{1r}^{} & p_{2r}^{}
\\[0.35ex]
p_{3r}^{} & p_{4r}^{}
\end{array}\!\!\right),
\ r=1,\ldots, \nu;
$ 
phase group $Ph ({\frak P}^8)$ of the covering foliation ${\frak P}^8$ 
\vspace{0.35ex}
is generated by the forming nondegenerate linear-fractional transformations 
\\[2ex]
\mbox{}\hfill                                           % (11.4)
$
Q_r^{}(x)=
\dfrac{q_{1r}^{}x+q_{2r}^{}}{q_{3r}^{}x+q_{4r}^{}}
$
\ \ for all
$x\in\overline{\R},
\quad 
r=1,\ldots, \nu,
$
\hfill (11.4)
\\[1.5ex]
to which we will put in accordance nondegenerate matrices 
$
Q_r^{}=
\left(\!\!
\begin{array}{cc}
q_{1r}^{} & q_{2r}^{}
\\[0.35ex]
q_{3r}^{} & q_{4r}^{}
\end{array}\!\!\right),
\ r=1,\ldots, \nu.
$
\\[3.75ex]
\centerline{
\large\bf  
12. Conjunctions of linear-fractional actions on $\overline{\R} $
}
\\[1.5ex]
\indent
Consider a problem about a finding of necessary and sufficient conditions 
\vspace{0.25ex}
of existence such homeomorphism (diffeomorphism, holo\-morphism) 
$f\colon\overline{\R} \to\overline{\R}$ that identities
\\[1.75ex]
\mbox{}\hfill                                           % (12.1)
$
f(P_r^{}(x))=Q_r^{}(f(x))
$ 
\ for all 
$
x\in\overline{\R},
$
\ for all 
$
r\in I, 
$
\hfill (12.1)
\\[1.75ex]
take place, where square matrices $P_r^{}\in GL(2, \R),\ Q_r^{}\in GL(2,\R)$ for all $r\in I.$ 
\vspace{0.5ex}

Group of linear-fractional actions on $\overline{\R},$ 
\vspace{0.25ex}
formed by matrices $P_r^{}$ for all $r\in I$ 
we will designate through $PL^3,$ and through $PL^4$ we will designate the similar group formed by matrices 
$Q_r^{}$ for all $r\in I.$

Consider at first a case of the Abelian real linear-fractional phase groups.

Let's prove some auxiliary statements on which basis we will receive criteria of topological, 
smooth and holomorphic conjunctions of Abelian real linear-fractional phase groups.
\vspace{0.35ex}

{\bf Lemma 12.1.} 
{\it 
Let $f$ is a homeomorphism, cojugating linear-fractional phase groups $PL^3$ and $PL^4.$
Then real normal Jordan forms of all matrices $P_r^{}$ and $Q_r^{},$ 
defining by nonidentical linear-fractional transformations, have the same structure, $r =1,\ldots,\nu.$
}
\vspace{0.25ex}

{\sl The proof} of the given statement is spent on the basis of that fact, that the quantity of fixed points 
of nonidentical linear-fractional maps (11.3) (maps (11.4)) coincides with number of eigenvectors of corresponding matrices 
$P_r^{}$ (matrices $Q_r^{})$ for all $r\in I.$ \k
\vspace{0.5ex}

{\bf Lemma 12.2.} 
{\it 
For a topological conjugation of Abelian linear-fractional groups $PL^3$ and $PL^4$ 
it is necessary, that real normal Jordan forms of all matrices $P_r^{}$ and $Q_r^{}$ for all $r\in I$ 
defining nonidentical linear-fractional transformations,  had the same structure.
}
\vspace{0.25ex}

{\sl The proof} of the given statement is spent on the basis of Lemma 12.1 and that fact that material normal Jordan forms 
of all permutable among themselves the matrices of the second order 
defining nonidentical linear-fractional transformations, 
\vspace{0.75ex}
have the same structure. \k

{\bf Theorem 12.1.} 
{\it 
Let matrices 
\vspace{0.5ex}
$
P_r^{}=S\;\!{\rm diag}\{p_{1r}^{},p_{2r}^{}\}\;\! S^{{}-1},
\ 
Q_r^{}=T\;\! {\rm diag}\{q_{1r}^{},q_{2r}^{}\}\;\!T^{{}-1}$ for all $r\in I.$
Then for a topological conjugation of linear-fractional groups $PL^3$ and $PL^4$ it is necessary and enough, that}
\\[1.75ex]
\mbox{}\hfill                                                     % (12.2)
$
\dfrac{q_{1r}^{}}{q_{2r}^{}}\, =\, 
\dfrac{p_{1r}^{}}{p_{2r}^{}}\, \left|\dfrac{p_{1r}^{}}{p_{2r}^{}}\right|^{\alpha}
$
\, \ for all 
$
r\in I,
\quad
\alpha\ne {}-1.
$
\hfill (12.2)
\\[2.25ex]
\indent
{\sl Proof.} 
\vspace{0.35ex}
Since groups $PL^3$ and $PL^4$ are topological conjugated, then identities (12.1) take place. 
With the help of replacement $\xi\colon x\to T^{{}-1}\circ f\circ S(x)$ for all $x\in\overline{\R},$ 
from identities (12.1) we pass to identities 
\\[1.75ex]
\mbox{}\hfill                                                     % (12.3)
$
\xi\bigl({\rm diag}\{p_{1r}^{},p_{2r}^{}\}(x)\bigr)=
{\rm diag}\{q_{1r}^{},q_{2r}^{}\}(\xi(x))
$
\ for all 
$
x\in\overline{\R},
$
\ for all 
$
r\in I.
$ 
\hfill (12.3) 
\\[1.75ex]
\indent
Hence, the topological conjunction of groups $PL^3$ and $PL^4$ is equivalent to performance of identities (12.3).
\vspace{0.25ex}

{\sl The necessity.} Let identities (12.3) take place.
\vspace{0.5ex}

If all $\dfrac{p_{1r}^{}}{p_{2r}^{}}=1$ 
\vspace{0.75ex}
for all $r\in I,$ then from (12.3) it is received, that 
$\dfrac{q_{1r}^{}}{q_{2r}^{}}=1$ for all $r\in I.$ Therefore in this case relations (12.2) are carried out at $\alpha=0.$
\vspace{0.75ex}

Let now $\dfrac{p_{1r}^{}}{p_{2r}^{}}\ne 1,\ r\in\{1,\ldots,\nu\}.$ 
\vspace{0.75ex}
Then on the basis of identities (12.3) we come to conclusion, that either $\psi (0) =0,$ or $\psi (0) = \infty.$
\vspace{0.35ex}

Let $\xi (0) =0.$ Directly from the theorem 9.1 we have relations (12.2) at $\alpha>{}-1.$
\vspace{0.35ex}

Let now $\xi (0) = \infty.$ Then by means of replacement $\zeta=\dfrac{1}{\xi}$ 
we come to the previous case and as a result we receive relations (12.2) at $\alpha <{}-1.$
\vspace{0.35ex}

{\sl The sufficiency} is proved by construction of conjugating homeomorphism 
\\[1.25ex]
\mbox{}\hfill
$
\xi\colon x\to\,  \gamma\;\! x\;\! |x|^{\alpha}
$ 
\ for all 
$
x\in\overline{\R}.\ \, \k
\hfill
$
\\[1.75ex]
\indent
{\bf Theorem 12.2.} 
{\it 
Let matrices 
\\[1.25ex]
\mbox{}\hfill
$
P_r^{}= S\;\! 
\left(\!\!
\begin{array}{cr}
\cos\alpha_r^{} & {}-\sin\alpha_r^{}
\\[0.35ex]
\sin\alpha_r^{} & \cos\alpha_r^{}
\end{array}\!\!\right)\;\! S^{{}-1}, 
\quad 
\alpha_r^{}\in ({}-\pi,\pi], 
\hfill
$
\\[1.75ex]
\mbox{}\hfill
$
Q_r^{}=T\;\! 
\left(\!\!
\begin{array}{cr}
\cos\beta_r^{} & {}-\sin\beta_r^{}
\\[0.35ex]
\sin\beta_r^{} & \cos\beta_r^{}
\end{array}\!\!\right)\;\! T^{{}-1},
\quad  \beta_r^{}\in ({}-\pi,\pi], 
$
\ for all $r\in I.
\hfill
$
\\[1.5ex]
Then for topological, smooth and holomorphic conjunctions of linear-fractional groups 
$PL^3$ and $PL^4$ it is necessary and enough, that either 
\\[1ex]
\mbox{}\hfill                                  % (12.4)
$
\beta_r^{}=\alpha_r^{}
$
\ for all 
$
r\in I, 
$
\hfill {\rm (12.4)}
\\[0.1ex]
or
\\[0.1ex]
\mbox{}\hfill                                  % (12.5)
$
\beta_r^{}={}-\alpha_r^{}
$
\ for all 
$
r\in I. 
$
\hfill {\rm (12.5)}
\\[1.5ex]
}
\indent
{\sl Proof.} 
Similarly, as well as at the proof of the theorem 12.1, we come to conclusion, that the topological conjunctions 
of groups $PL^3$ and $PL^4$ is equivalent to performance of identities 
\\[1.75ex]
\mbox{}\hfill
$
\xi\biggl(\dfrac{x\cos\alpha_r^{}-\sin\alpha_r^{}}{x\sin\alpha_r^{}+\cos\alpha_r^{}}\biggr)= 
\dfrac{\xi(x)\cos\beta_r^{} -\sin\beta_r^{}}{\xi(x)\sin\beta_r^{} +\cos\beta_r^{}}
$
\ for all 
$
x\in\overline{\R},
$
\ for all $r\in I.
\hfill
$ 
\\[1.75ex]
\indent
We will establish a homeomorphism between expanded real line $\overline{\R}$ and an unit circle $S^1$ 
by means of map 
\\[1.5ex]
\mbox{}\hfill
$
\zeta\colon x\to\  2\, {\rm arctg}\, x
$ 
\ for all $x\in\overline{\R},
\quad 
{\rm arctg}\;\! x\in \Bigl({}-\dfrac{\pi}{2}\,;\dfrac{\pi}{2}\Bigr].
\hfill
$ 
\\[1.5ex]
\indent
Now with the help of replacement 
\\[1.5ex]
\mbox{}\hfill
$
\varphi\colon t\to \zeta\circ\xi\circ\zeta^{{}-1}(t)
$ 
\ for all $t\in ({}-\pi,\pi],
\hfill
$ 
\\[1.5ex]
from last identities we come to equivalent identities 
\\[1.5ex]
\mbox{}\hfill
$
\varphi\bigl((t-2\alpha_r^{})({\rm mod}\;\! 2\pi)\bigr)= 
(\varphi(t)-2\beta_r^{})({\rm mod}\;\! 2\pi)
$ 
\ for all $t\in ({}-\pi,\pi],$ 
\ for all $r\in I.
\hfill
$
\\[1.5ex]
\indent
{\sl The necessity.} 
On the basis of that fact, that circle rotations on angles $\alpha$ and $\beta,$ 
where ${}-\pi<\alpha\leq\pi,\ {}-\pi<\beta\leq\pi,$ are topological conjugated by: 

1) positive oriented homeomorphism, if and only if $\alpha =\beta;$ 

2) oriented homeomorphism, if and only if $ \alpha ={}-\beta;$
\\ 
we come either to conditions (12.4), or to conditions (12.5), accordingly.

{\sl The sufficiency} is established by application of conjugating identical map at performance of conditions (12.4) 
and by application of conjugating holomorphism 
\\[1.5ex]
\mbox{}\hfill
$
\varphi\colon t ({\rm mod}\;\! 2\pi)\to ({}-t)({\rm mod}\;\! 2\pi)$ for all $t\in ({}-\pi,\pi]
\hfill
$ 
\\[1.5ex]
in case of performance of conditions (12.5). \k
\vspace{1ex}

{\bf Theorem 12.3.} 
{\it 
Let matrices 
$
P_r^{}\!=\!
S\;\! 
\left(\!\!\!
\begin{array}{cc}
1 & p_r^{}
\\[0.25ex]
0 & 1
\end{array}\!\!\!\right)\;\! S^{{}-1},
\ 
Q_r^{}\!=\!
T\;\! 
\left(\!\!\!
\begin{array}{cc}
1 & q_r^{}
\\[0.25ex]
0 & 1
\end{array}\!\!\!\right)\;\! T^{{}-1}$ for all $r\in I.$
Then for topological, 
\vspace{0.35ex}
smooth and holomorphic conjunctions of linear-fractional groups 
$PL^3$ and $PL^4$ it is necessary and enough, that}
\\[1.5ex]
\mbox{}\hfill                                 % (12.6)
$
q_r^{}=\lambda p_r^{}
$
\ \ for all 
$
r\in I,
\quad 
\lambda\ne 0.
$
\hfill (12.6)
\\[2ex]
\indent
{\sl Proof.} As well as earlier, we receive, that the topological conjunction of groups 
$PL^3$ and $PL^4$ is equivalent to performance of identities
\\[1.5ex]
\mbox{}\hfill                                 % (12.7)
$
\xi(x+p_r^{})=\xi(x)+q_r^{}
$ 
\ for all 
$
x\in\overline{\R},
$
\ for all 
$
r\in I. 
$
\hfill (12.7)
\\[1.5ex]
\indent
{\sl The necessity.} Let identities (12.7) take place.

If all $p_r^{}=0$ for all $r\in I,$ from (12.7) it is received that $q_r^{}=0$ for all $r\in I.$ 
\vspace{0.25ex}
Therefore in this case the relations (12.6) are carried out.
\vspace{0.35ex}

If $p_r^{}=0$ for all $r\in I,\ r\ne k,\ p_k^{}\ne 0,$ then from (12.7) it is had, that 
\vspace{0.5ex}
$q_r^{}=0$ for all $r\in I,\ r\ne k,\ q_k^{}\ne 0,$ and in quality of $\lambda$ it is possible to take number $\dfrac{q_r^{}}{p_r^{}}\,.$
\vspace{0.65ex}

Let now $p_{r_1^{}}^{}p_{r_2^{}}^{}\ne 0,\ r_1^{}\in I,\ r_2^{}\in I.$ 
From (12.7) it is received, as $q_{r_1^{}}^{}q_{r_2^{}}^{}\ne 0.$
\vspace{0.75ex}

Consider at first a case, when 
\vspace{0.5ex}
$\dfrac{p_{r_2^{}}^{}}{p_{r_1^{}}^{}}=\lambda=\dfrac{l}{n}\in\Q,\ l\in\Z, \ n\in\Z.$ 
Then owing to identities (12.7) at $r=r_1^{}$ it is had, that 
\vspace{0.5ex}
$\xi(x+lp_{r_1^{}}^{})=\xi(x)+lq_{r_1^{}}^{}$ for all $x\in\overline{\R},$ 
and at $r=r_2^{}$ it is received, that 
\vspace{0.5ex}
$\xi(x+lp_{r_1^{}}^{})=\xi(x+np_{r_2^{}}^{})=\xi(x)+nq_{r_2^{}}^{}$ for all $x\in\overline{\R}.$ 
Comparing the right parts of last expressions, we receive relations (12.6) at $r=r_1^{}$ and $r=r_2^{}.$
\vspace{0.5ex}

Let $\dfrac{p_{r_2^{}}^{}}{p_{r_1^{}}^{}}=\lambda\ne\Q.$ 
Owing to density of set of rational numbers in set of real numbers there are such sequences 
$\{l_s^{}\}$ and $\{n_s^{}\}$ of integers, that 
\\[1.25ex]
\mbox{}\hfill
$
|l_s^{}|\to {}+\infty,
\quad 
|n_s^{}|\to {}+\infty,
\quad 
\dfrac{l_s^{}}{n_s^{}}\to\lambda
$ 
\ at \ 
$
s\to {}+\infty.
\hfill
$ 
\\[1.25ex]
\indent
From identities (12.7) at $r=r_1^{}$ it is received, that 
\vspace{0.5ex}
$\xi(x+l_s^{}p_{r_1^{}}^{})-\xi(x)=l_s^{}q_{r_1^{}}^{}\!\!$ for all $\!x\!\in\!\overline{\R},$ 
and at $r=r_2^{}$ it is received, that 
\vspace{0.5ex}
$\xi(x+n_s^{}\lambda p_{r_1^{}}^{})-\xi(x)=n_s^{}q_{r_2^{}}^{}$ for all $x\in\overline{\R}.$ 
Pro\-ce\-ed\-ing from two last expressions, we receive the following chain of relations:
\\[2ex]
\mbox{}\hfill
$
1=\lim\limits_{s\to {}+\infty} \,
\dfrac{\xi(x+l_s^{}\;\!p_{r_1^{}}^{})-\xi(x)}{\xi(x+n_s^{}\;\!\lambda\;\! p_{r_1^{}}^{})-\xi(x)}=
\lim\limits_{s\to {}+\infty}\, l_s^{}\;\! q_{r_1^{}}^{}n_s^{{}-1}\;\!q_{r_2^{}}^{{}-1}= 
\lambda\;\! q_{r_1^{}}^{}q_{r_2^{}}^{{}-1}.
\hfill
$ 
\\[2ex]
From  here follows justice of relations (12.6) at $r=r_1^{}$ and $r=r_2^{}.$
\vspace{0.5ex}

Justice of relations (12.6) is similarly proved and at 
\vspace{0.5ex}
$r=r_3^{},\ p_{r_3^{}}^{}\ne 0,\ r_3^{}\ne r_2^{},\ 
r_3^{}\ne r_1^{}.$ If $p_{r_3^{}}^{}=0,$ then from (12.7) it is had $q_{r_3^{}}^{}=0$ 
and relations (12.6) at $r=r_3^{}$ take place.
\vspace{0.5ex}

{\sl The sufficiency} is established by application of conjugating holomorphism 
\\[1.25ex]
\mbox{}\hfill
$
\xi(x)=q_{r_1^{}}^{}p_{r_1^{}}^{{}-1}x
$ 
\ for all 
$
x\in\overline{\R}
\hfill
$ 
\\[1.5ex]
in case there is such index $r_1^{}\in I,$ that $p_{r_1^{}}^{}\ne 0,$ and identical map otherwise. \k
\vspace{0.5ex}

On the basis of the theorem 9.4 it is similar to the theorem 4.1 we receive the statement.
\vspace{0.5ex}

{\bf Theorem 12.4.} 
{\it 
Let the conditions of the theorem {\rm 12.1} are satisfied.
Then for smooth and holomorphic conjunctions of linear-fractional groups $PL^3$ and $PL^4$ it is necessary and enough, that} 
\\[1ex]
\mbox{}\hfill
$
\dfrac{q_{1r}^{}}{q_{2r}^{}}=\left(\dfrac{p_{1r}^{}}{p_{2r}^{}}\right)^{\!\varepsilon},
\quad
r=1,\ldots, \nu,
\quad 
\varepsilon^2=1.
\hfill
$
\\[2ex]
\indent
Consider now a case of non-Abelian linear-fractional groups $PL^3$ and $PL^4.$
\vspace{0.5ex}

{\bf Theorem 12.5.} 
{\it 
From a topological conjunction of non-Abelian linear-fractional groups $PL^3$ and $PL^4$ 
of general situation follows them holomorphic conjuction which is carried out 
by nondegenerate linear-fractional transformation.
}
\vspace{0.25ex}

{\sl Proof} of the given theorem is spent on the basis of following auxiliary statements from which 
Theorems 12.6 and 12.7 define constructive criteria of topological, smooth and holomorphic conjunc\-tions. \k
\vspace{0.5ex}

{\bf Theorem 12.6.} 
{\it 
Let conjugating of nondegenerate linear-fractional transformations {\rm (11.3)} and {\rm (11.4)} 
homeomorphism $f\colon\overline{\R}\to\overline{\R}$ is such that{\rm:}

{\rm 1)} the matrix $P_1^{}$ has pair of  complex conjugated roots 
$\cos\alpha\pm i\sin\alpha,\ \alpha\in \Bigl({}-\dfrac{\pi}{2}\,;\dfrac{\pi}{2}\Bigr];$

{\rm 2)} $\dfrac{\alpha}{\pi}\not\in\Q.$
\\[1ex]
Then this homeomorphism represents nondegenerate linear-fractional transfor\-mation.
}
\vspace{0.35ex}

{\sl Proof.} 
Let identities (12.1) are carried out.
Present the matrix $P_1^{}$ in a kind 
\\[1.5ex]
\mbox{}\hfill
$
P_1^{}=S\;\! 
\left(\!\!
\begin{array}{cr}
\cos\alpha & {}-\sin\alpha
\\[0.25ex]
\sin\alpha & \cos\alpha
\end{array}\!\!\right)\;\! S^{{}-1}.
\hfill
$ 
\\[1.55ex]
\indent
Owing to Lemma 12.1 and Theorem 12.2 for the matrix $Q_1^{}$ one of following two representations takes place: 
\\[1.25ex]
\mbox{}\hfill
$
Q_1^{}=T\;\! 
\left(\!\!
\begin{array}{rr}
\cos\alpha & {}\mp \sin\alpha
\\[0.25ex]
{}\pm \sin\alpha & \cos\alpha
\end{array}\!\!\right)\;\! T^{{}-1}.
\hfill
$ 
\\[1.5ex]
\indent
It is similar, as well as at the proof of the theorem 12.2, from identity (12.1) at $r=1$ we pass to identity 
\\[1.25ex]
\mbox{}\hfill                            % (12.8)
$
\varphi\bigl((t-2\alpha)({\rm mod}\, 2\pi)\bigr)= 
(\varphi(t)\mp 2\alpha)({\rm mod}\, 2\pi)
$ 
\ for all 
$
t\in ({}-\pi,\pi]. 
$
\hfill (12.8)
\\[1.25ex]
\indent
On the basis of (12.8) it is had following relations 
\\[1.5ex]
\mbox{}\hfill                            % (12.9)
$
\varphi\bigl((2l\alpha)({\rm mod}\, 2\pi)\bigr)=
(\varphi(0)\pm 2l\alpha)({\rm mod}\, 2\pi)
$ 
\ for all 
$
l\in\Z. 
$
\hfill (12.9)
\\[1.5ex]
\indent
Owing to a consequence from 
\vspace{0.35ex}
Kronecker theorem [31, p. 314 -- 315] and conditions 2 of these theorems it is had, 
\vspace{0.5ex}
that for everyone $t\in ({}-\pi, \pi] $ there is such sequence $\{l_s^{}(t)\} $ integers, that 
\vspace{0.25ex}
$\lim\limits_{s\to {}+\infty}(2\;\!l_s^{}(t)\;\!\alpha)({\rm mod}\, 2\pi)=t.$ 
From here, using relations (12.9), we come to conclusion, that homeo\-morphism, 
satisfying to identity (12.8), knows 
\\[1.5ex]
\mbox{}\hfill
$
\varphi(t)=(\varphi(0)\pm t)({\rm mod}\, 2\pi)$  for all $t\in ({}-\pi,\pi].
\hfill
$ 
\\[1.5ex]
\indent
Applying now in other replacement procedure, inverse already used, 
from last conjugating homeo\-morphism we come to required nondegenerate linear-fractional transformation. \k
\vspace{0.75ex}

{\bf Lemma 12.3.} 
{\it 
Let conjugating 
\vspace{0.35ex}
of nondegenerate linear-fractional transformations {\rm (11.3)} and {\rm (11.4)} 
homeomorphism $f\colon\overline{\R} \to\overline{\R}$ is such that{\rm:}
\vspace{0.35ex}

{\rm 1)} $f(0)=0,\ f(\infty)=\infty;$ 
\hfill {\rm (12.10)}
\\[1ex]
\indent
{\rm 2)} $f(\lambda x)=\lambda|\lambda|^{\alpha}f(x)$ for all $x\in\R,\ \alpha>{}-1;$ 
\hfill {\rm (12.11)}
\\[1.35ex]
\indent
{\rm 3)} $f\Bigl(\dfrac{ax+b}{cx+d}\Bigr)=\dfrac{Af(x)+B}{Cf(x)+B}$ \ for all $x\in\overline{\R},\ |b|+|c|>0;$ 
\hfill {\rm (12.12)}
\\[1.5ex]
\indent
{\rm 4)} the matrix 
$
P=\left(\!\!
\begin{array}{cc}
a & b
\\
c & d
\end{array}\!\!\right)=
S\;\! \left(\!\!
\begin{array}{cc}
\mu & 0
\\
0 & 1
\end{array}\!\!
\right)\;\! S^{-1}$  is such that 
$
|\mu|\ne 1,\ 
S=\left(\!\!
\begin{array}{cc}
a_{\ast}^{} & b_{\ast}^{}
\\
c_{\ast}^{} & d_{\ast}^{}
\end{array}\!\!
\right);
$
\\[1ex]
\indent
{\rm 5)} the subgroup of the group  $\R_+^{\ast}$ of positive real numbers on the multiplication, 
formed by numbers $|\lambda|$ and $\Bigl|\dfrac{b_{\ast}^{}}{d_{\ast}^{}}\Bigr|,$ 
is dense in the set $\R_+^{} $ of positive real numbers.
\\[0.5ex]
Then this homeomorphism looks like
\\[1.25ex]
\mbox{}\hfill
$
f(x)=x|x|^{\alpha}
$ 
\ for all 
$
x\in\overline{\R}. 
$
\hfill {\rm (12.13)}
\\[1.75ex]
}
\indent
{\sl Proof.} 
Let conditions (12.10) --- (12.12) are satisfied.
Owing to the condition 4 of lemma, relations (12.10), (12.11), Lemma 12.1 and Theorem 9.1 
\vspace{0.5ex}
we come to conclusion, that the matrix 
$
Q=\left(\!\!
\begin{array}{cc}
A & B
\\
C & D
\end{array}\!\!
\right)=
T\;\! \left(\!\!
\begin{array}{cc}
\theta & 0
\\
0 & 1
\end{array}\!\!
\right)\;\! T^{{}-1}
$ 
is such that 
$
\ln |\mu|\ln |\theta|>0,
\ 
T=\left(\!\!
\begin{array}{cc}
A_{\ast}^{} & B_{\ast}^{}
\\
C_{\ast}^{} & D_{\ast}^{}
\end{array}\!\!\right).
$ 
\vspace{0.75ex}

On the basis of (12.11) and (12.12) we receive relations 
\\[1.75ex]
\mbox{}\hfill
$
f\bigl(\lambda^k(P^m)^l(x)\bigr)=\lambda^k|\lambda|^{k\alpha}(Q^m)^l(f(x))
$
\ for all $x\in\R,$ 
\ for all $k, l, m\in\Z. 
\hfill
$ 
\\[1.75ex]
\indent
Passing in them to a limit at $m\to {}-\infty$ if $|\mu |> 1,$ and to a limit at $m\to {}+\infty $ if $|\mu | <1,$ 
we have 
\\[1.5ex]
\mbox{}\hfill                                               % (12.14)
$
f\biggl(\lambda^k\Bigl(\dfrac{b_{\ast}^{}}{d_{\ast}^{}}\Bigr)^l\;\!\biggr)=
\lambda^k|\lambda|^{k\alpha}
\biggl(\dfrac{B_{\ast}^{}}{D_{\ast}^{}}\biggr)^l
$ 
\, \ for all 
$k\in\Z,$
\ for all 
$
l\in\Z.
$
\hfill (12.14)
\\[2ex]
\indent
From the condition 5 of Lemma 12.3 follows, 
\vspace{0.35ex}
that for any positive real number $x> 0$ exist such sequences 
$\{k_s^{}(x)\}$ and $\{l_s^{}(x)\}$ of even integers, that
\\[1.75ex]
\mbox{}\hfill
$
\lim\limits_{s\to {}+\infty} |k_s^{}(x)|=
\lim\limits_{s\to {}+\infty} |l_s^{}(x)|={}+\infty,
\qquad 
\lim\limits_{s\to {}+\infty} \lambda^{k_s^{}(x)}
\Bigl(\;\!\dfrac{b_{\ast}^{}}{d_{\ast}^{}}\Bigr)^{l_s^{}(x)}=x.
\hfill
$ 
\\[1.75ex]
\indent
From here on the basis of relations (12.14) it is received, that 
\\[1.75ex]
\mbox{}\hfill
$
f(x)=x|x|^{\alpha}\lim\limits_{s\to {}+\infty}\,
\Bigl|\dfrac{b_{\ast}^{}}{d_{\ast}^{}}\Bigr|^{(\beta-\alpha)l_s^{}(x)}
$ 
\ \, for all $x>0,
\hfill
$ 
\\[2ex]
where 
$\Bigl|\dfrac{B_{\ast}^{}}{D_{\ast}^{}}\Bigr|=
\Bigl|\dfrac{b_{\ast}^{}}{d_{\ast}^{}}\Bigr|^{1+\beta}.
$ 
As from the condition 5 this lemma follow, that 
\\[1.75ex]
\mbox{}\hfill
$
\lim\limits_{s\to {}+\infty}\,
\Bigl|\dfrac{b_{\ast}^{}}{d_{\ast}^{}}\Bigr|^{(\beta-\alpha)l_s^{}(x)}=1.
\hfill
$ 
\\[1.75ex]
\indent
Therefore,
\\[1.25ex]
\mbox{}\hfill                                          % (12.15)
$
f(x)=x|x|^{\alpha}
$ 
\ for all 
$
x>0. 
$
\hfill (12.15)
\\[1.75ex]
\indent
Having executed replacement $\xi(x)={}-f({}-x)$ for all $x\in\overline{\R},$ 
for homeomorphism $\xi$ by the reasonings similar resulted above, we come to relations 
$
\xi(x)=x|x|^{\alpha}$ for all $x>0,$ 
and from them with the help of return replacement we come to relations 
\\[1.5ex]
\mbox{}\hfill                                          % (12.16)
$
f(x)=x|x|^{\alpha}
$
\ for all 
$
x<0. 
$
\hfill (12.16)
\\[1.5ex]
\indent
Uniting formulas (12.10), (12.15) and (12.16), we come to representation (12.13). \k
\vspace{0.5ex}

Now by direct calculations on the basis of identity (12.13) we do a conclusion, 
that in case of general situation at performance of conditions of the lemma 12.3 conjugating homeomorphism looks like 
$
f(x)=x
$ 
for all $x\in\overline{\R}.$ 

Now on the basis of the reasonings similar resulted at proof of the theorem 12.1, 
taking into account the theorem 12.6, we come to the statement of the theorem 12.5.

At first we will consider smooth and holomorphic conjunction of Abelian linear-fractional groups $PL^3$ and $PL^4.$
\vspace{0.5ex}

{\bf Theorem 12.7.} 
\vspace{0.25ex}
{\it 
Let conditions of the theorem {\rm 12.1}  are satisfied.
Then for smooth {\rm(}holomorphic{\rm)} conjunction of Abelian linear-fractional groups 
\vspace{0.5ex}
$PL^3$ and $PL^4$ it is necessary   and enough, that \,
$\dfrac{q_{1r}^{}}{q_{2r}^{}}=\Bigl(\dfrac{p_{1r}^{}}{p_{2r}^{}}\Bigr)^{\varepsilon}
$
for all $r\in I,\ \varepsilon^2=1.
$
}
\vspace{0.75ex}

{\sl Proof} of the given statement is similar to the proof of the theorem 12.1 and is based on the theorem 9.3. \k
\vspace{0.5ex}

In case of non-Abelian linear-fractional groups $PL^3$ and $PL^4$ the following statement takes place.
\vspace{0.5ex}

{\bf Theorem 12.8.} 
{\it 
From a smooth conjunction of non-Abelian linear-fractional groups 
$PL^3$ and $PL^4$ follows them conjunction, which is carried out by nondegenerate  linear-fractional transformation.
}
\vspace{0.25ex}

{\sl Proof} of this theorem is spent similarly to the proof of the theorem 12.5 
on the basis of following three auxiliary statements. \k
\vspace{0.75ex}

{\bf Lemma 12.4.} 
{\it 
Let diffeomorphism $f\colon\overline{\R}\to\overline{\R}$  is such that{\rm:}
\vspace{0.35ex}

{\rm 1)} relations {\rm (12.10)} are carried out{\rm;}
\vspace{0.5ex}

{\rm 2)} $f(\lambda x)=\lambda f(x)$ for all $x\in\R,\ \lambda\ne 0,\ \lambda\ne 1.$ 
\hfill {\rm (12.17)}
\\[0.75ex]
\indent
Then this diffeomorphism looks like
\\[1.5ex]
\mbox{}\hfill                                     % (12.18)
$
f(x)=ax
$ 
\ for all 
$
x\in\overline{\R}.
$ 
\hfill {\rm(12.18)}
\\[1.5ex]
}
\indent
{\sl Proof.} From (12.17) we come to conclusion about justice of relations 
\\[1.25ex]
\mbox{}\hfill
$
f(\lambda^l x)=\lambda^l f(x)$ 
\ for all $x\in\R,$ 
\ for all $l\in\Z,
\hfill
$ 
\\[1.5ex]
differentiating which on $x,$ we come to identities 
\\[1.25ex]
\mbox{}\hfill
$
f^{\;\!\prime}(\lambda^l x)=f^{\;\!\prime}(x)$ for all $x\in\R,$ for all $l\in\Z.
\hfill
$ 
\\[1.25ex]
\indent
Owing to that map $f $ is a diffeomorphism, we receive identity 
\\[1.25ex]
\mbox{}\hfill
$
f^{\;\!\prime}(x)=a$ for all $x\in\R.
\hfill
$ 
\\[1.25ex]
From here taking into account relations (12.17) we come to representation (12.18). \k
\vspace{0.5ex}

Similarly we prove the following statement.
\vspace{0.5ex}

{\bf Lemma 12.5.} 
{\it 
Let diffeomorphism $f\colon\overline{\R}\to\overline{\R}$  is such that{\rm:}

{\rm1)} $f (\infty) = \infty; $
\vspace{0.35ex}

{\rm2)} $f(x+p)=f(x)+q$ for all $x\in\R,\ pq\ne 0.$
\\[0.75ex]
Then this diffeomorphism looks like 
$
f(x)=qp^{{}-1}x+a$ for all $x\in\overline{\R}.$
}
\vspace{1ex}

{\bf Lemma 12.6.} 
{\it 
Let diffeomorphism $f\colon \overline{\R}\to\overline{\R}$ is such that 
\\[2ex]
\mbox{}\hfill
$
f\biggl(\dfrac{x\cos\alpha-\sin\alpha}{x\sin\alpha+\cos\alpha}\biggr)=
\dfrac{f(x)\cos\beta\mp \sin\beta}{\pm f(x) \sin\beta+\cos\beta}
$
\ \ for all 
$
x\in\overline{\R},
\quad 
\sin\alpha\ne 0.
\hfill
$
\\[2ex]
Then this diffeomorphism represents nondegenerate  linear-fractional transformation.
}
\vspace{0.35ex}

{\sl Proof.} 
It is similar, as well as at the proof of theorems 12.2 and 12.6, from identity from a condition of the lemma 12.6 
we pass to identities (12.8); and from it we pass to identities 
\\[1.5ex]
\mbox{}\hfill
$
\varphi\bigl((t-2\;\!l\alpha)({\rm mod}\, 2\pi)\bigr)= (\varphi(t)\mp 2\;\!l \alpha)({\rm mod}\, 2\pi)
$
\ for all $t\in ({}-\pi,\pi],$ \ for all $l\in\Z.
\hfill
$ 
\\[1.5ex]
\indent
Differentiating last identities on $t $, taking into account that map $\varphi $ is a diffeomorphism, 
in a similar way, as well as earlier, we receive representations 
\\[1.5ex]
\mbox{}\hfill
$
\varphi(t)=({}\pm t+a)({\rm mod}\, 2\pi)
$ 
\ for all $t\in ({}-\pi,\pi].
\hfill
$ 
\\[1.5ex]
\indent
Now it is similar to the proof of the theorem 12.6 we come to conclusion, that conjugating homeomorphism $f $ 
is nondegenerate linear-fractional transformation. \k
\vspace{0.5ex}

On the basis of the theorem 12.8 and a course of proofs of theorems 12.2 --- 12.4 it is had such statement.
\vspace{0.5ex}

{\bf Theorem 12.9.} 
{\it 
From a smooth conjunction of linear-fractional groups $PL^3$ and $PL^4$ follows 
them holomorphic conjunction which is carried out  by nondegenerate linear-fractional transformation.
}
\\[3.75ex]
\centerline{
\large\bf  
13. Applications to real nonautonomous Riccati equations 
}
\\[1.5ex]
\indent
Theorems 12.1 --- 12.9 and Lemmas 12.1 --- 12.6 allow on the basis of Theorems 1.1 --- 1.3 to spend topological, smooth and holomorphic classifications of real nonautonomous Riccati equations of a kind (11.1).
Besides, from the theorem 12.5 it is had such statement.
\vspace{0.5ex}

{\bf Theorem 13.1.} 
{\it 
From topological equivalence of real nonautonomous Riccati equations with non-Abelian phase groups 
of general situation follows them holomorphic equivalence.
}
\vspace{0.5ex}

From the given theorem, in particular, follows, that real nonautonomous Riccati equations with coefficients, 
holomorphic on path connected holomorphic varieties with non-Abelian fundamental groups, are structurally unstable.

}
\end{document}